\documentclass{amsart}
\usepackage{amssymb,latexsym,graphics,enumerate}
\usepackage[normalem]{ulem}
\usepackage[mathscr]{eucal}
\usepackage{fullpage}
\usepackage{color}
\usepackage{amsthm, amsopn, amsfonts}
\usepackage{multirow}
\usepackage{hyperref}
\hypersetup{linktocpage}

\newtheorem{thm}{Theorem}[section]
\newtheorem{cor}[thm]{Corollary}
\newtheorem{lem}[thm]{Lemma}
\newtheorem{prop}[thm]{Proposition}

\theoremstyle{remark}
\newtheorem{rem}[thm]{Remark}

\newtheorem{defn}[thm]{Definition}
\theoremstyle{definition}

\newtheorem{assum}[thm]{Assumption}
\theoremstyle{definition}

\newcommand{\la}{\langle}
\newcommand{\ra}{\rangle}

\newcommand{\lp}[2]{\Vert \, #1 \, \Vert_{#2}}

\newcommand{\td}{\widetilde}

\newcommand{\snabla}{ {\, {\slash\!\!\! \nabla  }} }

\newcommand{\jap}[1]{ {\langle #1 \rangle} }
\newcommand{\bL}{  {\underline{L}} }
\newcommand{\CE}{ {C\!E} }
\newcommand{\CH}{ {C\!H} }

\newcommand{\NLE}{ {N\!L\!E} }
\newcommand{\LE}{ {L\!E} }

\newcommand{\dint}{ {\int\!\!\!\!\int} }

\newcommand{\case}[2]{\noindent\textbf{Case #1:}(\emph{#2})}
\newcommand{\step}[2]{\noindent\textbf{Step #1:}(\emph{#2})}
\newcommand{\Part}[2]{\noindent\textbf{Part #1:}(\emph{#2})}

\begin{document}

\title[Vector Fields]
{A Vector Field Method for Non-Trapping, Radiating Spacetimes}
\author{Jes\'us Oliver}
\address{Department of Mathematics,
California State University, East Bay (CSUEB)
25800 Carlos Bee Blvd,
Hayward, CA 94542 USA}
\email{jesus.oliver@csueastbay.edu}
\thanks{The author was supported in part by NSF grant DMS-1001675
through Jacob Sterbenz.}
\date{}
\dedicatory{}
\commby{}

\maketitle
\begin{abstract}
We study the global decay properties of solutions to the linear wave
equation in 1+3 dimensions on time-dependent, weakly asymptotically
flat spacetimes. Assuming non-trapping of null geodesics and a local
energy decay estimate, we prove that sufficiently regular solutions to
this equation have bounded conformal energy. As an application we also show 
a conformal energy estimate with vector fields applied to the solution as well as a global $L^{\infty}$ decay
bound in terms of a weighted norm on initial data. For solutions to the wave equation
in these dynamical backgrounds, our results reduce the problem of establishing the classical
pointwise decay rate $t^{-3/2}$ in the interior and $t^{-1}$ along outgoing null cones
to simply proving that local energy decay holds.
\end{abstract}
\tableofcontents



\section{Introduction}

Let $(\mathcal{M},g)$ be a 4-dimensional, smooth, asymptotically 
flat, Lorentzian manifold. Assume $\mathcal{M}$ is
of the form $\mathbb{R}\times \mathbb{R}^3$ and that there exist
global coordinates $(t,x^i)$ such that the level sets of $t$ are uniformly
space-like. In this work we study the
dispersive properties of solutions to the wave equation:
 \begin{align}
 	 \Box_g\phi \ &= \ F(t,x) \ , \qquad 
	 (\phi , \partial_t\phi)(0,x) \ = \ (\phi_0, \phi_1) \ ,
	 \label{wave_eqn}
 \end{align}
where $(\phi_0,\phi_1)\in C_{c}^{\infty}(\mathbb{R}^{3})$
and,
\begin{equation*}
	\Box_g =  |\det g|^{-\frac{1}{2}}\partial_{\alpha} g^{\alpha\beta}
	|\det g|^{\frac{1}{2}}\partial_{\beta} \ ,
\end{equation*}	
 in local coordinates.  One can think of problem
\eqref{wave_eqn} with $F=\mathscr{N}(\phi,\nabla\phi)$
and $\mathscr{N}$ a non-linear function as a toy model for many important systems
of hyperbolic PDE including: Maxwell-Klein-Gordon, Yang-Mills, and Wave Maps.
A necessary first step in understanding the stability properties of these systems is
to prove that smooth solutions launched from sufficiently small, well-localized initial data
can be extended for all time. One particularly fruitful strategy for proving this type
of result is to show that solutions to the linear wave equation propagating in a fixed
background, evolving from $C_c^{\infty}$ data, have pointwise rates similar to
those in flat space. Provided that the method of proof is robust enough, one can
often leverage the linear bound into a small-data global existence result for some
perturbative non-linear problems. In this work we focus on one such method for
proving $L^{\infty}$ linear decay: the vector field method first introduced by S.
Klainerman in \cite{MR784477}. In the case of $1+3$ dimensions, this $L^2$-based
method has a proven record of success in dealing with many small-data semilinear
and quasilinear problems -- at least when the metric $g$ is uniformly close to Minkowski
and the nonlinearities have special structure. The most well-known demonstration
the power of this method is the proof of the global non-linear stability of the Minkowski
space by Christodoulou-Klainerman \cite{MR1316662}. Since the appearance of that
work there has been a concerted effort to extend the vector field method to backgrounds
far from Minkowski space. In many cases these efforts have been successful and have
yielded small data global existence results for non-linear problems in spacetimes
such as: Minkowski space with obstacles \cite{MR2217314}, exterior Kerr with
$|a|\ll M$ (\cite{Lindblad:2013ly}, \cite{Luk:2010fk}) and time-dependent,
inhomogeneous media (\cite{MR3056620}, \cite{Yang:2013kq})

In light of this connection with nonlinear stability results, the problem of pointwise decay
via vector field method for solutions to the linear wave equation on large, asymptotically
flat backgrounds has been intensely studied. In the last few years the main focus of research
activity has been on the Schwarzschild and Kerr metrics. This is due to the fact that it is widely
believed that this type of strategy has the best chance of proving
the non-linear stability of the (exterior) subextremal Kerr family. For cutting-edge results for
the linear problem in black hole spacetimes we mention, without being exhaustive,
the work of Dafermos-Rodnianski-Shlapentokh\_Rothman \cite{Dafermos:2014kq}
and Metcalfe-Tataru-Tohaneanu \cite{MR2921169} (see also \cite{Andersson:2013kx} for the
case of Maxwell field). In our work we will turn away
from black hole spacetimes to focus on a different, but related, problem
for which vector field methods have not yet been developed: pointwise
linear decay for solutions to Eq. \eqref{wave_eqn} on radiating, non-trapping spacetimes.
The motivation to look at this problem comes from the fact that these spacetimes provide a model for the far
exterior portion of a dynamic perturbation of a body emitting gravitational waves.
Heuristically speaking, one can think of gravitational waves as local disturbances
to the spacetime geometry propagating along (characteristic) null hypersurfaces.
Their presence perturbs the metric $g$ and leads it to decay towards flat space in
the null outgoing region at slower rates than the Schwarzschild or Kerr spacetimes.
This, in turn, precludes us from just repeating the proofs of the vector field
methods that are available near Minkowski space and forces us to produce a
new method for dealing with these weak asymptotics. In short, weak
decay for the metric takes away the classical VF proofs of $L^{\infty}$ decay
for solutions to Eq. \eqref{wave_eqn}.

In order to quantitatively define our weak decay for the metric
we take the work of Klainerman-Nicolo \cite{MR1946854} as a starting point (see also \cite{MR0147276} and \cite{MR0149908}).
In that work, the authors use a double null foliation
to derive a hierarchy of decay for all connection coefficients and curvature components
near null infinity for radiating spacetimes satisfying the Einstein vacuum equations.
In our work, instead of relying on a null frame, we {\it{assume the existence a coordinate system}} such
that metric coefficients, after subtracting the Minkowski metric, obey symbol bounds
in the null outgoing region which include those of \cite{MR1946854} as a special case
(see assumption \ref{ex_norm_coords} below). We chose a coordinate-dependent
condition mostly to simplify matters since the extensive geometric computations associated to null frames
are unnecessary for this problem. We also mention that our decay conditions in the null outgoing region
are general enough to include, as examples, the
spacetimes of Lindblad-Rodniaski \cite{MR2680391} and Bieri \cite{MR2531716}.
%
In the timelike region, i.e. within the domain of dependence of a compact set,
our metric is allowed to be time-dependent and may remain far from Minkowski space
for all time as long as a local energy decay estimate and non-trapping of null geodesics hold.
For spacetimes satisfying all these assumptions, we are able to
produce a novel vector field method yielding weighted $L^2$ and pointwise decay bounds
that are analogous to what is available near Minkowski space via classical vector field method.
Furthermore, the norms we impose on the initial data are suitable for
non-linear applications -- a topic which we will explore in subsequent work. We now give
a detailed description of the spacetimes we work with.


\subsection{Decay assumptions on the metric}

There are mainly two regions that need to be considered separately: a sufficiently
large compact set and its exterior. In general what happens outside of a compact
set only needs a detailed description where $r\sim t$. In view of this, we make the
following:

\begin{defn}
Let $r  =  \sum_{i=1}^{3}\sqrt{(x^{i})^{2}}$. Define the {\bf wave zone} to
be the set $\{(t,x)  |   \frac{2}{3}t < r < \frac{3}{2} t\}$ and the
{\bf interior region} to be the wedge $\{(t,x)  |  r< \frac{1}{2}t\}$.
\end{defn}

We let $i,j=1,2,3$ and $\alpha,\beta=0,1,2,3$ and make use
of the Einstein summation convention throughout.
The first condition on the metric is the following:

\begin{assum}[Existence of normalized coordinates]\label{ex_norm_coords}
There exists $u(t,x^i)\in C^{\infty}(\mathcal{M})$ with 
$du\neq 0$ satisfying the following conditions:

\begin{enumerate}[i)]
	\item $u=t-r$ on the set $\{t < 1\}\cup\{r<\frac{1}{2} t\}\cup\{r>2t\}$.
	\item (Asymptotics). There exists $\delta>0$ such that on the set $\{t > 1\}\cap\{\frac{1}{2} t<r< 2t\}$
	the function $u\approx t-r$ in the following sense:
		\begin{equation}
				|\partial_{t,x}^J (\partial_t u - 1, \ \partial_i u +\omega^i) | \ \lesssim \
				r^{-\delta-|J|}\tau_0^{-|J|} \ , \label{grad_u_bnd}
		\end{equation}
		where:
		\begin{equation*}
			\tau_-=\la u\ra \ , \qquad \tau_+=C+ u+2r \ , \qquad \tau_0=\tau_-\tau_+^{-1} \ , 
		\end{equation*}
		and $C$ is chosen large enough so that
		$\tau_+>1$ in $t\geqslant 0$. In particular $\tau_+ \approx \la t+r\ra$.
	\item(Renormalization). There exists $\gamma >0$ such that for any multi-index $J$
	the inverse metric coefficients $g^{\alpha\beta}$ in $(u,x^i)$ coordinates satisfy:
		\begin{subequations}\label{mod_coords}
		\begin{align}
			|\partial_u^k \td \partial_x^J( g^{ij}   -  \delta^{ij})  | \ &\lesssim \
			 \jap{r}^{-k-|J|-\delta}\tau_{0}^{-k} \tau^{-\gamma k}  \ ,  \label{mod_coordds1} \\			
			  |\partial_u^k \td \partial_x^J( d^{\frac{1}{2}} g^{ui}
			 + \omega^i)  | \ &  \lesssim \
			   \jap{r}^{-k-|J|-\delta}\tau_{0}^{\frac{1}{2}-k} \tau^{-\gamma k}
			\ , \label{mod_coordds2}  \\
			|\partial_u^k \td \partial_x^J ( g^{ui}   - \omega^i\omega_j g^{uj}  )  |
			 \ &\lesssim \
			   \jap{r}^{-k-|J|-\delta}\tau_{0}^{1-k} \tau^{-\gamma k}
			 \ , \label{mod_coordds3} \\
			 |\partial_u^k \td \partial_x^J g^{uu}  | \ &\lesssim \
			 \jap{r}^{-k-|J|-\delta}\tau_{0}^{2-k}\tau^{-\gamma k}
			 \ , \label{mod_coordds4}
		\end{align}
		where $d=|\det(g_{\alpha\beta})|$ is in $(u,x^{i})$ coordinates and:
		\begin{equation*}
			\td \partial_{x}=\partial_{x}|_{u=const} \ , \qquad 
			\omega^{i}=\frac{x^i}{r} , \qquad \tau=\frac{\tau_{+}}{\jap{r}} \ .
		\end{equation*}
		\end{subequations}
\end{enumerate}
\end{assum}

In the sequel we also refer to the $(u,x^i)$
coordinates as \emph{Bondi Coordinates}.

\begin{rem}[Decay outside the wave zone]
Let $g^{\alpha\beta}$ be in $(t,x^i)$ coordinates.
In the interior region a computation using
the chain rule yields:
	\begin{equation}
			|\partial_t^k  \partial_{x}^J (g^{\alpha\beta}-m^{\alpha\beta})| \ \lesssim \ 
			\langle r \rangle^{-k-|J|-\delta }\tau^{-\gamma k}
			\  , \label{local_nontrapping}
	\end{equation}
with $m^{-1} = m  =  diag(-1,1,1,1)$. By a similar reasoning, inside
$\{   \frac{1}{2} t < r< \frac{2}{3}t \} \cup \{\frac{3}{2}t < r\}$ we have:
\begin{equation}
	| \partial_{t}^{k} \partial_{x}^{J} (g^{\alpha\beta}-m^{\alpha\beta})| \ \lesssim \ 
	\jap{t+r}^{-k-|J|-\delta} \ . \label{spacel_timel_infty}
\end{equation}
In practice we use \eqref{mod_coords} within the wave zone, we switch
to \eqref{local_nontrapping} for the interior, and rely on \eqref{spacel_timel_infty} outside
these two regions.

\end{rem}

\begin{rem}
This hierarchy of decay for different components
of the metric along outgoing null directions is consistent with the radiation of gravitational waves
and with the peeling estimates in \cite{MR1946854}. These decay rates are also consistent
with the metrics constructed in the stability of Minkowski space in wave coordinates 
\cite{MR2680391}. Note, in particular, that the `shear' terms $g^{ui}-\omega^i\omega_j g^{uj} $ decay
only slightly faster than a solution to the wave equation on any fixed hypersurface $u=const$.
\end{rem}

\begin{rem}(Examples).
In Schwarzschild one can define $r^{*}=r+2M\log|\frac{r}{2M}-1|$ and let $u=t-r^{*}$.
The metric in $(u,x^{i})$ coordinates then satisfies assumption \ref{ex_norm_coords}
in the wave zone with $\delta=1$. The subextremal Kerr family also satisfies this assumption
in the wave zone (see the upcoming work \cite{JOJS} for details).
\end{rem}

\begin{rem}\label{u_rem}
Estimates \eqref{grad_u_bnd} and \eqref{mod_coords} also hold
if we trade $u$ for
\begin{equation*}
	\td{u}=\chi u + (1-\chi)(t-r) \ ,
\end{equation*}
where $\chi\equiv 1$ when
$\tau_0\leqslant c$ for any $0<c<1$ with bounds $|(\tau_{+}\partial)^{k} \chi |\lesssim 1$.
Thus the assumption that $u\equiv t-r$ in the set
$\{t < 1\}\cup\{r<\frac{1}{2} t\}\cup\{r>2t\}$ is not the weakest condition we can impose.
We chose this condition mostly to simplify matters dealing with estimates in the interior region. However, the reader should keep
in mind that one only needs to use the exact form of $u$ in a narrow wedge
$|t-r|\ll r$.
\end{rem}

\begin{rem}
By estimates \eqref{mod_coords} the $u=const$ hypersurfaces are approximately null.
We also note that, by construction, the $\td \partial_{x}$ derivatives are tangential to
these hypersurfaces.
%
\end{rem}

\begin{rem}Let $R_0>0$. Inside sets of the form $ r \leqslant R_0$,
our assumptions allow for $g$ to be a large perturbation of the Minkowski
metric for all $t>0$. In particular, the metric $g$ does not have to converge to a stationary
metric as $t\rightarrow +\infty$ in this region.
\end{rem}


\subsection{Non-trapping and local energy decay}
Next we introduce our non-trapping and local energy decay assumptions.
These are the two key ingredients that allow us to handle large deformation errors inside
regions of the form $r\leqslant R_0$.
Let us start by making the non-trapping assumption precise:

\begin{assum}[Quantitative non-trapping for null geodesics] \label{non-trapping} Let $R_0>0$ and
$\gamma(s)$ be a forward, affinely parametrized null geodesic satisfying:
\begin{equation*}
	\gamma(0)\in\{ r \leqslant R_0 \}, \qquad  \dot{\gamma}s\equiv 1 , \qquad
	\dot{\gamma}t\big|_{s=0}=1 \ .
\end{equation*}
Then, for any such $\gamma$, there exists
a uniform constant $C=C(R_0,g)$ such that $\gamma(s)\in \{ r > R_0\} $ for all $s\geq C$.

\end{assum}

For any norm $B[t_0,t_1]$ we set up the notation:
\begin{align*}
	\ell^p_r B[t_0,t_1] \ &= \ \lp{f}{\ell^p_r B[t_0,t_1]}^p=
	\sum_{i\geqslant 0} \lp{\chi_i(r) f}{B[t_0,t_1]}^p \ , \\
	\ell^p_t B[t_0,t_1] \ &= \ \lp{f}{\ell^p_t B[t_0,t_1]}^p=
	\sum_{k\geqslant 0} \lp{\chi_k(t) f}{B[t_0,t_1]}^p \ , \\
	\ell^p_u B[t_0,t_1] \ &= \ \lp{f}{\ell^p_u B[t_0,t_1]}^p=
	\sum_{j\geqslant 0} \lp{\chi_j(u) f}{B[t_0,t_1]}^p \ ,
\end{align*}
where $\chi_i(r)$ is a series of dyadic cutoffs on
$\jap{r}\approx 2^i$, $\chi_k(t)$ is a series of dyadic cutoffs on
$t\approx 2^k$ covering $[t_0,t_1]$, and $\chi_j(u)$ is a series of dyadic cutoffs on
$\jap{u}\approx 2^j$. We also make the obvious modification
for $p=\infty$. Using this we can define the \emph{Local Energy Decay} (LED) norms:
\begin{align*}
		\lp{\phi}{\LE[t_0,t_1]}  \ &=  \ \lp{\jap{r}^{-\frac{1}{2}}\nabla\phi}{\ell_{r}^{\infty}L^{2}[t_0,t_1]}
		\ + \ \lp{\jap{r}^{-\frac{3}{2}}\phi}{\ell_{r}^{\infty}L^{2}[t_0,t_1]}   \ , \\
		\lp{F  }{\LE^{*}[t_0,t_1]} \ &= \   \lp{\jap{r}^{\frac{1}{2}}F}{\ell_{r}^{1}L^{2}[t_0,t_1]}  \ .
\end{align*}
These norms allow us to state the final decay assumption on the metric:

\begin{assum}[Local energy decay estimate]\label{ls_est} For all values
$0 \leqslant t_0 < t_1$ the evolution \eqref{wave_eqn} satisfies:
\begin{equation}
	\lp{\phi}{\LE[t_0,t_1]} \ \lesssim \  \lp{\nabla\phi(t_0)}{L_{x}^{2}} \ + \  
	\lp{\Box_{g}\phi}{ {\LE}^{*}[t_0,t_1]}	 \ . \label{LS1}
\end{equation}
\end{assum}

Non-trapping of null geodesics is a necessary condition for \eqref{LS1}
to hold in this form. If there's trapping, the work of Ralston \cite{MR0254433}
and some additional geometric optics considerations can be used to show
that the LED estimate, if it holds at all, must lose derivatives in a neighborhood
of the trapping set. In the non-trapping case, estimates such as \eqref{LS1}
date back to work of Morawetz \cite{MR0151712}, \cite{MR0234136} and are known to hold in
a variety of settings. In the case of Minkowski space, Keel-Smith-Sogge
proved a limiting version of this estimate \cite{MR1945285} (see also \cite{MR2128434}).
For uniformly small, time-dependent perturbations of Minkowski, Alinhac \cite{MR2666888}
and Metcalfe-Tataru \cite{MR2944027} both established this result. The work of
Bony-Hafner \cite{MR2748617} extended the validity of this estimate to the case
of large, stationary, non-trapping metrics of the form $ds^2=-dt^2+h_{ij}(x)dx^idx^j$
with $h$ Riemannian (see also \cite{Rodnianski:2011kq}, \cite{MR2762995}). 
As the Schwarzschild and Kerr solutions have trapped null geodesics,
this work will not apply to full perturbations of such spacetimes. However, a suitable modification in
the upcoming work \cite{JOJS} will do so.

Since the LED estimate is generally expected to hold for a large class of spacetimes,
it is a natural assumption to include in our problem. In particular, the LED estimate 
should hold for time-dependent, non-trapping, asymptotically flat spacetimes satisfying
a smallness condition for $ \partial_{t}g$ and it should also hold (with loss of derivatives) for the
domain of outer communications of a small time-dependent perturbation of the sub-extremal Kerr family.
In the non-trapping, time-dependent case, some work is already under way to prove this result
(upcoming work of Sterbenz-Tataru). In the trapping case, the work \cite{MR2921169} already established the result
for fast decaying perturbations of Kerr with $|a|\ll M$. 
In view of this, we will take the estimate
as given and focus on developing a precise understanding of the asymptotic properties of the solution
via vector fields.


\subsection{Weighted norms}\label{normsect}

We define the \emph{Weighted LED} norms by:
\begin{align*}
		\lp{\phi}{\LE^{a,b}[t_{0},t_{1}]} \ &= \ \lp{\tau_{+}^{a}\tau_{0}^b\phi }{\LE[t_{0},t_{1}]}
		\ , \qquad
		\lp{\phi}{\LE^a[t_{0},t_{1}]} \ = \ \lp{\phi}{\LE^{a,0}[t_{0},t_{1}]} \ , 
\end{align*}
with analogous definitions for $\lp{F}{\LE^{*,a,b}} $. Next we set up the fixed-time and 
null energies that we will use in the sequel. Define the vector $\nabla\phi=(\partial_{0}\phi,...,\partial_{3}\phi)$ where 
$\partial_{0},...,\partial_{3}$ denotes any basis which can be written as a bounded linear combination of $(t,x)$ coordinate derivatives.
For $t\geq 0$ we define the fixed-time \emph{Conformal Energy}:
\begin{align*}
	\lp{\phi(t)}{\CE} \ = \ \lp{\tau_{+}\tau_{0}\nabla\phi (t)}{L_{x}
	^{2}} \ +\  \lp{\tau_{+}(\td \nabla_{x}\phi, r^{-1}\phi )  (t)}{L_{x}^{2}} \ . \notag
\end{align*}
For a smooth, positive weight function $\Omega$ we also have the \emph{Conjugated Conformal Energy}:
\begin{align}
	 \lp{\phi(t)}{{}^{\Omega}\CE}  \ =  \  \lp{\tau_{+}\tau_{0} \ \Omega^{-1} \nabla(\Omega\phi)  (t)}{L_{x}
	^{2}} \ + \  \lp{\tau_{+}\Omega^{-1}\td \nabla_{x}(\Omega\phi)(t)}{L_{x}^{2}} \label{conj_conf_en1} \ .
\end{align}
In this work we only use the conformal weights ${}^{I}\Omega=\jap{r}, {}^{I\!I}\Omega=\tau_{-}\tau_{+}$.
We set up the notation:
\begin{align*}
	\lp{\cdot}{{}^{I} C\!E} \ &= \ \lp{\cdot}{{}^{{}^{I}\Omega} \CE} \ , \qquad
	\lp{\cdot}{{}^{I\!I} C\!E} \  = \ \lp{\cdot}{{}^{{}^{I\!I}\Omega} \CE} \ .
\end{align*}
Note that $\lp{\cdot}{\CE}$ and $ \lp{\cdot}{{}^{\Omega}\CE}$ are fixed-time
norms and do not contain null energies. To introduce these we define the scale of spaces:
\begin{align*}
		\lp{\phi}{\NLE[t_{0},t_{1}]}  &=  \lp{\jap{r}^{-\frac{1}{2}}\phi}
		{\ell_{u}^{\infty}\ell_{r}^{\infty}L^2[t_{0},t_{1}]} \ , \qquad
		\lp{F}{\NLE^*[t_{0},t_{1}]}  =  \lp{\jap{r}^{\frac{1}{2}}F}
		{\ell_{u}^{1}\ell_{r}^{1}L^2[t_{0},t_{1}]}\ .
\end{align*}
Our (weighted) null energies are defined to be:
\begin{align*}
	\lp{\phi}{{}^{\Omega}C\!H [t_{0},t_{1}]} \ &= \ \lp{ \jap{r}^{\frac{1}{2}}  \Omega^{-1}\td
	\partial_r(\Omega \phi)}{\NLE^{\frac{1}{2},-\frac{1}{2}}[t_{0},t_{1}]} \ ,\\
	\lp{\cdot}{{}^{I} C\!H[t_0,t_1]} \ &= \ \lp{\cdot}{{}^{{}^{I}\Omega} \CH[t_0,t_1]} \ ,\\
	\lp{\cdot}{{}^{I\!I} C\!H[t_0,t_1]} \  &= \ \lp{\cdot}{{}^{{}^{I\!I}\Omega} \CH[t_0,t_1]} \ .
\end{align*}
Using the norms above as building blocks we then define the $S$ norm:
\begin{align*}
	\lp{\phi}{S[t_{0},t_{1}]} \ = \ \lp{\phi}{{}^{I} C\!H[t_0,t_1]}  + \lp{\phi}{{}^{I\! I}\CH[t_0,t_1]}
	+ \lp{\chi_{r<\frac{1}{2}t}\phi}{\ell_{t}^{\infty}\LE^{1}[t_{0},t_{1}]}  \ .
\end{align*}
Associated to this is the source term norm $N$:
\begin{align*}
	&\lp{F}{N[t_{0},t_{1}]}  =  \lp{(\chi_{r<\frac{1}{2}t}+\chi_{r>2t})F}{\LE^{*,1,\frac{1}{2}}[t_{0},t_{1}]}
	  +  \lp{ \chi_{\frac{t}{2}<r<2t}F}{\NLE^{*,1,\frac{1}{2}}[t_{0},t_{1}] } \ .
\end{align*}


\subsection{Vector fields and associated norms}

We now define versions of the previous norms which incorporate additional
decay. These will be stated in terms of modifications of the usual
Lorentz vector fields which are adapted to the null geometry of our spacetime as dictated by the
function $u(t,x)$. First we define the Lie algebra in $(t,x^i)$ coordinates (see also Eq.
\eqref{vf_bondi}):
		\begin{equation*}
			\mathbb{L} \  =  \ \big\{  T=\frac{1}{u_t}\partial_t  , \  
			S=\frac{u-ru_r}{u_t}\partial_t + r\partial_r  , \  \Omega_{ij}=
			\Omega_{ij, mink} - 
			\frac{\Omega_{ij, mink} u}{u_t}\partial_t \big\} \ , \label{mod_fields}
		\end{equation*}
where $\Omega_{ij,mink}=x^{i}\partial_{j}-x^{j}\partial_{i}$. With $\Gamma\in \mathbb{L}$ we define the 
higher order norms to be:
\begin{align*}
	\lp{\phi}{\LE_{1}[t_{0},t_{1}]}  &=   \sum_{|J|\leqslant 1}
	\lp{\Gamma^{J} \phi }{{\LE}[t_{0},t_{1}]}  \ ,  \\
	\lp{\phi(t)}{\CE_{1}[t_{0},t_{1}]}   &=   \sum_{|J|\leqslant 1}
	\lp{\Gamma^{J} \phi(t)}{\CE[t_{0},t_{1}]} \ , 
\end{align*}
with analogous definitions for the weighted $\LE$ and $\CH$ norms respectively. Using
these as building blocks we define the higher order $S$ norm by:
\begin{align*}
	&\lp{\phi}{S_{1}[t_{0},t_{1}]} \ = \ 
	 \lp{\phi}{{}^{I} C\!H_{1}[t_0,t_1]}  + \lp{\phi}{{}^{I\! I}\CH_{1}[t_0,t_1]}  + 
	\lp{\chi_{r<\frac{1}{2}t}\phi}{\ell_{t}^{\infty}\LE_{1}^{1}[t_{0},t_{1}]}  \ .
\end{align*}
Associated to this are the higher-order source term norms:
\begin{align*}
	\lp{F}{N_{1}[t_{0},t_{1}]} \  &=  \ \sum_{|J|\leqslant 1}
	\lp{\Gamma^{J} F }{{N}[t_{0},t_{1}]}
	+  \lp{\chi_{r<\frac{1}{2}t} F}{\LE^{2-\frac{\gamma}{4}}[t_{0},t_{1}]}  \ , \\
	\lp{F}{M_{1}[t_{0},t_{1}]} \  &=  \ \sum_{k+|J| \leqslant 1}
	\lp{(\tau_{+}\tau_0 \partial_u)^k(r \td \partial_{x})^J F }{{N}[t_{0},t_{1}]}
	+  \lp{\chi_{r<\frac{1}{2}t} F}{\LE^{2-\frac{\gamma}{4}}[t_{0},t_{1}]}  \ ,
\end{align*}
Finally, we have the following initial data spaces 
which will applied to $\nabla \phi(0)$:
\begin{equation}
		\lp{f}{H^{s,a}_{k}} \ = \ \sum_{\substack{|I|\leqslant s,\ |J|\leqslant k}}
		\lp{\jap{r}^{a+|J|}\nabla^{I+J}_x f}{L^2_x} \ . \notag
\end{equation}


\subsection{Statement of the main results}

\begin{thm}[Main Theorem] \label{main_thm}

Assume that $\mathcal{M}$ is of the form $\mathbb{R}\times \mathbb{R}^3$ and
that the metric $g$ satisfies the decay assumption \ref{ex_norm_coords}
with $0< \gamma< \delta$, the non-trapping assumption \ref{non-trapping}, and
the local energy decay assumption \ref{ls_est}. We then have the following uniform
estimates for all $T>0$:
\begin{enumerate}[I)]
	\item(Conformal Energy Estimate)
\begin{align}
	&\sup_{0\leqslant t\leqslant T}\lp{\phi(t)}{\CE}  \ + \ \lp{\phi}{S[0,T]} \ \lesssim \
	\lp{\nabla\phi(0)}{H_{0}^{0,1}} \ + \ \lp{\Box_{g}\phi}{{\ell_{t}^{1}N[0,T]}}  \ . \label{conf_energy_est}
\end{align}

\item(Conformal Energy Estimate With Vector fields)
\begin{align}
	&\sup_{0\leqslant t\leqslant T}\lp{\phi(t)}{\CE_{1}} \ + \ \lp{\phi}{S_{1}[0,T]} \  \lesssim \
	\lp{\nabla\phi(0)}{H^{0,1}_{1}}   +   \lp{\Box_{g}\phi}{\ell_{t}^{1}N_{1}[0,T]} \ . \label{conf_energy_est_vf}
\end{align}

\item(Global Pointwise Decay)
\begin{align}
	\lp{\tau_{+}^{\frac{3}{2}}\tau_{0}^{\frac{1}{2}}\phi}{L^{\infty}[0,T]} \  \lesssim \
	\lp{\nabla\phi(0)}{H^{0,1}_{1}} \ + \  \lp{\Box_{g}\phi}{\ell_{t}^{1}M_{1}[0,T]} \ . \label{point1}
\end{align}
\end{enumerate}
\end{thm}

One can think of this of this as a \emph{conditional} linear stability result for the
wave equation in non-trapping backgrounds satisfying our weak asymptotics. In other words:
for such backgrounds, as long as the LED estimate \eqref{LS1} holds,
then the pointwise decay rates $t^{-\frac{3}{2}}$ in the interior and $t^{-1}$ along
light cones for solutions to the linear wave equation \eqref{wave_eqn}
will also hold (see \cite{MR2266993}, \cite{MR2764864} for a proof of LED estimates in curved backgrounds).
In principle, the decay estimates above should be suitable
for non-linear applications since the norm $N_1$ on the source term
$\Box_{g}\phi$ should be able to handle quadratic derivative non-linearities
with a null condition. We will explore this in future work.

For our results below we will assume the initial data are smooth and
compactly supported. However, this condition can be relaxed so that
the only regularity requirements for the data are that they belong to
the weighted Sobolev spaces discussed in the statement of the main
theorem. This relaxation can be achieved by standard approximation
arguments which we omit.


\subsection{Comparison with previous literature}

To the best of the author's knowledge, this is the first work dealing
with outgoing metrics via vector fields since the stability of Minkowski space
\cite{MR1316662}, \cite{MR1946854}
and its extensions \cite{MR2531716}. We also note that not only is our
metric large and time-dependent, but we actually relax the conditions
on the causal structure of our spacetime significantly. In particular,
the function $u(t,x)$, can deviate appreciably from being a true optical function
for $(\mathcal{M},g)$ and from the Minkowski analogue $u_{mink}=t-r$.
We believe this type of setup may have some useful applications since in
practice constructing an exact optical function is a laborious process and
using a suitable replacement might be desirable. Our result also differs from previous
work in one crucial way: we build our estimates using conformal energy
instead of the Dafermos-Rodnianski p-weighted estimates as in
\cite{MR2730803}, \cite{Yang:2013kq},  \cite{MR3056620} or the fundamental solution of
Minkoswki space as in \cite{MR2921169}. Boundedness of conformal
energy is crucial since it is precisely the reason we are able to prove that solutions
in the wave zone decay pointwise at the sharp rate of $t^{-1}$. Another key difference,
at least with respect to the work of S. Yang above, is that we commute once with the scaling vector field $S$.
As mentioned in that work, it is a commonly held belief that one needs to have $t\partial_{t}g=O(1)$
in the set $r\leqslant 1$ for all time in order to use a scaling vector field as a commutator. We show here that in fact
$|t\partial_{t}g|\lesssim t^{1-\gamma}$ will suffice and that therefore the classical methods involving commuting the scaling vector field
into conformal energy still apply in this general setting. Commuting with $S$
is desirable since it leads to the higher interior decay rate of $t^{-\frac{3}{2}}$.
Since this rate of decay is integrable in time, we hope that it is useful for
some non-linear applications.

For pointwise decay via vector fields on large, time-dependent, asymptotically
flat spacetimes the only previous results are those mentioned above: \cite{Yang:2013kq},
\cite{MR3056620} and \cite{MR2921169} (see also \cite{Baskin:2012oz} for a non-vector field proof
of radiation field asymptotics for weakly decaying metrics with a full asymptotic expansion, as
well as \cite{MR2113761}, \cite{MR2259204}, \cite{MR2527808}, \cite{MR2736525}, \cite{Luk:2010fk}, \cite{MR2994507} for related applications in Black Hole spacetimes).
The work \cite{MR3056620} establishes an $L^{\infty}$ decay rate of $\jap{t+r}^{-1}$
for compactly supported metrics satisfying a non-sharp version of \eqref{LS1} for large perturbations.
One of the main differences with our results is that $\partial_{t}g$ only needs to be small
in the interior leading to more general metrics. However, the $L^{\infty}$ decay proved in that work is weaker
than ours in the interior and the metric equals Minkowski in the exterior. As an application
of his method the author also shows a small data global existence result for semilinear
equations satisfying the null condition. In the more recent work \cite{Yang:2013kq} the same author
proves a pointwise decay rate of $\jap{r}^{-\frac{1}{2}}\jap{t-r}^{-\frac{1}{2}}$
for time-dependent metrics which are uniformly close to Minkowski and decay weakly in the null outgoing region.
The main difference with our work again is that both the interior decay and the wave zone
decay achieved for the solution are weaker. Additionally, the outgoing conditions assumed for
the metric are inhomogeneous and demand more decay on the undifferentiated terms $g^{\alpha\beta}$
as well as $\Omega_{ij} g^{\alpha\beta}$. On the other hand, we point out that once again the
assumptions on the metric in the interior are slightly more general than ours and that global existence for
quasilinear equations satisfying a null condition is again shown as an application.
Lastly, in \cite{MR2921169} the authors prove a decay rate of $\jap{t+r}^{-1}\jap{t-r}^{-2}$
(Price Law) for non-trapping spacetimes with $t\partial_{t}g=O(1)$. The authors
assume a sharp LED estimate with norms similar to ours as well as wave zone decay rates $|\partial^{k}(g-m)|\lesssim \jap{r}^{-1-k}$
-- which are more restrictive than ours. Despite the fact that a lot of decay is achieved for the linear problem
the norms for the source term $\Box_{g}\phi$ involved in getting that decay do not allow for applications to
non-linear problems. However, we mention that \cite{MR2921169} also proves the Price Law
for the black hole case.


\subsection{Organization of the paper} 

In section \ref{energ_form} we recall the standard
energy formalism for the wave equation. In section \ref{conf_chan_mult}
we set up and prove a generalization of the conformal method of Lindblad-Sterbenz
\cite{MR2253534} (see also \cite{Bieri:2014kq}). This method is a general framework
for proving weighted energy estimates arising from asymptotically conformal
Killing vector fields in curved spacetimes. This framework is central to our work
since it is the foundation upon which our exterior proof of conformal energy is built. 
In a curved spacetime there's three advantages to using this
method versus the classical proof of conformal energy\footnote{See \cite{MR1820023}
for the classical proof of the conformal energy in a curved background close to Minkowski.}:
firstly, since the identities are already in divergence form we avoid having to perform
several integrations by parts in order to take advantage of special cancellations
for the boundary terms. Secondly, this method is robust enough to prove other
useful weighted energy estimates such as the fractional conformal
energy bounds we see in \cite{MR2253534}. Lastly, the method is capable of handling
the weak decay of our metrics in the wave zone. To the best of the author's knowledge,
no other method has proved capable of proving conformal energy bounds with these
types of conditions.

In the case of Minkowski space, which is the only case covered in \cite{MR2253534},
this method is motivated by the observation that
the Morawetz vector field $K_{mink}=(t^2+r^2)\partial_t+2tr\partial_r$ is conformal Killing.
Therefore it is desirable to understand how the energy formalism for multipliers changes under conformal
maps: $g\rightarrow \Omega^{-2}g$. The choices ${}^{I}\Omega=r$ and
${}^{I\!I}\Omega =t^2-r^2$ make $K_{mink}$ a \emph{Killing field} in these backgrounds.
Since the deformation errors vanish, it is a simple matter to then use the energy formalism corresponding
to the conformal metrics to prove two conjugated Morawetz estimates which,
together, combine to yield the conformal energy bound. To extend this method
to curved spacetimes we once again look at the conformal wave equation and use it to
develop a general formalism for multipliers. Inspired by the Minkowski case, we choose
smooth positive weights ${}^{I}\Omega=\jap{r}, {}^{I\!I}\Omega=\tau_{-}\tau_{+}$ which asymptotically
behave like $r$ and $t^2-r^2$. We use these in combination with a modified Morawetz vector field
$K_{0}=u^{2}\partial_{u}+2(u+r)\td \partial_{r}$ which is asymptotically Killing with respect to
these conformal backgrounds. Given this input, the generalized Lindblad-Sterbenz machinery
established in our work effectively reduces the bulk of the proof of \eqref{conf_energy_est}
to a multiplier bound modulo error terms. It is important to use both of these
weights in in our method since ${}^{I}\Omega$ degenerates where $r\sim 0$ but
${}^{I\!I}\Omega$ is well-behaved there and (locally) controls the bulk of the conformal energy.
In the wave zone the opposite behavior takes place and it is the weight
${}^{I}\Omega$ that is responsible for the bound on conformal energy.

The preceding method requires a positive-definite energy
density associated to $\partial_t$. This is addressed in section \ref{nosuper} by proving
a general result stating that non-trapping, plus smallness of $\partial_{t}g$, plus
asymptotic flatness implies that the vector field $\partial_{t}$ is uniformly timelike.
For our types of metrics this implies that $\partial_t$ is uniformly timelike in
the asymptotic region $t\gg 1$. The proof of this fact is by contradiction: if $\partial_{t}g$
is small and $\partial_t$ is close to null then the inner product
$\langle \partial_t, \partial_t \rangle \approx 0$ and is almost conserved along
the null geodesic flow. However, by the non-trapping condition and
asymptotic flatness $\langle \partial_{t}, \partial_{t} \rangle\approx -1$ in the far exterior -- a fact
which contradicts the previous claim. In section \ref{bondi_coord_sect} we derive all the
identities for error terms in Bondi coordinates as well the corresponding asymptotic bounds.
In section \ref{reduction_sec} we reduce to proving our results in the asymptotic
region $t\gg 1$ and set up some notation for absorbing small errors there.
In section \ref{conf_ener_sect} we use of our Lindblad-Sterbenz machinery
to prove \eqref{conf_energy_est}. The proof of the error bounds for this estimate rely
on three main building blocks: firstly, an upgrade of \eqref{LS1} to a t-weighted LED
estimate in the interior which is used to control all large deformation
tensor errors supported within this region. Secondly, the fact that $\partial_t$
is timelike plus the dominant energy condition gives us a coercive bound for the weighted energies
we wish to control. Thirdly, in the wave zone the $K_0$ field is set-up so that the deformation tensors
yield better spacetime errors compared to the standard Minkowski Morawetz field $K_{mink}$.
In short, the Lindblad-Sterbenz formalism coupled with the hierarchy of decay \eqref{mod_coords}
suffices to control these error terms.

In section \ref{comm_sect} we prove the higher conformal bound \eqref{conf_energy_est_vf}.
We do this by commuting the equation once with the Lie algebra
$\mathbb{L}$ -- in particular we avoid the use of the Lorentz boosts
$\Omega_{0i}=t\partial_{i}-x^{i}\partial_{t}$. In the wave zone the
desired estimates follow from the fact that the modified vector
fields have favorable errors that work well with the renormalization \eqref{mod_coords}.
In the interior the main problem is that \eqref{local_nontrapping} implies $|t\partial_{t}g|\lesssim t^{1-\gamma}$, thus
commuting with the scaling vector field is non-trivial. We fix this in stages: we
start by proving some commutator estimates. After applying Hardy estimates
to the ensuing lower order terms, the main errors arising from commuting with $S$
consist of T-weighted $L^{2}$ terms with two derivatives
supported inside $ r<c t $ with $c\ll \frac{1}{2}$. We control these by proving a t-weighted LED estimate with
vector fields and use it to trade $t\partial_t\phi$ for $S\phi$ plus small errors. This leaves only terms with two spatial derivatives
$t\nabla^{2}_{x}\phi$ to be bounded. Thanks to the global weighted $L^{2}$ elliptic estimate \eqref{elliptic},
we are able to trade
two space derivatives for the elliptic part of the wave operator $P=|g|^{-\frac{1}{2}}\partial_{i}|g|^\frac{1}{2}g^{ij}\partial_{j}$.
We then trade $tP$ for $\nabla\Gamma\phi+t\Box_{g}\phi+ t\nabla\partial_{t}\phi$. This method is somewhat
reminiscent of the work of Klainerman-Sideris \cite{MR1374174} and relies crucially on the global weighted $L^{2}$ elliptic
estimate \eqref{elliptic} which is shown in the appendix.
In this procedure it is convenient to use the $L^{2}$ norm (instead of the dyadic LED norms) when commuting with $S, \Omega_{ij}$
because we ultimately need to resort to an $L^{2}$ Hardy estimate in order to deal with
the lower order terms. Only then can we apply the elliptic estimate \eqref{elliptic}
to close the argument outlined above. Additionally, using the $L^{2}$ norms for
these terms is advantageous since it also sets up the estimates so we can re-use them
in the proof of the global pointwise decay in section 6 which follows by a similar
type of argument.


\subsection{Basic notation}

The following notation will be used in the sequel:
\begin{itemize}
	\item We denote $A \lesssim  B$ (resp. ``$A\ll B$"; ``$A\approx B$") if $A \ \leqslant \ C B$ for some fixed $C>0$ which may 
	change from line to line (resp. $A\leqslant \epsilon B$ for a small $\epsilon$; both $A\lesssim B$ and $B\lesssim A$).
	\item By default, any norms involving a range for the $t$ variable have $t\in [t_{0},t_{1}]$ with $0\leqslant t_0<t_1$.
	\item Given norms $\lp{\cdot}{A},   \lp{\cdot}{B}$ and a weight $\omega$, the notation $B\subseteq \omega A$
	means $\lp{\omega f}{A}\lesssim \lp{f}{B}$.
	\item The notation $m=m^{-1}=diag(-1,1,1,1)$ denotes the Minkowski metric in $(t,x)$ coordinates.
	\item The notation $\eta^{-1}$ denotes the inverse Minkowski metric in $(u_{mink},x^i)$ coordinates with $u=u_{mink}= t-r$:
		\begin{equation}
			\eta^{uu}\ = \ 0 \ , \qquad
			\eta^{ui} \ = \ -\omega^i \ , \qquad
			\eta^{ij} \ = \ \delta^{ij} \ , \qquad \omega^i=\frac{x^i}{r} \ .\notag
		\end{equation}
	\item We denote $\td \partial_{x}=\partial_{x}|_{u=const}$
	in Bondi coordinates. We also denote $\td \partial_r = \omega^{i}\td \partial_{i}$.
	\item $\Omega_{ij}:=\{\Omega_{ij}\}_{ i<j }$ denotes the modified rotations.
	\item The weights used in our norms are:
	\begin{equation*}
			\tau_-=\la u\ra = (1+u^2)^{\frac{1}{2}}\ , \qquad \tau_+=C+ u+2r \ , \qquad \tau_0=\frac{\tau_-}{\tau_+}
			\qquad \tau=\frac{\tau_{+}}{\jap{r}} \ , 
	\end{equation*}
	\item The conformal weights are ${}^{I}\Omega=\jap{r},{}^{I\!I}\Omega=\tau_{-}\tau_{+}$
	and their corresponding rescaled solutions are: $\psi=\Omega\phi$, ${}^{I}\psi={}^{I}\Omega\phi$,
	and ${}^{I\!I}\psi={}^{I\!I}\Omega\phi$, respectively.
	\item $D$ will denote the covariant derivative corresponding to the metric $g$.
	\item $N$ will denote the future-directed unit normal vector to the level sets $t=const$:
		\begin{equation}
			N^{\alpha}=\frac{-g^{\alpha\beta}\partial_{\beta}t}{(-g^{\alpha\beta}
			\partial_{\alpha}t\partial_{\beta}t)^{\frac{1}{2}}} \label{normalv}\ . 
		\end{equation}
\end{itemize}


\section{Preliminary Setup}\label{prelim_setup}


\subsection{Energy formalism}\label{energ_form}

Here we recall the basic energy setup for vector
fields multipliers and commutators for problem \eqref{wave_eqn}.

\subsubsection{Vector field multipliers}

Define the \emph{Energy-Momentum Tensor associated to $(g,\phi)$} by:
\begin{equation*}
	T_{\alpha\beta} \ = \ \partial_{\alpha}\phi\partial_{\beta}\phi  - 
	\frac{1}{2}g_{\alpha\beta}
	\partial^{\mu}\phi\partial_{\mu}\phi \ .
\end{equation*}
$T_{\alpha\beta}$ is related to the wave operator $\Box_g\phi$ by the identity $D^{\alpha}T_{\alpha\beta} =
\Box_{g}\phi\cdot \partial_{\beta}\phi$, with $D$ denoting the covariant derivative. Given a smooth vector field $X$ we define the 1-form $
{}^{(X)}\!P_{\alpha}  =  T_{\alpha\beta}X^{\beta}$. Taking the divergence of this we arrive at the well-known formula:
\begin{equation}
	D^{\alpha}{}^{(X)}\!P_{\alpha} \ = \ \Box_{g}\phi\cdot X\phi 
	  +  \frac{1}{2}{}^{(X)} \pi^{\alpha\beta}T_{\alpha\beta} \ , \label{stokes_1}
\end{equation}
with $\mathcal{L}_{X}g_{\alpha\beta} =  {}^{(X)}\pi_{\alpha\beta}$. The symmetric 2-tensor ${}^{(X)}\pi$ is the \emph{Deformation Tensor} of 
$g$ with respect to $X$ and measures the change of $g$ under the flow generated by $X$. Integrating \eqref{stokes_1} over the time slab $
\{(t,x)| \ t_0\leq t\leq t_1\}$ and using Stokes' theorem we get the \emph{Multiplier Identity}:
\begin{equation}
		\int_{t=t_0} \!{}^{(X)}\!P_{\alpha}N^{\alpha} \,
		d^{\frac{1}{2}}dx   - \int_{t=t_1} \!{}^{(X)}\!P_{\alpha}N^{\alpha} \,
		d^{\frac{1}{2}}dx \ = \  \int_{t_0}^{t_1}\!\!\!\int_{\mathbb{R}^{3}}
		\big( \Box_{g}\phi\cdot X\phi  +  \frac{1}{2}{}^{(X)}
		\pi^{\alpha\beta}T_{\alpha\beta}\big)\, dV_{g} \ , \label{stokes_2}
\end{equation}
where $dV_g=d^{\frac{1}{2}}dtdx$ and $N$ is the vector defined in Eq. \eqref{normalv}.
The integrand on the left hand side is the \emph{Energy Density} associated to $X$
through the foliation by spacelike hypersurfaces $t=const$. We also
recall that $T_{\alpha\beta}$ obeys:

\emph{The Dominant Energy Condition}: For any two timelike,
future-directed vector fields $X,Y$ the energy momentum tensor $T_{\alpha\beta}$ satisfies $|\nabla\phi|^{2} \ \lesssim \ T(X,Y)$.

\subsubsection{Formulae for commutators and multipliers}\label{form_comm}

Given a vector field $X$, we define the \emph{Normalized Deformation Tensor} of $X$
to be:
\begin{equation*}
	{}^{(X)} \widehat{\pi} = {}^{(X)} \pi - \frac{1}{2}g \cdot trace({}^{(X)}\!\pi) \ .
\end{equation*}
This tensor is present in some of the most important formulas dealing with vector fields for the wave equation.

\begin{lem}[Basic formulas involving ${}^{(X)} \widehat{\pi}$]\label{basic_iden_lemma}
Let $\phi$ be a smooth function and $X$ a vector field. The following identities hold:
\begin{subequations}
\begin{align}
		{}^{(X)} \widehat{\pi}^{\alpha\beta} \ &= \ -d^{-\frac{1}{2}} X(d^{\frac{1}{2}} g^{\alpha\beta})
		- g^{\alpha\beta }\partial_\gamma X^\gamma + g^{\alpha\gamma}\partial_\gamma X^\beta
		+ g^{\beta\gamma}\partial_\gamma X^\alpha \ , \label{pi_hat_formula} \\
		[\Box_g,X]\phi \ &= \ D_\alpha \! {}^{(X)} \widehat{\pi}^{\alpha\beta}D_\beta\phi
		\ + \ (D_\gamma X^\gamma)\Box_g\phi \ , \label{X_comm_formula}\\
		D^\alpha {}^{(X)}\! P_{\alpha} \ &= \ \Box_g\phi\cdot X\phi\ + \ \frac{1}{2}{}^{(X)} \widehat{\pi}^{\alpha\beta}
		\partial_\alpha\phi \partial_\beta\phi \ , \label{divergence_formula}
\end{align}
where \eqref{pi_hat_formula} is computed in local coordinates.
\end{subequations}
\end{lem}

\begin{proof}
We'll prove each of these formulas separately.

\Part{1}{The identity \eqref{pi_hat_formula}}
In local coordinates:
	\begin{align}
		{}^{(X)}\pi^{\alpha\beta} \ = \ -X(g^{\alpha\beta})
		+g^{\alpha\gamma}\partial_\gamma X^\beta + 
		g^{\beta\gamma}\partial_\gamma X^\alpha \ . \label{well_known}
	\end{align}
Subtracting the expression $\frac{1}{2} g^{\alpha\beta}  (trace{}^{(X)}{\pi})
=g^{\alpha\beta}d^{-\frac{1}{2}}\partial_{\gamma}(X^{\gamma}d^{\frac{1}{2}})$ 
from both sides gives the result.

\Part{2}{The identity \eqref{X_comm_formula}}
In local coordinates we have the following well-known formula (see \cite{MR2666888}):
\begin{equation*}
	[\Box_g,X]\phi\ = \ {}^{(X)}\pi^{\alpha\beta}D_{\alpha\beta}^{2}\phi+(D_{\alpha}\!\! {}^{(X)}\pi^{\alpha\beta})
	\partial_{\beta}\phi-\frac{1}{2}\partial^{\alpha}(trace{}^{(X)}{\pi})\partial_{\alpha}\phi \ .
\end{equation*}
Using the definition of ${}^{(X)}\widehat{\pi}$:
\begin{align*}
	{}^{(X)}\pi^{\alpha\beta}D_{\alpha\beta}^{2}\phi  \ &=  \ {}^{(X)}\widehat\pi^{\alpha\beta}D_{\alpha\beta}^{2}\phi+
	\frac{1}{2} (tr{}^{(X)}{\widehat\pi}) \Box_{g}\phi \ , \\
	(D_{\alpha}\!\!{}^{(X)}\pi^{\alpha\beta})\partial_{\beta}\phi
	\ &= \ (D_{\alpha}\!\!{}^{(X)}\widehat\pi^{\alpha
	\beta})\partial_{\beta}\phi+\frac{1}{2}\partial_{\alpha}(tr{}^{(X)}\widehat\pi)\partial^{\alpha}\phi \ .
\end{align*}
Identity \eqref{X_comm_formula} follows by combining these two statements.

\Part{3}{The identity \eqref{divergence_formula}} Combining \eqref{stokes_1} with the identity:
\begin{align*}
	T_{\alpha\beta}{}^{(X)}{\pi}^{\alpha\beta} \ = \ \partial_\alpha\phi\partial_\beta\phi
	{}^{(X)}{\pi}^{\alpha\beta} - \frac{1}{2} g^{\alpha\beta}\partial_\alpha\phi
	\partial_\beta\phi \cdot trace ({}^{(X)}{\pi}) \ = \  \partial_\alpha\phi
	\partial_\beta\phi {}^{(X)}\widehat{\pi}^{\alpha\beta} \ , 
\end{align*}
gives the result.
\end{proof}


\subsection{Conformal changes for vector field multipliers}\label{conf_chan_mult}

As mentioned in the introduction, our goal in this section is to record how all the formulae
associated with the vector field multiplier method change under conformal deformations of the metric.


\subsubsection{The conformal wave equation} Let $ g_{\alpha\beta} $ be a Lorentzian
metric on an $1+3$ dimensional spacetime. We consider 
a conformally equivalent metric $ \td{g}_{\alpha\beta} $ where
$\Omega^2 \td{g} =g $ for some smooth weight function $\Omega>0$.
Let $\td D$ denote the covariant derivative of $\td{g}$
and $\Box_{\td{g}}=\td D^\alpha\td D_\alpha$ the corresponding wave
operator. We then have the following standard formula from geometry:

\begin{lem} [Conformal wave equation]\label{conf_dividen_lem} Let $\psi  =  \Omega\phi$ and $\Box_{g}\phi=F$.
Then, for the wave operator $\Box_{\td{g}}$ of the conformal metric $\td{g}=\Omega^{-2}g$ we have:
\begin{align}
		\Box_{\td{g}}\psi + V\psi \ = \ \Omega^{3} F 
		\ , \qquad V \ = \ \Omega^{3} \Box_g\Omega^{-1} \ . \label{trans_eq}
\end{align}
\end{lem}

\begin{proof} Since $\Box_{\td{g}} \ = \ 
\Omega^{2}( \Box_g - 2g^{\alpha\beta}\partial_\alpha (\ln\Omega)\partial_\beta )$,
rescaling $\Box_{g}\psi=\Box_g\Omega\phi$ yields:
\begin{align}
	\Box_g\psi  
	\ = \  \Omega F+2 g^{\alpha\beta}\partial_{\alpha}(\ln\Omega)
	\partial_{\beta}\psi-W\phi \ , \label{first_form_box}
\end{align}
with:
\begin{equation*}
	W  \ = \ -\Box_g\Omega +2\Omega\partial_{\alpha}(\ln
	\Omega)\partial^{\alpha}(\ln\Omega)
	\ = \  \Omega^{2}\Box_g \Omega^{-1} \ .
\end{equation*}
\end{proof}


\subsubsection{Conformal vector field multipliers} Let $\chi(t,x)$ a smooth cutoff function
and $\psi=\Omega\phi$. Using equation \eqref{trans_eq}
we define the \emph{Conformal Energy-Momentum Tensor for $(\psi,\chi)$}:
\begin{equation}
		\td{T}_{\alpha\beta}^{\chi} \ = \ \partial_\alpha\psi\partial_\beta\psi  - 
		\frac{1}{2}\td{g}_{\alpha\beta}(\td{g}^{\mu\nu}\partial_\mu
		\psi\partial_\nu\psi - \chi V\psi^2) \ , \qquad V \ = \ \Omega^{3} \Box_g\Omega^{-1}  \ .
		\label{em_tens_cut} 
\end{equation}
This satisfies the divergence law:
\begin{equation}
		\td{D}^\alpha \td{T}_{\alpha\beta}^{\chi} \ = \ ((\chi-1)V\psi+\Omega^{3} F)\partial_\beta\psi
		 +   \frac{1}{2}\partial_\beta (\chi V)
		\psi^2 \ . \label{div_law_conf}
\end{equation}
Let ${}^{(X)}\! \td{P}^{\chi}_\alpha  = \td{T}^{\chi}_{\alpha\beta}
X^\beta$. We have the following identities:
\begin{lem} [Conformal multiplier identity] \label{confdivlem}Let $X$ be a vector field, $\chi(t,x)$ be a smooth cutoff
function and $\Box_{g}\phi=F$. Then, ${}^{(X)}\! \td{P}^{\chi}_\alpha$ satisfies the identity:
	\begin{align}
		\int_{ t=t_0}\!\! {}^{(X)} \td{P}^{\chi}_\alpha N^{\alpha} \Omega^{-2} \
		d^{\frac{1}{2}} dx  -  \int_{ t=t_1} \!\! {}^{(X)} \td{P}^{\chi}_\alpha N^{\alpha} \Omega^{-2} \
		d^{\frac{1}{2}} dx =   \int_{t_0}^{t_1}\!\!\!\int_{\mathbb{R}^{3}} 
		 \td{D}^\alpha {}^{(X)} \td{P}^{\chi}_\alpha \Omega^{-4} \ dV_{g}
		\ , \label{div_identity1_cut} 
	\end{align}
	where $dV_g=d^{\frac{1}{2}}dtdx$ and $N$ given by Eq. \eqref{normalv}.
	Additionally, for the divergence on the RHS we have the identity:
\begin{equation}
		\td D^\alpha {}^{(X)}\!\td{P}^{\chi}_\alpha \Omega^{-4}  = 
		F \cdot \Omega^{-1}X\psi + \Omega^{-2} A^{\alpha\beta}\partial_\alpha\psi
		 \partial_\beta\psi +  B^{\chi}\phi^2 
		 +  C^{\chi}  \phi\, \Omega^{-1}X\psi \ , \label{td_div_iden_cut3} 		
\end{equation}
with:
\begin{align}
		&{}^{(\Omega,X)}A  =  \frac{1}{2}  \big({}^{(X)} \widehat{\pi}
		+2X\ln (\Omega) g^{-1}\big)  , \label{AB_formulas1} \\
		&{}^{(\Omega,X)} B^{\chi}   =     \frac{1}{2\Omega^{2}}
		\big(  X(\chi V)- \hbox{trace}(A)
		\chi V\big)  \ ,   \label{AB_formulas2} \\
		&{}^{(\Omega,X)} C^{\chi}  =   \frac{1}{\Omega^{2}}(\chi - 1)V \ . 
		\label{AB_formulas3}
\end{align}
On the RHS of the last four lines above all contractions are computed with respect to $g$.
\end{lem}

\begin{proof}
Identity \eqref{div_identity1_cut} follows immediately by integrating $\td{D}^\alpha {}^{(X)}\! \td{P}^{\chi}_\alpha$ with 
respect to the volume form $dV_{\td g}=\sqrt{|\td{g}|}dtdx=\Omega^{-4}dV_g$ over the set
$t_0\leqslant t\leqslant t_1$ and applying Stokes' theorem. It remains to compute the divergence $\td D^{\alpha}{}^{(X)} \!
\td{P}^{\chi}_\alpha$. Using Eq. \eqref{divergence_formula}
we have:
\begin{equation}
		\td{D}^\alpha {}^{(X)}\!\td{P}^\chi_\alpha  =  
		\frac{1}{2} ( \widehat {\mathcal{L}_X \td{g}})^{\alpha\beta} \partial_\alpha\psi\partial_\beta\psi
		- \frac{1}{4}\hbox{trace} ( \widehat {\mathcal{L}_X \td{g}})\chi V\psi^2 
		+\big( (\chi-1)V\psi
		+\Omega^3 F\big)X\psi + \frac{1}{2}X (\chi V)
		\psi^2 \ , \label{div_P_iden}
\end{equation}
where all the contractions are computed with respect to $\td{g}$.
To compute the first two RHS terms we use the identities:
\begin{equation}
		\mathcal{L}_{X}\td{g}  
		 =  \Omega^{-2}\big(\mathcal{L}_X g - 2X\ln(\Omega)g \big) \ ,
		\qquad 
		 \widehat {\mathcal{L}_X \td{g}}   =  
		\Omega^{-2}\big(\widehat{\mathcal{L}_X g} 
		+ 2 X\ln (\Omega) g\big) \ , \label{lind_ster1}
\end{equation}
Substituting the last  line into RHS \eqref{div_P_iden} gives us
\eqref{td_div_iden_cut3} and \eqref{AB_formulas1} -- \eqref{AB_formulas3}.
\end{proof}

\begin{rem}
Choosing $\Omega=1$ and $\chi\equiv 0$ we recover the classical
multiplier identities in section \ref{energ_form}. In the Minkowski
space, choosing $g=m$, $\chi\equiv 0$ and weights
${}^{I}\Omega=r,{}^{I\!I}\Omega=t^{2}-r^{2}$ gives $\mathcal{L}_{X}
\td{g}\equiv 0$ in line \eqref{lind_ster1} and $V\equiv 0$ in Eq.
\eqref{trans_eq}. Thus we recover the conformal energy
formalism in the work of Lindblad-Sterbenz \cite{MR2253534}.
\end{rem}

\begin{rem}
The quantity $\td{P}_\alpha N^{\alpha}$ will denote
$\td{P}^{\chi}_\alpha N^{\alpha}$ with $\chi\equiv0$.
\end{rem}




\subsection{No superradiance}\label{nosuper}

In order to produce the coercive bound
$|\nabla\psi|^{2}\lesssim {}^{(\partial_t)} \td{P}_\alpha N^{\alpha}$
inside the set $r\leqslant 1$, the conformal multiplier method
above requires the vector field $\partial_{t}$ to be uniformly timelike.
Since $(g-m)$ could remain large inside this set
as $t\rightarrow +\infty$, we have no reason, a priori, to expect
this condition to hold. To address this issue, we will prove
below that given some mild conditions on
the metric, the vector field $\partial_t$ is uniformly timelike everywhere.
We start with some preliminary lemmas:

\begin{lem}[Coercive bound for energy on null-geodesics] Let $T$ be a past-directed, uniformly timelike vector field
and $L$ be a future-directed null vector field with $L^{\alpha}=g^{\alpha\beta}\xi_{\beta}$. The following uniform bound holds:
		\begin{align}
			\langle T,L\rangle \ \approx \ \lp{\xi}{\ell^{2}(\mathbb{R}^{4})} \ . \label{normal_null}
		\end{align}
\end{lem}

\begin{proof}
It suffices to show $\langle T,L\rangle  \gtrsim \lp{\xi}{\ell^{2}(\mathbb{R}^{4})}$.
Since $T$ is uniformly time-like we can construct a (local) system of coordinates $\{\partial_0,\partial_i\}$ such that:
\begin{equation*}
	\partial_{0}=T \ , \qquad g_{00}\ =\ -1 \ , \qquad g_{0i} \ = \ 0 \ ,  \qquad g_{ij}Y^{i}Y^{j} \ \approx \ \delta_{ij}Y^{i}Y^{j}
	\ =\ \lp{Y}{\ell^{2}(\mathbb{R}^{3})}^2 \ .
\end{equation*}
By our hypotheses we have $\langle T,L\rangle >0$. Therefore, in this system of coordinates $\langle T,L\rangle=g^{0\beta}\xi_{\beta}=-\xi_0>0$. 
Since $L$ is null we have:
\begin{align*}
	g^{00}(\xi_{0})^{2} \ + \ g^{ij}\xi_{i}\xi_{j} \ = \ 0 \qquad
	\Rightarrow \qquad -\xi_{0} \ = \ (g^{ij}\xi_{i}\xi_{j})^{\frac{1}{2}}
	\ \approx \ (\sum_{i}\xi_{i}^{2})^{\frac{1}{2}} \ .
\end{align*}
The lemma follows.
\end{proof}

\begin{lem} [Exponential bounds for null-geodesic coefficients] \label{est_nullgeo} Let $(\mathcal{M},g_{\alpha\beta})$
be asymptotically flat and suppose the vector field $N$ given by Eq. \eqref{normalv}, is uniformly
timelike and future-directed. Additionally, assume that $g_{\alpha\beta}$ satisfies the quantitative
non-trapping assumption \ref{non-trapping}.
Let $\dot \gamma=\xi$ be an outgoing, future-directed null-geodesic given in $(t,x^i)$ coordinates by $\xi$ with affine parameter $s$
satisfying $\dot{\gamma}t\big|_{s=0}=1$ with $\dot{\gamma}s\equiv 1$, $\gamma(0)\in\{ |x|\leqslant R_0 \}$.
Then, there exists a constant $A(R_0)>0$ such that for all $t>0$:
		\begin{align}
			\lp{\xi(0)}{\ell^{2}(\mathbb{R}^{4})} A^{-1} \ \lesssim \ \lp{\xi(t)}{\ell^{2}(\mathbb{R}^{4})} \ \lesssim \ 
			\lp{\xi(0)}{\ell^{2}(\mathbb{R}^{4})} A \ . \label{nullgeo_bounds}
		\end{align}
\end{lem}

\begin{proof}
Let $\lambda >0$ and $\lp{\xi(0)}{\ell^{2}(\mathbb{R}^{4})} =1$ and $s$ be an affine parameter.
By assumption \ref{non-trapping} there exists $C_\lambda>0$ such that for all $s> C_{\lambda}$,
$\gamma(s)\in \{ |g-m| < \lambda\}$. Therefore, by choosing $\lambda$ small it suffices to prove the result for
the range $s\in [0,C_{\lambda}]$. For this we will use $(t,x^i)$ coordinates. We claim that there exists a constant $k(R_0)>0$ such that for all $t>0$:
		\begin{align}
			\lp{\xi(0)}{\ell^{2}(\mathbb{R}^{4})} \exp(-kt) \ \lesssim \ \lp{\xi(t)}{\ell^{2}(\mathbb{R}^{4})} \ \lesssim \ 
			\lp{\xi(0)}{\ell^{2}(\mathbb{R}^{4})} \exp(kt) \ . \label{exp_nullgeo_bounds}
		\end{align}
To prove the claim we consider the Hamiltonian formulation for the geodesic flow. Let $p(x,\xi)=g^{\alpha\beta}\xi_{\alpha}\xi_{\beta}$ be the 
principal symbol for $\Box_g$. The Hamiltonian flow obeys the equations:
\begin{align*}
	\frac{dx^{\alpha}}{ds} \ = \ \partial_{\xi^{\alpha}} p\ , \qquad \frac{d\xi_{\alpha}}{ds} \ = \ -\partial_{x^{\alpha}}p \ .
\end{align*}
Since $-N$ is uniformly timelike, past-directed and $L^{\alpha}=g^{\alpha\beta}\xi_{\beta}$ is null, line \eqref{normal_null} implies:
\begin{align*}
	\frac{dt}{d s} \ = \ \partial_{\xi^0}p \ = \ 2g^{0\beta}\xi_{\beta} \ 
	= \ 2\langle -N,L\rangle  \ \approx \ \lp{\xi}{\ell^{2}(\mathbb{R}^{4})} \ .
\end{align*}
By chain rule: $\frac{ds}{d t}  =  (2 \langle-N,L\rangle)^{-1} 
\approx  \lp{\xi}{\ell^{2}(\mathbb{R}^{4})}^{-1}$. This leads to:
\begin{align*}
	&\frac{dx^{0}}{dt} \ = \ 1 \ , \qquad \frac{dx^{i}}{dt} \ = \ \frac{dx^i}{d s} \frac{ds}{d t}
	\ =\ \frac{\partial_{\xi^{i}}p}{2g^{0\beta}\xi_{\beta}} \ = \ O(1)\ , \\
	&\frac{d\xi_{\alpha}}{dt} \ = \ -\frac{\partial_{x^{\alpha}}p}{2g^{0\beta}\xi_{\beta}} \ \approx \
	O( \lp{\xi}{\ell^{2}(\mathbb{R}^{4})}) \ .
\end{align*}
Therefore the $x^{\alpha}$ change with constant speed. For the frequency $\xi_{\alpha}$, since the last statement holds for all $\alpha$,
there exists a $k>0$ sufficiently large so that:
\begin{align*}
	- k\lp{\xi}{\ell^{2}(\mathbb{R}^{4})} \ \leqslant \ \frac{d}{dt} \lp{\xi}{\ell^{2}(\mathbb{R}^{4})}
	\ \leqslant \ k\lp{\xi}{\ell^{2}(\mathbb{R}^{4})} \ .
\end{align*}
Integrating these bounds finishes the proof of the claim. Combining
with our initial remarks yields the lemma.
\end{proof}

\begin{lem}
Let $(\mathcal{M},g)$ satisfy all the hypotheses of Lemma \ref{est_nullgeo}.
Choose $R_0>0$ and denote $\xi$, $\lambda$, $A(R_0)$ as above. Then, for
$\lambda$ sufficiently small, there exists $C_1>0$ such that in the exterior region $\{ |g-m| < \lambda\}$:
\begin{align}
	|\langle\partial_t,\xi \rangle| \ > \ C_{1}A^{-1} > 0  \ . \label{almost_cons_low}
\end{align}
\end{lem}

\begin{proof}
Since $\xi$ is a null vector, choosing $\lambda$
sufficiently small and using asymptotic flatness and non-trapping gives us
$\xi =  \left(\partial_t+\omega^{i}\partial_{i}+O(\lambda)\right)\lp{\xi}{\ell^{2}(\mathbb{R}^{4})}$.
By lemma \ref{est_nullgeo} this implies
$|\langle \partial_t,\xi\rangle| > C_1 \lp{\xi}{\ell^{2}(\mathbb{R}^{4})}\geqslant C_1A^{-1}$
for some constant $C_1>0$.
\end{proof}

\begin{rem}
In general we have $C_{1}A^{-1}\ll 1$ due to the exponential nature of $A$.
\end{rem}

We now prove the main result of this section:

\begin{prop}[$\partial_t$ is uniformly timelike]\label{T_timelike} Let $(\mathcal{M},g)$ satisfy all the
hypotheses of Lemma \ref{est_nullgeo} together with $|\mathcal{L}_{\partial_{t}}g|< \epsilon $. Choose $R_0>0$ and denote
$A, C_1 > 0$ as above. Then there exists $\epsilon_{ntrap}(A)>0$ sufficiently small 
such that for all $0<\epsilon<\epsilon_{ntrap}$ the bound $\langle\partial_t,\partial_{t}\rangle \ \leqslant  \ - \frac{1}{2}C_{1}A^{-1}$
holds everywhere.
\end{prop}

\begin{proof}
Let $\epsilon>0$, choose $\lambda$ sufficiently small (as in the previous lemma)
and let $C_{\lambda}$ be the constant given in assumption \ref{non-trapping}.
The proof is by contradiction and breaks up into two cases:\\
\case{1} {$\partial_t$ becomes null} Assume
for a contradiction that there exists a $p\in\mathcal{M}$ such that $\langle\partial_t,
\partial_{t}\rangle_{p}=0$.
Let $\gamma(s)$ be the unique, affinely parametrized forward null geodesic with $\dot\gamma(0)=\partial_t \in T_p\mathcal{M}$
and $\dot{\gamma}t\big|_{s=0}=1$. For sufficiently small $\tau\in(-\tau_0,\tau_0)$
we can define a smooth 1-parameter family of curves $\gamma_{\tau}(s)$ with $
\gamma_{0}(s)=\gamma(s)$. Let $\frac{d}{ds}=\xi$, then along the null geodesic $\gamma$ we have:
\begin{align*}
	\xi \langle \partial_t, \xi \rangle \ = \  \langle \nabla_{\xi}\partial_t,\xi\rangle
	 +   \langle \partial_t,\nabla_{\xi}\xi\rangle \  &= \ -\langle 
	[\partial_t,\xi] , \xi \rangle  +   \langle \nabla_{\partial_t}\xi,
	\xi \rangle \ = \ \frac{1}{2}\mathcal{L}_{\partial_t}g(\xi,\xi) \ .
\end{align*}
Integrating this along $\gamma$ from $s=0$ to $s=s_1>0$:
\begin{align}
	\langle \partial_t, \xi \rangle_{s=s_1} \ = \ \langle \partial_t, \xi \rangle_{s=0}
	\ + \ \frac{1}{2} \int_{0}^{s_1}\mathcal{L}_{\partial_t}g(\xi,\xi)\ ds \ .  \label{almost_cons_int}
\end{align}
Since $\xi(0)=\partial_t$, 
the first term on the RHS
\eqref{almost_cons_int} satisfies $\langle \partial_t, \xi \rangle_{s=0}=0$.
Using lemma \ref{est_nullgeo}
together with the hypotheses yields the almost conservation law:
\begin{align}
	|\langle \partial_t, \xi \rangle_{s_1}| \ \leqslant \ \frac{1}{2}A^{2} \ \epsilon \cdot s_1 \ . \label{almost_cons}
\end{align}
Choosing $\epsilon_{ntrap}  < 2 C_{1}A^{-1}
(A^{2} \cdot 2C_\lambda)^{-1}$ contradicts line \eqref{almost_cons_low}
when $s_{1}= 2C_\lambda$.

\case{2} {$\partial_t$ is close to null} By asymptotic flatness, together with continuity of
$g_{\alpha\beta}$ and the results in case 1, we may assume that $\langle \partial_t,
\partial_t \rangle<0$. Suppose for a contradiction that there exists a point $p\in\mathcal{M}$
such that $0>\langle\partial_t, \partial_{t}\rangle_{p}>-\frac{1}{2}C_{1}A^{-1}$. 
By computing in an appropriate local coordinate system near $p$ we can find
an outgoing null vector $\overline{\xi}\in T_{p}\mathcal{M}$
such that $\langle \partial_t, \overline{\xi} \rangle_{p}= \langle\partial_t, \partial_{t}\rangle_{p}$.
Proceeding in the same manner as before we then get:
\begin{align}
	\langle \partial_t, \xi \rangle_{s=s_1} \ = \ \langle \partial_t, \xi \rangle_{s=0}
	 +  \frac{1}{2} \int_{0}^{s_1}\mathcal{L}_{\partial_t}g(\xi,\xi)\ ds \ .  \label{almost_cons_int2}
\end{align}
Since $\xi(0)=\overline\xi$, our initial remarks imply:
\begin{align*}
	|\langle \partial_t, \xi \rangle_{s=0}| \ = \ |\langle \partial_{t},
	\overline\xi \rangle_{p}  | \ \leqslant \ \frac{1}{2}C_{1}A^{-1} \ .
\end{align*}
Combining this with line \eqref{almost_cons_int2} and using lemma \ref{est_nullgeo}
together with our hypotheses:
\begin{align}
	|\langle \partial_t, \xi \rangle_{s_1}| \ \leqslant \ \frac{1}{2}C_{1}A^{-1}
	\ + \ \frac{1}{2}A^{2} \ \epsilon \cdot s_1 \ . \label{almost_cons2}
\end{align}
Choosing: $\epsilon_{ntrap}  < C_{1}A^{-1}(A^{2} \cdot 2C_\lambda)^{-1}$
contradicts \eqref{almost_cons_low} when $s_{1}= 2C_\lambda$.

\end{proof}

This proposition leads to two important consequences:

\begin{cor}\label{P_cor}
Let $(\mathcal{M},g)$ satisfy the conditions of the main theorem \ref{main_thm}. Then:
\begin{enumerate}[I)]
	\item There exists a $T_{ntrap}\gg 1$ such that $\partial_{t}$
	is uniformly timelike for all $t\in [T_{ntrap},+\infty)$.
	\item In $(t,x^i)$ coordinates, the operator $P(t,x,\nabla_{x})=\partial_i h^{ij}\partial_j$ with $h^{ij}=d^\frac{1}{2}g^{ij}$ is uniformly elliptic. Furthermore, inside 
		the region $\{ r< \frac{1}{2}t\} \cap [T_{ntrap},+\infty)$, $P(t,x,\nabla_{x})$ satisfies the uniform estimate:
		\begin{equation}
			|P\phi| \ \lesssim \ |\Box_g\phi|+ |\partial^{2}_t \phi|+\jap{r}^{-\delta}|\partial_{t}\nabla_{x}\phi|
			+\jap{r}^{-1-\delta}\big(|\partial_t\phi| + \tau^{-\gamma}|\nabla_{x}\phi|\big)
			\ , \label{pbound}
		\end{equation}
		with $h^{ij}=\delta^{ij}+O(\jap{r}^{-\delta})$.
\end{enumerate}
\end{cor}

\begin{proof}
\Part{1}{Statement for $\partial_{t}$} Choose $T_{ntrap}\gg 1$ sufficiently large so that
$\jap{T_{ntrap}}^{-\gamma}C<\frac{1}{10}\epsilon_{ntrap}$ holds with $C=$ the maximum
of the implicit constants in estimates \eqref{mod_coords}.
In the interior this yields $|\mathcal{L}_{\partial_{t}}g|<
\epsilon$ with $\epsilon=\frac{1}{10}\epsilon_{ntrap} $ .
In the exterior this bound follows by \eqref{mod_coords}
and the fact that $\gamma<\delta$. An application of proposition
\ref{T_timelike} then yields the result.
 
 \Part{2}{Ellipticity} By Cramer's rule we have $\det(M^{0})=\langle \partial_{t},\partial_{t}
 \rangle \det(g^{-1})$ with $M^{0}=(g^{ij})$ the $3\times 3$ matrix arising from eliminating the first row and column from the matrix 
 coefficient of $g^{-1}$. By part I we have $g_{00}<C<0$ and thus $\det(g^{-1})<C<0$ therefore it follows that
 $\det (M^{0})$ is always positive and has uniform lower bounds. Since $M^{0}$ is a $3\times 3$ matrix and $g$ has Lorentzian 
 signature it follows that $M^{0}$ has three positive eigenvectors. This proves the first claim. For \eqref{pbound} we note that:
\begin{equation}
	\Box_{g} \ = \ g^{00}\partial^{2}_{t} \ + \ d^{-\frac{1}{2}}\big(\partial_{i}g^{0i}d^{\frac{1}{2}}\partial_{t}
	 \ + \ \partial_{t}g^{0i}d^{\frac{1}{2}}\partial_{i} \ + \  P\big) \ .
\end{equation}
This, together with the asymptotic form \eqref{local_nontrapping} of the metric 
inside $r<\frac{1}{2}t$, combine to give the result.
\end{proof}



\section{Bondi Coordinates} \label{bondi_coord_sect}

\subsection{Algebraic formulas involving Bondi coordinates}

The vector fields $X\in\mathbb{L}$ written in $(u,x^i)$ 
coordinates have a very simple form:

\begin{lem}[Lie algebra property] \label{lie_lem}
In $(u,x^i)$ coordinates the vector fields defined in Eq. \eqref{mod_fields} are given by:
\begin{equation}
	\mathbb{L} \ = \ \big\{  T=\partial_u \, , \,  
	S=u\partial_u + r\td\partial_r \, , \,  \Omega_{ij}= \ x^i\td\partial_j-x^j\td\partial_i\} \ , \label{vf_bondi}
\end{equation}
and $\mathbb{L}$ forms a Lie algebra on $\mathcal{M}$.
\end{lem}
\begin{proof} 
By chain rule:
	\begin{align}
		\partial_u \ = \ (u_{t})^{-1}\partial_t \ , \qquad \td\partial_{i} \ = \ 
		\partial_{i}-(u_{t})^{-1}u_i\partial_{t} \ , \qquad \td\partial_{r} \ = \ 
		\partial_{r}-(u_{t})^{-1}u_r\partial_{t}  \label{coords_chain_rule} \ .
	\end{align}
Line \eqref{grad_u_bnd} implies $u_{t}\approx 1$ everywhere;
thus, the equations in \eqref{coords_chain_rule} are also valid everywhere.
Applying this to \eqref{mod_fields} yields \eqref{vf_bondi}.
From this we can compute the commutators to be:
\begin{align}
	[\partial_{u},S]  =   \partial_{u} \ , \qquad [\partial_{u},\Omega_{ij}]   =   0 \ , \qquad 
	[S, \Omega_{ij}]   =   0 \ , \qquad   [\Omega_{ij}, \Omega_{kl}]  =   
	-\delta_{ (ik}\Omega_{jl )} \ . \label{comm_bondi_polar}
\end{align}
The Lie algebra property follows.
\end{proof}

\begin{rem}
Since $u_{t}\approx 1$ the first identity in \eqref{coords_chain_rule} implies that
$\partial_{t}\approx\partial_{u}$. This will be used often in the sequel.
\end{rem}

Next we compute the key quantities from lemmas \ref{basic_iden_lemma} 
and \ref{confdivlem} of the previous section.

\begin{lem}[Formulas for deformation tensors]

Let $\Omega$ be a smooth function and $X$ be
a vector field in Bondi coordinates. We have the following formula for the
contravariant tensor $2A$ in line \eqref{AB_formulas1}:
\begin{align}
		 &{}^{(X)} \widehat{\pi}
		+ 2 X\ln (\Omega) g^{-1} 
		  \ = \ 
		-d^{-\frac{1}{2}}\big( \mathcal{L}_X \eta^{-1} +
		(\partial_u X^u + \td \partial_r X^r + \td \partial_i\overline X^i 
		+2(\frac{X^{r}}{r}-X \ln \Omega) )\eta^{-1}
		\big) + 
		\mathcal{R}_{X}  \ , \label{polar_bondi_defiden1}
\end{align}
where $\overline X^i=X^i-\omega^i \omega_j X^j$ denotes the angular portion of $X$
and $X^r=\omega_i X^i$ the radial portion. The remainder tensor  $\mathcal{R}$ is given
by the covariant formula:
\begin{equation}
		\mathcal{R}_{X}  \ =\  -d^{-\frac{1}{2}}\big( \mathcal{L}_X(d^\frac{1}{2}g^{-1}-\eta^{-1}) 
		+(\partial_u X^u + \td \partial_r X^r + \td \partial_i\overline X^i 
		+2(r^{-1}X^{r}-X \ln \Omega) )
		 (d^\frac{1}{2} g^{-1}-\eta^{-1}) \big)
		\ . \label{polar_bondi_defiden2}
\end{equation}
\end{lem}

\begin{proof}
We start with the identities:
\begin{align}
	{}^{(X)}\widehat{\pi}^{\alpha\beta} \ &= \ -d^{-\frac{1}{2}}\mathcal{L}_{X}
	(d^{\frac{1}{2}}g^{-1})-g^{\alpha\beta }\partial_\gamma X^\gamma \ , \label{piandL}\\
	\td \partial_\gamma X^\gamma \ &= \  2 r^{-1}X^{r} +
	\partial_u X^u + \td \partial_r X^r+ \td \partial_i\overline{X}^i \ . \label{eucdiv}
\end{align}
Applying this to formula \eqref{pi_hat_formula} in Bondi coordinates:
\begin{equation}
	{}^{(X)} \widehat{\pi}
	+ 2 X\ln (\Omega) g^{-1}
	=  -d^{-\frac{1}{2}}\big( \mathcal{L}_X (d^\frac{1}{2} g^{-1})
	+ (\partial_u X^u + \td \partial_r X^r
	+ \td \partial_i\overline{X}^i+2(\frac{X^{r}}{r}-X \ln \Omega))
	d^\frac{1}{2}g^{-1}\big)   \ . \label{pol_bondi5}
\end{equation}
Adding and subtracting $d^{-\frac{1}{2}}(\mathcal{L}_X \eta^{-1}+\eta^{-1})$
on the last line gives Eq. \eqref{polar_bondi_defiden1} and Eq. \eqref{polar_bondi_defiden2}.
\end{proof}

\begin{lem}[Formulas for commutators]

Let $d=|\det (g_{\alpha\beta})|$ and $X\in \{\partial_u , \Omega_{ij}\}$ be in
Bondi coordinates. The following identities hold:
\begin{equation}
		[\Box_g,X] \  = \  D_\alpha \mathcal{R}_{X}^{\alpha\beta}D_\beta
		  +    \frac{1}{2} X\ln (d)   \Box_g \ ,  \label{X_comm_form}
\end{equation}
where:
\begin{equation}
		\mathcal{R}_{X}  \ =\ -d^{-\frac{1}{2}} \mathcal{L}_X(d^\frac{1}{2}g^{-1}-\eta^{-1}) 
		\ . \label{usual_R}
\end{equation}
For $S$ in Bondi coordinates we have:
\begin{equation}
		[\Box_g,S] \  = \  D_\alpha \mathcal{R}_{S}^{\alpha\beta}D_\beta
		\ + \  \frac{1}{2}( 4+ S \ln (d))  \Box_g \ , \label{S_comm_form}
\end{equation}
where $\mathcal{R}_{S}$ is given by formula:

\begin{equation}
		\mathcal{R}_{S}   =  -d^{-\frac{1}{2}}\big( \mathcal{L}_X(d^\frac{1}{2}g^{-1}-\eta^{-1}) 
		+(\partial_u X^u + \td \partial_r X^r + \td \partial_i\overline X^i )
		 (d^\frac{1}{2} g^{-1}-\eta^{-1}) \big)
		\ . \label{polar_bondi_defiden3}
\end{equation}


\end{lem}

\begin{proof}
For $X\in \{\partial_u , \Omega_{ij}\}$ we use formula \eqref{polar_bondi_defiden1}
with $\Omega=1$:
\begin{equation}
		{}^{(X)} \widehat{\pi}
		  =  -d^{-\frac{1}{2}} \big( \mathcal{L}_X \eta^{-1} +
		(\td \partial_\alpha X^\alpha )\eta^{-1}
		 + \mathcal{L}_X(d^\frac{1}{2}g^{-1}-\eta^{-1}) 
		+(\td \partial_\alpha X^\alpha)
		 (d^\frac{1}{2} g^{-1}-\eta^{-1}) \big)
		  \ . \label{aux_lie1}
\end{equation}
For these two vector fields we have $\mathcal{L}_X
\eta^{-1} =0$ and $\partial_u X^u=\td \partial_i X^i=0$. 
Applying Eq. \eqref{X_comm_formula} gives us \eqref{X_comm_form} and \eqref{usual_R}.
For $X=S$ one can compute:
\begin{equation*}
	\mathcal{L}_{S} \eta^{-1}  +  (\partial_u S^u + \td \partial_r S^r
	+ \td \partial_i\overline {S} {}^{i})\eta^{-1}   =  0 \ , \qquad r^{-1}S^{r} =   1
	\ , \qquad \partial_u S^u+\td \partial_i S^i  =  4 \ .
\end{equation*}
Substituting this together with \eqref{eucdiv}
into \eqref{aux_lie1} then applying \eqref{X_comm_formula}
finishes the proof of \eqref{X_comm_form} -- \eqref{S_comm_form}.
\end{proof}



\subsection{Asymptotic estimates involving Bondi coordinates}

Our first task here is to compute the decay rates for the Lie derivatives $\mathcal{L}_{X}g$ with $X\in \mathbb{L}$.
%


\begin{lem}[Basic Lie derivative estimates]\label{basic_lie_lem}
Let $X=X^\alpha\partial_\alpha$ be in Bondi coordinates.
\begin{enumerate}[I)]
	\item \label{main_LR_estimate} Suppose that $X$ satisfies the symbol-type bounds:
	\begin{equation}
		\big| (\tau_-\partial_u)^k (\jap{r}\td \partial_{x})^J X^u\big|
		\ \lesssim \ \tau_- \ , \qquad\qquad
		\big| (\tau_-\partial_u)^k (\jap{r}\td \partial_{x})^J X^i\big|
		\ \lesssim \ \jap{r} \ , \label{vect_sym_bnds}
	\end{equation}
	and obeys the conditions:
	\begin{equation}
		  \td \partial_{i} X^u \ = \ \partial_u X^i \ = \ 
		  \td \partial_{r} r^{-1}(X^i-\omega^i\omega_jX^j) 
		  \ \equiv\  0 \ .
		\label{good_X_conds}
	\end{equation}
	Let $\mathcal{R}^{\alpha\beta}$ be a contravariant two tensor
	satisfying the bounds:
	\begin{subequations}\label{R_sym_bnds}
	\begin{align}
			| \partial_u^k \td \partial_{x}^J \mathcal{R}^{ij}  | \ &\lesssim \
			 \jap{r}^{-k-|J|-\delta}\tau_{0}^{-k} \tau^{-\gamma k}
			 \label{R_sym_bndsa} \ , \\
			|\partial_u^k \td \partial_{x}^J \mathcal{R}^{ui}| \ &\lesssim \
			  \jap{r}^{-k-|J|-\delta}\tau_{0}^{\frac{1}{2}-k} \tau^{-\gamma k} 
			\ , \label{R_sym_bndsb}\\
			| \partial_u^k \partial_{x}^J( \mathcal{R}^{ui}
			\! \! - \!\! \omega^i\omega_j \mathcal{R}^{uj})| \ &\lesssim \
			  \jap{r}^{-k-|J|-\delta}\tau_{0}^{1-k} \tau^{-\gamma k}
			  \label{R_sym_bndsc}  ,  \\
			 | \partial_u^k \td \partial_{x}^J \mathcal{R}^{uu}  | \  &\lesssim \
			  \jap{r}^{-k-|J|-\delta}\tau_{0}^{2-k} \tau^{-\gamma k} 
			 \ , \label{R_sym_bndsd}
	\end{align}
	\end{subequations}
	with similar estimates for $\mathcal{R}^{iu}$. Then, the Lie derivative
	$\mathcal{L}_X\mathcal{R}$ satisfies the bounds \eqref{R_sym_bnds}
	with the exponent $-\gamma k$ above replaced by $1-\gamma (1+ k)$. 
	\item \label{R_alt2} Alternatively, if we
	substitute the condition $\partial_u X^r=0$ with $\partial_u(X^i-\omega^i X^r)=0$ and
	keep the rest of  \eqref{vect_sym_bnds} and \eqref{good_X_conds} the same,
	then the result of part \ref{main_LR_estimate}
	holds with the bound on line  \eqref{R_sym_bndsa} replaced by:
	\begin{align}
		| \partial_u^k \td \partial_{x}^J( \omega_i\omega_j \mathcal{R}^{ij})  | \ &\lesssim \
		\jap{r}^{-k-|J|-\delta}\tau_{0}^{-\frac{1}{2}-k} \tau^{1-\gamma (1+ k)}  \ ,
		\label{alt_ij_1}\\
		| \partial_u^k \td \partial_{x}^J( \mathcal{R}^{ij} - 
		\omega^i\omega^j \omega_k\omega_l \mathcal{R}^{kl})| \ &\lesssim \
		\jap{r}^{-k-|J|-\delta}\tau_{0}^{\frac{1}{2}-k} \tau^{1-\gamma (1+ k)}   \ . 
		\label{alt_ij_2}
	\end{align}
	\item \label{X0_part} Let $X=\partial_u$
	and suppose $\mathcal{R}$ satisfies \eqref{R_sym_bnds}. Then $\jap{r} \tau_0\mathcal{L}_X\mathcal{R}$
	satisfies \eqref{R_sym_bnds} as well.
\end{enumerate}
\end{lem}

\begin{proof}[Proof of Lemma \ref{basic_lie_lem}]

\Part{1}{The  $\mathcal{R}$ bounds involving condition \eqref{vect_sym_bnds}}
We begin with the proof of estimates \eqref{R_sym_bnds} for $\mathcal{L}_X\mathcal{R}$
assuming  conditions \eqref{vect_sym_bnds} and \eqref{good_X_conds} or the
alternative listed in item \ref{R_alt2} above.
The formula for the Lie derivative is:
\begin{equation*}
	\mathcal{L}_X\mathcal{R}^{\alpha\beta} = X(\mathcal{R}^{\alpha\beta})
	-\partial_\gamma (X^\alpha)\mathcal{R}^{\gamma\beta}-\partial_\gamma
	(X^\beta)\mathcal{R}^{\alpha \gamma} \ .
\end{equation*}
We check each component:\\
\case{1}{The $uu$ component}
Here we have:
\begin{equation}
		\mathcal{L}_X\mathcal{R}^{uu} \ = \ X(\mathcal{R}^{uu})-2\partial_u(X^u)\mathcal{R}^{uu}
		-\td \partial_i (X^u)(\mathcal{R}^{ui})
		-\td \partial_i (X^u)(\mathcal{R}^{iu}) \ . \notag
\end{equation}
Since $\td \partial_{i} X^u=0$, the estimate on line \eqref{R_sym_bndsd} for $\mathcal{L}_X\mathcal{R}^{uu}$
is immediate from estimates \eqref{vect_sym_bnds}--\eqref{R_sym_bnds}.

\case{2}{The $ui$ and $iu$ components}
By symmetry of the estimates on lines \eqref{R_sym_bnds}
it suffices to treat the $ui$ case. We have:
\begin{equation}
		\mathcal{L}_X\mathcal{R}^{ui} \ = \ X(\mathcal{R}^{ui})
		- \partial_u(X^u)\mathcal{R}^{ui} - \partial_u(X^i)\mathcal{R}^{uu} 
		-\td \partial_j (X^u) \mathcal{R}^{ji}
		-\td \partial_j (X^i) \mathcal{R}^{uj}  \ . \notag
\end{equation}
Using estimates \eqref{vect_sym_bnds}--\eqref{R_sym_bnds} we get a symbol bound
on the order of $\jap{r}^{-\delta}\tau_0^\frac{1}{2}\tau^{1-\gamma}$ for this term. In addition one sees that for all 
parts of the above formula 
save the expression $X(\omega^i \mathcal{R}^{ur})-\td \partial_{r} (X^i) \mathcal{R}^{ur}$ the
bound is on the order of $\jap{r}^{-\delta}\tau_0\tau^{1-\gamma}$. To see the improvement for 
$ \mathcal{L}_X \mathcal{R}^{ui}-\omega^i\omega_j \mathcal{L}_X\mathcal{R}^{uj} $
we note that the worst term is absent once one subtracts off the radial part since:
\begin{equation}
		( X(\omega^i) -\omega^i\omega_j X(\omega^j)) \mathcal{R}^{ur} 
		-\td \partial_{r} (X^i-\omega^i\omega_j X^j) \mathcal{R}^{ur} 
		=  
		-r \td \partial_{r} r^{-1} (X^i \! - \! \omega^i\omega_j X^j) \mathcal{R}^{ur}  
		 =  0 \ , \notag 
\end{equation}
where we used \eqref{good_X_conds}  for the last identity above.

\case{3a}{The $ij$ components assuming $\partial_u X^i=0$}
Here we have:
\begin{equation}
		\mathcal{L}_X\mathcal{R}^{ij} 
		\ = \ X(\mathcal{R}^{ij})-\td \partial_k (X^i)\mathcal{R}^{k j}
		-\td \partial_k (X^j)\mathcal{R}^{i k} \ . \notag
\end{equation}
The bound on line \eqref{R_sym_bndsa} for $\mathcal{L}_X\mathcal{R}^{ij} $
follows by multiplying together \eqref{vect_sym_bnds} and \eqref{R_sym_bndsa}.

\case{3b}{The $ij$ components assuming $\partial_u(X^r)\neq 0$} 
In this case we are still assuming $\partial_u(X^i-\omega^i\omega_j X^j)=0$.
Therefore:
\begin{equation}
		\mathcal{L}_X\mathcal{R}^{ij} 
		\ = \ X(\mathcal{R}^{ij})-\omega^i \partial_u (X^r)\mathcal{R}^{u j}
		-\omega^j \partial_u (X^r)\mathcal{R}^{i u}
		-\td \partial_k (X^i)\mathcal{R}^{k j}
		-\td \partial_k (X^j)\mathcal{R}^{i k} \ . \notag
\end{equation}
By \eqref{vect_sym_bnds}--\eqref{R_sym_bnds} this has a
symbol bound of order $\jap{r}^{-\delta}\tau_0^{-\frac{1}{2}}\tau^{1-\gamma}$.
On the other hand all but the second and third terms above yield a bound of order $\jap{r}
^{-\delta}\tau^{1-\gamma}$. Subtracting the radial part yields:
\begin{align}
		&\mathcal{L}_X\mathcal{R}^{ij} -\omega^i\omega^j \mathcal{L}_X\mathcal{R}^{rr}
		\notag \ = \  -\omega^i \partial_u (X^r)( \mathcal{R}^{u j} - \omega^j \mathcal{R}^{u r})
		-\omega^j \partial_u (X^r)( \mathcal{R}^{i u} - \omega^i \mathcal{R}^{r u} )
		+O(\jap{r}^{-\delta}\tau^{1-\gamma}) \ . \notag
\end{align}
By line \eqref{R_sym_bndsc} we have $O(\jap{r}^{-\delta}\tau^{1-\gamma})$ symbol bounds for the first two terms on the
RHS above as well.

\Part{2}{Estimates involving $\partial_{u}$}
Fix a dyadic region $r \tau_0\approx 2^k$.
Then if $X=\partial_{u}$ the vector field $2^k X$ satisfies conditions \eqref{vect_sym_bnds}
and all of \eqref{good_X_conds}. From the above calculations
one immediately has all of \eqref{R_sym_bnds} for $r \tau_0\mathcal{L}_X\mathcal{R}$.
\end{proof}

Along a similar vein we can derive the following set of bounds which will be needed when employing the conformal
multiplier method.

\begin{lem} \label{detandpot}Let $g_{\alpha\beta}$ be in Bondi coordinates.
	\begin{enumerate}
		\item (Estimates for the determinant) Let $d=|\det(g_{\alpha\beta})|$
		be computed in Bondi
		coordinates $(u,x^i)$. We have the symbol bounds:
			\begin{align}
				|(\jap{r}\tau_{0}\partial_u)^{k}(\jap{r}\td\partial_{x})^{J}(d^{\frac{1}{2}}-1)|
				\ \lesssim \ \jap{r}^{-\delta} \label{det_bondi_decay} \ .
			\end{align}
		\item (Estimates for the conformal potential)
		Let ${}^{I}\Omega=\jap{r}, {}^{I\!I}\Omega=\tau_{-}\tau_{+}$ and
		${}^{\Omega}V=\Omega^3 \Box_g \Omega^{-1}$. The
		potentials ${}^{I}V$, ${}^{I\!I}V$ satisfy the following
		symbol bounds:
		\begin{subequations}\label{V_bound}
			\begin{align}
				\big| (\jap{r}\tau_{0}\partial_u)^k (\jap{r}\td \partial_x)^J
				({}^{I}V)\big| \ &\lesssim \ \jap{r}^{ -\delta}
				\tau_0^{-\frac{1}{2}} \ , \label{V_bound1}\\
				\big| (\jap{r}\tau_{0}\partial_u)^k (\jap{r}\td \partial_x)^J
				({}^{I\!I}V)\big| \ &\lesssim   \
				(\tau_{-}\tau_{+})^{\frac{3}{2}}\jap{r}^{ -1-\delta} \ . \label{V_bound2}  
			\end{align}
		\end{subequations}
	\end{enumerate}
\end{lem}

\begin{proof}
\Part{1}{Determinant bounds} Follows from estimates \eqref{mod_coords}
since the determinant is a continuous function of the metric components
$g_{\alpha\beta}$.\\
\Part{2}{Potential bounds}
Let $\mathcal{R}=d^\frac{1}{2}g^{-1}-\eta^{-1}$ and write the wave operator in
Bondi coordinates as: $\Box_g =   d^{-\frac{1}{2}}(\Box_\eta  + \td \partial_\alpha 
\mathcal{R}^{\alpha\beta}\td \partial_\beta) $, where $\Box_\eta$ is the Minkowski
wave equation in Bondi coordinates. Expanding $\Box_\eta$ yields:
\begin{equation}
		\Box_\eta \ =  \ -2\partial_u\td \partial_r 
		+ \td \partial_{r}^2 - 2r^{-1}\partial_u + 2r^{-1}\td \partial_{r} 
		+ r^{-2}\sum_{i<j}(\Omega_{ij})^2 \ . \label{eta_idens}
\end{equation}
Using this identity we can compute:
\begin{equation}
	{}^{I}\Omega^3 \Box_{\eta} ({}^{I}\Omega^{-1})  = -\frac{3}{\jap{r}^{2}} ,
	\qquad {}^{I\!I}\Omega^3 \Box_{\eta} ({}^{I\!I}\Omega^{-1})  = 
	-\frac{4}{ \tau_{-}\tau_{+}} +\frac{u}{r}O(\frac{1}{ \tau_{-}\tau_{+}})  \ . \label{mink_V_lot}
\end{equation}
Thus:
\begin{align*}
	{}^{I}V   &=  -\frac{1}{d^{\frac{1}{2}}}\big(   r \partial_u  \mathcal{R}^{u r}
	+   r \td \partial_i   \mathcal{R}^{i r}-   2   \mathcal{R}^{r r}(\frac{r}{\jap{r}})^2 +\frac{3}
	{\jap{r}^{2}} \big) \ , \\
	{}^{I\!I}V   &=  -\frac{1}{d^{\frac{1}{2}}}\big( {}^{I\!I}\Omega \partial_{\alpha}
	\mathcal{R}^{\alpha\beta}\partial_{\beta}{}^{I\!I}\Omega- 
	2 \mathcal{R}^{\alpha\beta}(\partial_{\alpha}{}^{I\!I}\Omega)
	(\partial_{\beta}{}^{I\!I}\Omega)+\frac{4}{ \tau_{-}\tau_{+}}
	-\frac{u}{r}O(\frac{1}{ \tau_{-}\tau_{+}}) \big) \ . 
\end{align*}

The estimate \eqref{V_bound1} follows from this and the fact that $\mathcal{R}$
satisfies the estimates \eqref{R_sym_bnds} by assumption \ref{ex_norm_coords}. For ${}^{I\!I}V$ the two worst
cases are when $\alpha=\beta=u$ and $\alpha=u,\beta=r$. A quick count of the weights in these two cases
and an application of Eqs. \eqref{det_bondi_decay} and \eqref{R_sym_bnds} gives us the result.
\end{proof}

\begin{lem}[Formulas for boundary terms]
Let $X^r$ and $Y^u$ be positive and set $X=X^r\td \partial_r$
, $Y=Y^u\partial_u$. Then there exist a constant $C>0$ such that
the following pointwise estimates for boundary terms on the
divergence identity \eqref{div_identity1_cut} hold:
\begin{align}
		& \frac{{}^{(X)}\!\td{P}_\alpha^\chi N^{\alpha}}
		{\Omega^{2}} \ \gtrsim  \
		X^r  \big|\frac{\td\nabla_x 
		(\Omega \phi)}{\Omega}\big|^2   -C\big( \frac{X^r}{\jap{r}^{\delta}\tau_0^2}\big)
		\big|\frac{\nabla 
		(\Omega \phi)}{\Omega}\big|^2
		- C(|{}^{\Omega}V|\chi X^r ) \big|\frac{\phi}{\Omega}\big|^{2}
		  \ , \label{boundary1}\\
		&\frac{{}^{(Y)}\!\td{P}_\alpha^\chi N^{\alpha}}
		{\Omega^{2}} \  \gtrsim \  Y^u \big|\frac{\nabla 
		(\Omega \phi)}{\Omega}\big|^2
		-C(|{}^{\Omega}V|\chi Y^u)  \big|\frac{\phi}{\Omega}\big|^{2} \ . \label{boundary2} 
\end{align}
\end{lem}

To prove estimate \eqref{boundary1} we will need the following elementary result:

\begin{lem}[Approximate null frame] \label{app_null}
Let $X$ and $Y_A$, $A=1,2$ be vector fields satisfying
$\sup_\alpha |X^\alpha|\approx 1$ and $\sup_\alpha |Y_A^\alpha|\approx 1$
for their components. Suppose that there exists $\omega>0$
such that $\la X,X\ra =O(\omega^2)$, $\la X,Y_A\ra=O(\omega)$, and 
in addition $| \la Y_A,Y_B\ra -\delta_{AB}|\ll 1$. Then there exists
a null frame $\{L,\bL,e_A\}$ with $\la L,L\ra=\la L,e_A\ra=\la \bL, \bL
\ra=0$, $\la L,\bL\ra=-1$,  and 
$\la e_A,e_B\ra=\delta_{AB}$, and coefficients $c_X^A$ and $c_A^B$ for $A,B=1,2$, 
and $\gamma$, such that:
\begin{equation}
		X \ = \ \theta L +   c^A_X e_A + \gamma \bL \ , \qquad
		Y_A \ = \   c_A^B e_B \ , \notag
\end{equation}
where,
\begin{equation}
		\theta \gamma=O(\omega^2)
		\ , \qquad
		c^X_A = O(\omega) \ ,
		\qquad | c^A_B - \delta^A_B|\ll 1 \ . \notag 
\end{equation}
\end{lem}

\begin{proof}
Choose $e_A$ to be an orthonormal basis for the spacelike two plane spanned by $Y_A$ and
let $c^B_A$ be the change of basis. Then $|c^A_B - \delta^A_B|\ll 1$. Let $L,\bL$
be, respectively, the outgoing and incoming null generators over span $e_A$ with $\la L,\bL\ra=-1$. We have:
$X=\theta L + c^A_X e_A + \gamma \bL$ for some set of coefficients $\theta, c_X^A,\gamma$.
From $\la X,Y_A\ra=O(\omega)$ we have $\la X,e_A\ra=O(\omega)$ and 
so $c^A_X=O(\omega)$. Thus $\la X,X\ra = -2\theta \gamma + O(\omega^2)$ and 
so $\theta \gamma=O(\omega^2)$ follows from $\la X,X\ra=O(\omega^2)$.
\end{proof}

\begin{proof}[Proof of \eqref{boundary1}]
By choosing $C$ sufficiently large it suffices to prove the result in the wave zone
$t\approx r$ with $r\gg 1$. Consider the local basis $\{\td \partial_r,Y_A\}$
where $Y_A$ is a (local) euclidean ONB on the spheres $r=const,u=const$. We now check the hypotheses of
the preceding lemma. Since the metric $g$ is asymptotically 
flat $|\la Y_A, Y_B\ra -\delta_{AB}|\ll 1$. On the other hand
by the asymptotic formulas \eqref{mod_coords} and Cramer's 
rule we have  $\langle \td \partial_r , \td \partial_r \rangle=O(\jap{r}^{-\delta}\tau_0^2)$ and
$\langle \td \partial_r,Y_A \ra = O(\jap{r}^{-\delta}\tau_0)$. Additionally, inside this region,
assumption \ref{ex_norm_coords} together with Prop. \ref{T_timelike}
imply $|\langle\td \partial_{r},\underline{L}
\rangle|=|\theta|\approx 1$. An application of the previous lemma then gives us,
with $\omega=\jap{r}^{-\frac{{\delta}}{2}}\tau_0$:
\begin{equation}
		\td \partial_r \ = \ L + O(\jap{r}^{-\frac{{\delta}}{2}}\tau_0)\td \snabla_x
		+ O(\jap{r}^{-\delta}\tau_0^2)\nabla
		  \ , \label{dr_null}
\end{equation}
where $L$ is (outgoing) null and $\td \snabla_x$ denotes derivatives  tangent to $u=const,r=const$ which
also lie in the null plane generated by $L$. Let $\td{T}$, $\psi$ be as in Eq. \eqref{em_tens_cut}
with $\chi\equiv 0$. Since $N$ is uniformly timelike and future-directed:
\begin{align*}
		\td{T}(L,N) \ &\approx  \ |L\psi |^2 + |\td\snabla_x\psi|^2 \ , \\
		|\td{T}(\td \snabla_x,N)| \ &\lesssim \ |\td \snabla_x\psi| \cdot |\nabla \psi| + 
		|g^{\alpha\beta}\partial_\alpha\psi \partial_\beta\psi| \ , \\
		|\td{T}(\nabla,N)| \ &\lesssim \ |\nabla \psi|^2 \ .
		 \notag
\end{align*}
Using the bounds \eqref{mod_coords} and Young's inequality with $c>0$:
\begin{equation}
		\jap{r}^{-\frac{{\delta}}{2}}\tau_0 |g^{\alpha\beta}\partial_\alpha\psi
		\partial_\beta\psi| \ \lesssim \ c|\td\nabla_x \psi|^2 + 
		c^{-1} \jap{r}^{-\delta}\tau_0^2|\nabla\psi|^2 \ . \label{quad_boot}
\end{equation}
Choosing $c\ll 1$, applying Eq. \eqref{dr_null} together with the last three
inequalities and absorbing the small $|\td\nabla_x \psi|^2$ term gives:
\begin{equation}
		\td{T}(\td \partial_r,N ) \ \gtrsim \ |\td\nabla_x\psi|^2
		+O( \jap{r}^{-\delta}\tau_0^2)
		|\nabla\psi|^2 \ . \notag
\end{equation}
Adding the undifferentiated terms and using the definition
of $\td T^{\chi}$ finishes the proof.\\
\end{proof}

\begin{proof}[Proof of \eqref{boundary2}]
The vector field $N$ is uniformly timelike and future-directed by our initial assumptions.
The vector field $\partial_t$ is uniformly timelike and future-directed by
proposition \ref{T_timelike}. The proof of inequality \eqref{boundary2} then follows from  
$\partial_u\approx \partial_t$ and the dominant energy condition.
\end{proof}


\section{Additional Notation and Preliminary Reduction}\label{reduction_sec}

In this section we reduce to the asymptotic
region where $t\gg 1$. We first set up some notation
which will be used below to absorb small errors
inside this region.

\begin{defn}[Description of $\mu$, $\epsilon$, $T^{*}$ and $I^{*}$] \label{bspara}
We make the following definitions:
\begin{enumerate}[a)]
	\item Let $0<\mu <\frac{1}{2}$ be sufficiently small so that $\mu\cdot 
		\sup_{j}\! C_{j} \! \ll \!\frac{1}{2}$ where $C_{j}$ are the implicit constants
		in all the estimates in Lemma \ref{lemma_two_der} and lines
		\eqref{dtdx_l2}, \eqref{dxdx_l2}.
	\item Choose $0<\epsilon\ll \min \{\gamma, \mu\}$ satisfying the
		following property: for any estimate in the sequel, of the form
		$A   \leqslant  C_0 (B   +   \epsilon A) $ with
		absolute constant $C_0>0$, the number $\epsilon$ is 
		small enough that we can absorb $C_0\epsilon A$ on the LHS
		to yield the bound $A  \leqslant  2C_0 B$.
	\item Let $T^{*}(\epsilon)>T_{ntrap}\gg 1$ be sufficiently large so that the following holds:
		\begin{align}
	 		0 < (\epsilon T^{*})^{-\frac{\gamma}{100}}  <  \epsilon \ .  \label{tstar_cond}
		\end{align}
	\item Let $I^{*}=[T^{*},\infty)$ with $T^{*}$ as above.
\end{enumerate}

\end{defn}

The constant $\epsilon$ depends only on
$\{\lp{g}{C^{2}}, \gamma,\delta, \mu \}$ and
on the implicit constants in the assumptions
of the main theorem. In principle, we can
choose an explicit $\epsilon$ satisfying the
property above. However, as we only use
this constant to close a finite number of estimates below,
it is neither necessary nor particularly useful to keep track of its size.
We also note that $\epsilon$ will usually arise from a small gain in $t$ power in our
estimates, with the only exception being the small interior wedge
in the proof of estimate \eqref{error_density_int2}. On the other hand,
the purpose of $T^{*}$ is to give us an explicit lower bound for
$t$ which help us produce $\epsilon$ via inequality \eqref{tstar_cond} when
$t\gg 1$.

In order to take advantage of the setup above we now reduce to
the asymptotic region $t\in I^*$. In the sequel it suffices to show
that estimates \eqref{conf_energy_est}, \eqref{conf_energy_est_vf},
and \eqref{point1} hold for all $t\in I^{*}$ with constants that do not depend on $t$.
This is a straightforward consequence of local energy estimates.

\begin{lem}[Reduction to the asymptotic region $t\in I^*$] \label{reduction}
Estimates \eqref{conf_energy_est}, \eqref{conf_energy_est_vf}, and \eqref{point1} hold for all $t\in [0,T^{*}]$
\end{lem}
\begin{proof}
By local energy estimates both \eqref{conf_energy_est} and \eqref{conf_energy_est_vf} follow 
in the range $t\in [0,T^{*}]$ with constants that depend on $T^{*}$. Estimate \eqref{point1} follows similarly
after an application of the $L^{\infty}-L^{2}$ Sobolev embedding.
\end {proof}



\section{Conformal Energy Estimate} \label{conf_ener_sect}
In this section we prove the conformal energy estimate \eqref{conf_energy_est}.
This bound will form the basis for the higher regularity estimate \eqref{conf_energy_est_vf}
as well as the global $L^{\infty}$ decay \eqref{point1}. In the asymptotic region
$t\in I^{*}$ the conformal energy estimate will follow from:

\begin{thm}\label{conf_energy_th}
Assume the hypotheses of the main theorem hold. Then, for all $[t_0, t_1]\subset I^{*}$
the following inequalities hold:
\begin{enumerate}[I)]
	\item(T-weighted $L\!E\!D$ estimate in timelike regions) 
	\begin{align}
		& 	\lp{\chi_{r<\frac{1}{2}t}\phi}{\ell^{\infty}_t \LE^1}  \  \lesssim \
		  \lp{ \chi_{r<\frac{1}{2}t}\Box_{g}\phi}{\ell^{\infty}_t N}
		  +    \sup_{t_0\leqslant t\leqslant t_1}\lp{\phi(t)}{\CE}  \ . \label{t_weight_ls}
	\end{align}
	\item(Uniform boundedness)
		 \begin{equation}
			\sup_{t_0\leqslant t\leqslant t_1}\lp{\nabla\phi(t)}{L_{x}^{2}} + 
			\lp{ \td \nabla_x   \phi }{\NLE^{0,-\frac{1}{2}}}
			 +  \lp{\phi}{ {\LE}}\ \lesssim \  \lp{\nabla\phi(t_0)}{L_{x}^{2}}   + 
			\lp{\Box_{g}\phi}{ {\LE}^{*}}\ . \label{unif_bound}
		\end{equation}
	\item(Conformal energy estimate with interior error) 
	For any $0 < c < 1$:
		\begin{equation}
			\sup_{t_0\leqslant t\leqslant t_1}\lp{\phi(t)}{\CE}  +  
			 \lp{\phi}{{}^{I} C\!H}  + \lp{\phi}{{}^{I\! I}\CH}
			\ \lesssim \ \lp{\phi(t_0)}{\CE}  +   c^{-1}\lp{\Box_{g}\phi}{\ell^{1}_tN}
			+ c\lp{\chi_{r<\frac{1}{2}t}\phi}{\ell^{\infty}_t \LE^1} \ .    \label{k_0_output}
		\end{equation}
		
\end{enumerate}
\end{thm}

Let's show how the conformal energy estimate follows from this:

\begin{proof} [Proof of estimate \eqref{conf_energy_est}]  By definition \ref{bspara}
we can choose a $c$ small enough satisfying $\epsilon \ll c$, yet smaller
than the implicit constants in estimates \eqref{t_weight_ls} -- \eqref{k_0_output}.
Taking an appropriate linear combination of \eqref{t_weight_ls} and \eqref{k_0_output}:
	\begin{multline}
		\sup_{t_0\leqslant t\leqslant t_1}\lp{\phi(t)}{\CE}+
		 \lp{\phi}{{}^{I} C\!H}  + \lp{\phi}{{}^{I\! I}\CH}
		+ c\lp{\chi_{r<\frac{1}{2}t}\phi}{\ell^{\infty}_t \LE^1}\\
		\ \lesssim \ \lp{\phi(t_0)}{\CE} + 
		c^{-1}\lp{\Box_{g}\phi}{\ell^{1}_t N}
		+c\sup_{t_0\leqslant t\leqslant t_1}\lp{\phi(t)}{\CE} \ . \label{conf_raw}
	\end{multline}
Since $c\ll 1$ we can bootstrap the
last term above on to the LHS and close the estimate. Inside $[0,T^{*}]$ the result
follows from Lemma \ref{reduction}.
\end{proof}


The proof of the t-weighted LED estimate is modular and does not depend on anything other than assumption 
\ref{ls_est}. Let's prove this estimate right now.

\begin{proof}[Proof of estimate  \eqref{t_weight_ls}]
We apply the LED bound \eqref{LS1}
to $2^{k} \chi_k(t)\chi_{r<\frac{1}{2}t}\phi$ where $\chi_{k}$ are
a series of dyadic cutoffs supported where $t\approx 2^{k}$ and $\chi_k(t_0)=0$.  
Commuting with $\Box_{g}$ gives us:
\begin{equation}
		\lp{\tau_{+} \chi_k\chi_{r<\frac{1}{2}t} \phi}{ \LE}
		\ \lesssim \ \lp{\td{\chi}_k\td{\chi}_{r<\frac{1}{2}t}(\nabla\phi,\jap{r}^{-1}\phi)}{\LE^{*}}
		 +  \lp{\tau_+ \chi_k\chi_{r<\frac{1}{2}t}\Box_g\phi }{\LE^{*}}\   , \label{t-LE_raw}
\end{equation}
where $\td{\chi}$ denote cutoffs with slightly larger support.
For the first  RHS term and for fixed   $N,C>0$ there exists an
$M=M(C,N)>0$ and a uniform implicit constant such that:
\begin{multline}
		\lp{\td{\chi}_k \td{\chi}_{r<\frac{1}{2}t}(\nabla \phi,\jap{r}^{-1}\phi) }{\LE^*}
		\ \lesssim \ M\lp{\td{\chi}_k \jap{r}^{\frac{1}{2}+N}\! \tau_+^{-N}\! 
		\td{\chi}_{r<\frac{1}{2}t} (\nabla\phi, \jap{r}^{-1} \phi) }{L^2} \\
		+ \sum_{j<k-C} 2^{j-k}
		\lp{\chi_j(r) \td{\chi}_k(t)\tau_+ 
		\td{\chi}_{r<\frac{1}{2}t} \phi}{\LE} \ .
		\notag
\end{multline}
Plugging this into RHS \eqref{t-LE_raw} and taking $\ell_t^\infty$ for the resulting bound over 
a collection of finitely overlapping $\chi_k(t)$ and choosing $C$ large enough to absorb the LE
error from the last line above gives us:
\begin{equation}
		 \lp{ \chi_{r<\frac{1}{2}t}\phi}{\ell^\infty_t \LE^{1}} 
		\  \lesssim  \  \lp{\jap{r}^{\frac{1}{2}+N}\! \tau_+^{-N}\! \chi_{r<\frac{1}{2}t} 
		 (\nabla\phi,  \jap{r}^{-1}\phi) }{\ell^\infty_t L^2 }
		 +  \lp{ \chi_{r<\frac{1}{2} t}\Box_g\phi}{\ell^\infty_t \LE^{*,1}} \ .  \notag
\end{equation}

Choosing $N=\frac{1}{2}$, using the
support property and the definition of the norms:

\begin{equation*}
	\lp{\jap{r} \tau_+^{-\frac{1}{2}}\! \chi_{r<\frac{1}{2}t} 
		 (\nabla\phi, \jap{r}^{-1} \phi) }{\ell^\infty_t L^2 }
		 \ \lesssim \sup_{t_0 \leqslant t\leqslant t_1}
		 \lp{\phi(t)}{\CE} \ . 
\end{equation*}
This finishes the proof.

\end{proof}

The rest of this section is devoted to the proof of the estimates \eqref{unif_bound} and \eqref{k_0_output}.
This will be done over the course of the next three subsections.

\subsection{Some preliminary estimates}

Here we establish a number of technical estimates needed in the proof of
Theorem \ref{conf_energy_th}. Each argument below is self-contained.

\begin{lem}[Hardy estimates]  For test functions $\phi$ we have the following fixed-time bound:
\begin{align}
	\lp{r^{a-1} \phi(t)}{L^2_x }  \ &\leqslant \frac{2}{2a+1}\    \lp{r^a \partial_r \phi (t)}{L^2_x }
	\ , & -\frac{1}{2}< a< \infty \ . \label{hardy0}
\end{align}
Additionally let $\chi_{\jap{r}< 2}$, $\chi_{\jap{u}< 2}$ be smooth cutoff functions 
supported on the sets $ \jap{r}< 2 $, $ \jap{u}< 2 $ respectively. For all $t\in I^{*}$ one has the fixed-time bounds:
\begin{align} 
	\lp{\chi_{\jap{r}< 2} \tau_{+}\frac{\phi(t)}{r}}{L^{2}_{x}}
	 \ &\lesssim \ \lp{\phi(t)}{{}^{I\!I}\CE}  +  \lp{\nabla\phi(t)}{L^{2}_{x}}
	  +  \lp{\chi_{\jap{r}\sim 2} \tau_{+}\frac{\phi(t)}{r}}{L^{2}_{x}}  \ , 
	 \label{hardy_int}\\
	\lp{  \chi_{\jap{u}< 2}\phi(t)}{L_{x}^{2}}  \ &\lesssim \  \lp{\phi(t)}{{}^{I}\CE} 
	 +  \lp{\nabla\phi(t)}{L^{2}_{x}}
	 +   \lp{\chi_{\jap{u}\sim 2}\phi(t)}{L_{x}^{2}}\label{hardy_near_cone} \ .
\end{align}
\end{lem}

\begin{proof}[Proof of estimate \eqref{hardy0}]
For a fixed value of the angular variable we have the integral identity:
\begin{equation}
		(2a+1)\int_{0}^\infty r^{2a} \phi^2 dr  \ = \ 
		-2\int_{0}^\infty r^{2a+1}\phi\partial_r\phi dr \ . \label{int_identity}
\end{equation}
As long as $2a+1>0$ estimate
\eqref{hardy0} follows from integration of this identity in the angular
variable  and Cauchy-Schwartz.

\end{proof}


\begin{proof}[Proof of estimate \eqref{hardy_int}]
We apply \eqref{hardy0} with $a=0$ to $\lp{\chi_{\jap{r}< 2} r^{-1}\tau_{+}\phi}{L^{2}_{x}}$.
Since $t\approx \tau_{+}\tau_{0}$
inside this set we get:
	\begin{align}
		\lp{\chi_{\jap{r}< 2} \tau_{+}\frac{\phi(t)}{r}}{L^{2}_{x}} \ \lesssim \ \lp{\chi_{\jap{r} < 2}
		(\tau_{+}\tau_{0} \partial_r\phi) (t)}{L^{2}_{x}}
		\ + \  \lp{\chi_{\jap{r}\sim 2}\tau_{+}\frac{\phi (t)}{r}}{L^{2}_{x}} \ . \label{hardy_int_half}
	\end{align} 
Within this region we have $\big|({}^{I\!I}\Omega)^{-1}\partial_r({}^{I\!I}\Omega)\big|  
\leqslant 2 (\tau_{+}\tau_{0})^{-1}$. Combining this with Young's inequality:
	\begin{align}
		|\tau_{+}\tau_{0}\partial_r\phi|^{2} \ \leqslant \  24 
		\big( |\tau_{+}\tau_{0}({}^{I\!I}\Omega)^{-1}\partial_r
		({}^{I\!I}\Omega\phi)|^{2} \ + \ \phi^{2}\big) \ . \label{poincare2}
	\end{align}
Applying this inequality to the first term on the RHS\eqref{hardy_int_half} yields:
\begin{align}
	\lp{\chi_{\jap{r} < 2} (\tau_{+}\tau_{0} \partial_r\phi)(t)}{L^{2}_{x}}  \ \lesssim \
	\lp{\chi_{\jap{r} < 2}\tau_{+}\tau_{0} \frac{\partial_r({}^{I\!I}\Omega\phi)(t)}{ {}^{I\!I}\Omega}}{L^{2}_{x}}  +  \lp{\chi_{\jap{r} < 2} \phi(t)}{L^{2}_{x}}
	\ . \label{hardy_int_half2}
\end{align}
For the last term on RHS\eqref{hardy_int_half2} 
we use the support property followed by
\eqref{hardy0} with $a=0$ to get:
\begin{equation*}
	\lp{\chi_{\jap{r} < 2} \phi}{L^{2}_{x}} \ \lesssim  \
	\lp{\nabla\phi}{L^{2}_{x}} \ + \  \lp{\chi_{\jap{r}\sim 2} r^{-1}\phi}{L^{2}_{x}} \ .
\end{equation*}
This finishes the proof of estimate \eqref{hardy_int}.
\end{proof}

\begin{proof}[Proof of estimate \eqref{hardy_near_cone}]
For a fixed value of the angular variable we have the integral identity:
\begin{equation}
	-\int_{0}^{\infty} \!\! u_{r}\phi^2\chi_{\jap{u}< 2} r^2dr  \ = \ \int_{0}^{\infty}
	\! \big(2u\partial_r\phi \cdot \phi\chi_{\jap{u}< 2}+ 2ur^{-1}\phi^2\chi_{\jap{u}< 2}+
	 \phi^2\chi_{\jap{u}\sim 2} \big) r^2dr  
	\ = \ \mathcal{E}_{1}\!+\mathcal{E}_{2}+
	 \mathcal{E}_{3} \ .  \label{123steps}
\end{equation}
Integrating Eq. \eqref{123steps} in the angular variable and using inequality
\eqref{grad_u_bnd} on the LHS above gives us:
\begin{multline} 
	(1\!+\!O(r^{-\delta}))\int_{\mathbb{S}^{2}}\! \int_{0}^{\infty} \!\! \phi^2\chi_{\jap{u}< 2} \ r^2dr d\omega
	\ = \ -\int_{\mathbb{S}^{2}}\! \int_{0}^{\infty} u_{r}\phi^2\chi_{\jap{u}< 2} \ r^2dr d\omega \\
	 = \int_{\mathbb{S}^2} |\mathcal{E}_{1}|  d\omega+\int_{\mathbb{S}^2} |\mathcal{E}_{2}| 
	 d\omega+\int_{\mathbb{S}^2} |\mathcal{E}_{3}|  d\omega \label{123steps_aux}
\end{multline}
Since $r\gg 1$ in this set, we may absorb the term with $O(r^{-\delta})$ as a
small bootstrap error. Next we take absolute values on RHS\eqref{123steps_aux} and bound each term
separately.

For $\int | \mathcal{E}_{1} | d\omega$ we go back to Bondi derivatives via
$\partial_{r}=\td \partial_{r}+u_r\partial_{u}$, then apply $\big|({}^{I}\Omega)^{-1}
\partial_r({}^{I}\Omega)\big|  \lesssim  \jap{r}^{-1}$ to conjugate by ${}^{I}\Omega$.
Using Young's inequality then gives us, with $0<c\ll1$:
\begin{equation*}
	\int_{\mathbb{S}^2} | \mathcal{E}_{1}| \ d\omega  \
	 \lesssim \  c^{-1}  \big(\lp{\chi_{\jap{u} <2} \tau_{+}\tau_{0} {}^{I}\Omega^{-1}(\partial_u({}^{I}\Omega\phi),\td \partial_{r}({}^{I}
	\Omega\phi))} {L^{2}_{x}}^{2}  
	+   \lp{\chi_{\jap{u}< 2} r^{-1}\tau_{+}\tau_{0}\phi}{L_{x}^{2}}^{2}\big)+
	c \lp{\chi_{\jap{u}< 2}\phi}{L_{x}^{2}}^{2} \ .
\end{equation*}
The last term above can be bootstrapped to LHS\eqref{123steps_aux} by choosing
$c$ sufficiently small. For the next-to-last term on the RHS above we observe that
inside the set where $\jap{u}<2$, the condition $\tau_{+}\tau_{0}=O(1)$ holds.
After an application of estimate \eqref{hardy0} with $a=0$ we get:
\begin{equation}
	\lp{\chi_{\jap{u}< 2} r^{-1}\tau_{+}\tau_{0}\phi}{L_{x}^{2}}^{2} \ \lesssim \
	\lp{\nabla\phi}{L^{2}_{x}}^2+  \lp{\chi_{\jap{u}\sim 
	2}\phi}{L_{x}^{2}} \ . \label{123steps_aux2}
\end{equation}
For $\int |\mathcal{E}_{2}| d\omega$
we apply Young's inequality and get:
\begin{equation*}
	\int_{\mathbb{S}^2} |\mathcal{E}_{2}|  d\omega\ \lesssim \ c^{-1} \lp{\chi_{\jap{u} <2}
	r^{-1}\tau_{-}\phi}{L_{x}^{2}}^{2} \ + \ c \lp{\chi_{{\jap{u}< 2}}\phi}{L_{x}^{2}}^{2} \ .
\end{equation*}
Choosing $c$ sufficiently small allows us to bootstrap
the small error term to LHS\eqref{123steps_aux}, while the other term
is addressed directly via \eqref{123steps_aux2}.

Since the term $\int |\mathcal{E}_{3}| d\omega$ is acceptable as part of the RHS,
we combine the last few lines and take square roots to finish the proof of estimate
\eqref{hardy_near_cone}.
\end{proof}

\begin{lem}[Estimates for undifferentiated boundary terms] \label{estbdr2} Let $\chi_{r<\frac{1}{2}t}(r/t)$
be a smooth cutoff supported on the wedge $ r<\frac{1}{2}t $. Let ${}^{I}V$, ${}^{I\!I}V$
denote the potentials in Lemma \ref{detandpot}.2. For test functions $\phi$
and for all values $t\in I^{*}$ one has the fixed-time bounds:
\begin{align}
	 \lp{ \tau_{-}^{2}\chi_{r<\frac{1}{2}t}  ({}^{I}\Omega)^{-2} \cdot  {}^{I}V \phi^{2}(t)}
	 {L^{1}_{x}}  \  &\lesssim  \  \lp{\phi(t)}{{}^{I\!I}\CE}^2 \ + \  \epsilon \cdot
	 \lp{\phi(t) }{\CE}^2 \ ,  \label{error_density_int1} \\
	 \lp{  \tau_{-}^{2}\chi_{r<\frac{1}{2}t}  ({}^{I\!I}\Omega)^{-2} \cdot  {}^{I\!I}V \phi^{2}(t)}
	 {L^{1}_{x}} \ &\lesssim  \  \epsilon \cdot \lp{\phi(t)}{\CE}^2 \ .
	 \label{error_density_int2}
\end{align}
\end{lem}

\begin{proof}

\Part{1}{Proof of estimate \eqref{error_density_int1}} 
Let $\chi_{r<    c t}(r/t)$
be a smooth cutoff to the region $r<  c t $ with $\chi'_{r<  c t}
=r^{-1}\chi_{r\sim  c t}$ and $c \ll 1$ small. We do separate proofs for the regions
$r<  c t $ and $ \frac{  c t}{2}<r< \frac{t}{2} $. Substituting $a=-\delta/2$
into estimate \eqref{hardy0} gives us $C=2/(1-\delta)$. Using this together with
the bound \eqref{V_bound1} yields:
	\begin{multline}
	 	 \lp{ \tau_{-}^{2} \chi_{ r<   ct}  ({}^{I}\Omega)^{-2} \cdot  {}^{I}V \phi^{2}(t)}
		 {L^{1}_{x}} \  \lesssim  \ \lp{  \chi_{ r<   ct}
		r^{-1-\frac{\delta}{2}}t\phi}{L^{2}_{x}}^{2} \\
		\  \leqslant  \  \frac{2}{1-\delta} \big(\lp{  \chi_{ r<  ct} 
		r^{-\frac{\delta}{2}} t \partial_{r}\phi}{L_{x}^{2}} +C\lp{\chi_{r\sim c t}
		r^{-1-\frac{\delta}{2}}t\phi }{L^{2}_{x}}^{2} \big)   
		 \ = \    \mathcal{E}_{1}  +  \mathcal{E}_{2} \ . \label{cutoff_err1} 
	\end{multline}
For $\mathcal{E}_{1}$ we multiply line \eqref{poincare2} by
$\chi_{r<  ct} r^{-\frac{\delta}{2}}$ and integrate to get:
\begin{align}
	\mathcal{E}_{1} \  &\leqslant  \  \frac{48}{1-\delta}
	\big(\lp{ \chi_{r<  ct} r^{-\frac{\delta}{2}}t({}^{I\!I}\Omega)^{-1}\partial_r
	({}^{I\!I}\Omega\phi)}{L^{2}_{x}}^{2} \ + \   \lp{ (r/t) \chi_{r<  ct}
	r^{-1-\frac{\delta}{2}}t\phi}{L^{2}_{x}}^{2} \big) \notag \\
	& \leqslant  \ C\lp{t \chi_{r<  ct}  ({}^{I\!I}\Omega)^{-1}
	\partial_r({}^{I\!I}\Omega\phi)}{L^{2}_{x}}^{2} \ + \
	\frac{48}{1-\delta}\cdot c  \lp{ \chi_{r<  ct} r^{-1-\frac{\delta}{2}}t\phi}{L^{2}_{x}}^{2} \ , \notag
\end{align}
where in the last line we used $r/t < c$. By choosing $c\ll 1$ sufficiently small
we can bootstrap out last term above onto the second term of the first line in estimate
\eqref{cutoff_err1}.

For $\mathcal{E}_{2}$ we note that the support is where $r\sim ct $, so we can use \eqref{tstar_cond} to get
$r^{-\frac{\delta}{2}} \lesssim \big(ct \big)^{-\frac{\delta}{2}}\lesssim \epsilon$ which implies:
	\begin{align}
		\mathcal{E}_{2} \ \lesssim \ t^{-\frac{\delta}{2}} \lp{ r^{-1}t\phi}{L^{2}_{x}}^{2}
		 \ \lesssim \ \epsilon \cdot \lp{\phi(t)}{\CE}^2 \ . \notag 
	\end{align}
This proves estimate \eqref{error_density_int1} inside $r< ct $. In the
complement we apply \eqref{V_bound} plus the estimate \eqref{tstar_cond}
to get $r^{-\frac{\delta}{2}} \lesssim \big( c t \big)^{-\frac{\delta}{2}}\lesssim \epsilon$
and finish the proof.
 
\Part{2}{Proof of estimate \eqref{error_density_int2}} we do separate proofs
 for the regions $r<\epsilon t$ and $\epsilon t/2<r<t/2 $. Inside
 $r<\epsilon t$ we use the bound \eqref{V_bound2} and $r/t<\epsilon$ to get:
\begin{align*}
	\lp{  \tau_{-}^{2}\chi_{r< \epsilon t}  ({}^{I\!I}\Omega)^{-2} \cdot  {}^{I\!I}V \phi^{2}(t)}
	 {L^{1}_{x}}\
	\lesssim \ \lp{(r/t)^{\frac{1}{2}}\chi_{r<\epsilon t} r^{-1-\frac{\delta}{2}}t\phi}{L^{2}_{x}}^{2}
	\  \lesssim  \ \epsilon  \lp{\phi}{\CE}^2 \ .
\end{align*}
For the term supported where $\epsilon t/2 <r< t/2 $ the bound follows by using
\eqref{V_bound2}, then \eqref{tstar_cond} to produce $r^{-\frac{\delta}{2}} \lesssim
(\epsilon t)^{-\frac{\delta}{2}}\lesssim \epsilon$.
\end{proof}

\begin{lem}[Conjugation removal]
For all $t\in I^{*}$ we have the fixed-time bounds:
	\begin{align}
		\lp{\phi(t)}{\CE} \ \approx \   \lp{\phi(t)}{{}^{I}\CE} \ + \
		\lp{\phi(t)}{{}^{I\!I}\CE}\ + \ \lp{\nabla\phi(t)}{L_{x}^{2}}  
		\label{conjlemma} \ .
	\end{align}
\end{lem}

\begin{proof}
It suffices to show $\lp{\phi(t)}{\CE}  \lesssim
\lp{\phi(t)}{{}^{I}\CE}+\lp{\phi(t)}{{}^{I\!I}\CE}+\lp{\nabla\phi(t)}{L_{x}^{2}}$.
This reduces to proving:
\begin{align}
	\lp{\tau_{+} (\td\partial_r\phi, \tau_{0} \partial_u\phi)(t)}{L_{x}^{2}}^2 \ 
	&\lesssim \ \lp{\tau_{+} \big( \frac{\td\partial_r({}^{I}\Omega\phi)}{{}^{I}\Omega},  
	\tau_{0} \frac{\partial_u({}^{I}\Omega\phi)}{{}^{I}\Omega},  \frac{\phi}{r}\big) (t)}{L_{x}^{2}}^{2}
	   \ , \label{rem_conf_claim1}\\
	\lp{\frac{\tau_{+}\phi(t)}{r}}{L^{2}_{x}}^{2} \ &\lesssim \ \sum_{J=I,I\!I}
	\lp{ \tau_{+}(\frac{\td\partial_r({}^{J}\Omega\phi)}{{}^{J}\Omega}, 
	\tau_{0} \frac{\partial_u({}^{J}\Omega\phi)}{{}^{J}\Omega} )(t)}{L_{x}^{2}}^{2}+\lp{\nabla\phi(t)}{L_{x}^{2}}^{2} \ . 
	\label{rem_conf_claim2}
\end{align}
To prove \eqref{rem_conf_claim1} we start with the identities:
\begin{align}
	{}^{I}\Omega^{-1} \td\partial_r({}^{I}\Omega\phi) \ = \ \td\partial_r\phi+\jap{r}^{-2}r\phi 
	\ , \qquad  {}^{I\!I}\Omega^{-1} \td\partial_r({}^{I\!I}\Omega\phi) \ = \ \td\partial_r\phi+2\tau_{+}^{-1}\phi \ . \label{remove_conj_12}
\end{align}
We multiply the first identity by $\tau_{+}$ then square.
Applying Young's inequality with $c \ll 1$ to the resulting estimate,
absorbing the small error term on the LHS,
then adding $|\tau_{+}\tau_{0}\partial_u\phi|^2$ yields:
\begin{align*}
	|\tau_{+}\td\partial_r\phi|^2+|\tau_{+}\tau_{0}\partial_u\phi|^2 \ \lesssim \
	|\tau_{+}{}^{I}\Omega^{-1} \td\partial_r({}^{I}\Omega\phi)|^2
	+|\tau_{+}\tau_{0}{}^{I}\Omega^{-1} \partial_u({}^{I}\Omega\phi)|^2
	+|r^{-1}\tau_{+} \phi|^{2} \ .
\end{align*}
Integrating the last line gives us \eqref{rem_conf_claim1}. Estimate \eqref{rem_conf_claim2} will follow
directly from the claim:
\begin{align}
	 \lp{r^{-1}\tau_{-}\phi(t)}{L^{2}_{x}}^{2}  \ + \ \lp{\phi(t)}{L^{2}_{x}}^{2} \ + \
	\lp{r^{-1}\tau_{+}\phi(t)}{L^{2}_{x}}^{2} \ \lesssim \ RHS \eqref{rem_conf_claim2}
	\label{addthree} \ . 
\end{align}
To prove the claim we bound each of the terms in the LHS above from left-to-right. For the first
term we subtract the two identities in line \eqref{remove_conj_12} then multiply by $\tau_{+}$ to get:
\begin{align}
	\tau_{+}({}^{I\!I}\Omega^{-1} \td\partial_r({}^{I\!I}\Omega\phi)-{}^{I}\Omega^{-1} \td\partial_r({}^{I}\Omega\phi))
	\ &= \ -\jap{r}^{-2}r(C+u)\phi+2\jap{r}^{-2}\phi \label{rweight} \ , 
\end{align}
where we have used $\tau_{+}=C+u+2r$. Next observe that $(2r)^{-1}\leqslant \jap{r}^{-2}r\leqslant r^{-1}$ holds in the set where
$1\leqslant r$. Therefore re-arranging \eqref{rweight}, squaring and using Young's inequality yields:
\begin{align}
	|r^{-1}\tau_{-}\phi|^{2} \ \lesssim \  |\tau_{+}{}^{I\!I}\Omega^{-1} \td\partial_r({}^{I\!I}\Omega\phi)|^{2}
	 \ + \ |\tau_{+}{}^{I}\Omega^{-1} \td\partial_r({}^{I}\Omega\phi)|^{2}  \ + \ |r^{-1}\phi|^{2} \ . 
	\label{conj_rem_interm2}
\end{align}
which is valid for $1\leqslant r$. Multiplying this by a smooth cutoff $\chi_{\jap{r}>2}$, integrating,
and taking the resulting bound in a linear combination with estimate \eqref{hardy_int} then applying the Hardy
bound \eqref{hardy0} gives us $\lp{r^{-1}\tau_{-}\phi(t)}{L^{2}_{x}}^{2} \ \lesssim \ RHS \eqref{rem_conf_claim2}$.
This bounds the first term on LHS \eqref{addthree}.

To control $\lp{\phi(t)}{L^{2}_{x}}^{2}$ we define the linear combination of $\partial_u,\td \partial_r$
derivatives:
\begin{equation*}
	\underline\partial \ = \ 2\partial_u-\td\partial_r \ .
\end{equation*}	
This satisfies:
\begin{equation*}
	\underline\partial(\tau_{+})  =  0  , \qquad \underline\partial(u)  =  2
	, \qquad \underline\partial(r)  =  -1 \  .
\end{equation*}
Collectively these imply:
\begin{multline}
	u({}^{I\!I}\Omega^{-1} \underline\partial({}^{I\!I}\Omega\phi)-
	{}^{I}\Omega^{-1} \underline\partial({}^{I}\Omega\phi)) \ = \
	2 (1-\!\!\jap{u}^{-2})\phi+ ( 1-\!\!\jap{r}^{-2} ) r^{-1}u\phi \\
	 \ =  \ (u\!+\!2r)r^{-1}\phi - (2\jap{u}^{-2} + \jap{r}^{-2}r^{-1}u ) \phi \ .
	\label{conj_rem_id2}
\end{multline}
Rearranging the first equation above, squaring and using Young's inequality
gives us:
\begin{align}
	|\big(1- \jap{u}^{-2}\big)\phi|^{2} \ \lesssim \ u^{2}\big(|
	\frac{\underline\partial({}^{I\!I}\Omega\phi)}{{}^{I\!I}\Omega}|^{2}+|\frac{\underline
	\partial({}^{I}\Omega\phi)}{{}^{I}\Omega}|^{2}\big)+ |( 1-\jap{r}^{-2} ) r^{-1}\tau_{-}\phi|^{2}\label{conj_rem_interm5} \ . 
\end{align}
Next we apply estimate \eqref{conj_rem_interm2} to the last term on the RHS\eqref{conj_rem_interm5},
we multiply the resulting bound by $\chi_{\jap{u}>1}$, then integrate and use estimate \eqref{hardy0}
for the undifferentiated term. This gives us control of $\lp{\phi(t)}{L_x^{2}}$ within the region $\jap{u}>1$.
Taking this resulting bound in a linear combination with estimate \eqref{hardy_near_cone}
then yields $\lp{\phi(t)}{L^{2}_{x}}^{2} \lesssim
RHS \eqref{rem_conf_claim2}$. Finally, the result for
$\lp{r^{-1}\tau_{+}\phi(t)}{L^{2}_{x}}^{2} \lesssim RHS \eqref{rem_conf_claim2}$ follows by using the last line
in equation \eqref{conj_rem_id2} and applying the bounds above.

\end{proof}



\subsection{Core multiplier estimates} 

In this section we list and prove two multiplier bounds which will be
the core constituents of estimates \eqref{unif_bound} and \eqref{k_0_output}.

\begin{prop}[Output of $\partial_u$]\label{du_prop1}
For any interval $[t_0, t_1]\subset I^{*}$ we have the uniform estimate:
 \begin{multline}
		\sup_{t_0\leqslant t\leqslant t_1}\lp{\nabla\phi(t)}{L_{x}^{2}} + 
		\lp{ \td \nabla_x   \phi }{\NLE^{0,-\frac{1}{2}}} 
		  \lesssim   \lp{\chi_{\frac{1}{4}t<r}r^{-\frac{1}{2}-\frac{\delta}{2}}(\nabla\phi, \tau_{0}
		^{-\frac{1}{2}}\td \nabla_{x}\phi)}{L^{2}}  
		\label{NLE_bound} \\
		+  \lp{\frac{\chi_{r<\frac{1}{2}t}}{\tau_{+}^{\frac{\gamma}{2}}\jap{r}
		^{\frac{1}{2}-\frac{\gamma}{2}-\frac{\delta}{2}}}\nabla\phi}{L^{2}}+  \lp{\nabla\phi(t_0)}{L_{x}^{2}}
		+ \sup_j \Big|\dint_{t_0\leqslant t\leqslant t_1} 
		 \Box_g\phi \cdot Y_j \phi \ dV_g \Big|^\frac{1}{2}   \ ,
\end{multline}
where  the vector fields
$Y_j=q_j(u)\partial_u$  are indexed by $j\in\mathbb{Z}$ with 
uniform bounds $|\partial_u^l q_j|\lesssim \tau_-^{-l}$.
\end{prop}

\begin{prop}[Output of  $K_0$]\label{K0_prop}
Let $[t_0, t_1]\subset I^{*}$. For each of the weights
$\Omega ={}^{I}\Omega, {}^{I\!I}\Omega$ have the estimate:
 \begin{multline}
		\sup_{t_0\leqslant t\leqslant t_1} \!\!\lp{\phi(t)}{{}^{\Omega}\CE}
		 +\lp{\phi}{{}^{\Omega}C\!H [t_{0},t_{1}]} \ \lesssim \ \lp{ \chi_{\frac{1}{8}t<r} \frac{
		\tau_{+}}{r^{\frac{1}{2}+\frac{\delta}{2}}}\big(\frac{\tau_{0}\nabla (\Omega\phi)}
		{\Omega},   \frac{\td \nabla_{x}(\Omega\phi)}{\Omega}, \frac{\phi}{r}\big)}{L^{2}}    \\
		+\lp{\chi_{r<\frac{1}{2}t}\frac{\tau_{+}^{1-\frac{\gamma}{2}}}{\jap{r}^{\frac{1}{2}
	 -\frac{\gamma}{2}+\frac{\delta}{2}}}(\nabla\phi,\frac{\phi}{\jap{r}})}{L^{2}} 
	 + \sup_{t_0\leqslant t\leqslant t_1}  \lp{ \tau_{-}^{2}\chi_{r<\frac{1}{2}t}
	  {}^{\Omega}V \big(\frac{\phi}{\Omega}\big)^{2}(t)}{L^{1}_{x}}^{\frac{1}{2}} \\
	+  \lp{\phi(t_0)}{{}^{\Omega}\CE} 
	+ \sup_j \Big|\dint_{t_0\leqslant t\leqslant t_1} 
	\Box_g\phi \cdot  \frac{X_j(\Omega\phi)}{\Omega} \ dV_g \Big|^\frac{1}{2}  , \label{K0_est}
\end{multline}
where $X=q_j(u)K_0$ are indexed by $j\in\mathbb{Z}$,
with $K_0$  given by the formula  $K_0=\tau_{-}^{2} \partial_u + 2(u+r)r\td \partial_r $,  and 
where  $q_j$ has the uniform bounds $|\partial_u^l q_j|\lesssim \tau_-^{-l}$.
\end{prop}

\begin{proof}[Proof of proposition \ref{du_prop1}]
This is a classical multiplier calculation  using formulae  \eqref{pi_hat_formula} and 
\eqref{divergence_formula}. For convenience we employ these in the form of 
Lemma \ref{conf_dividen_lem} with $\Omega\equiv 1$ and $\chi\equiv 0$
in formulae \eqref{div_identity1_cut}--\eqref{AB_formulas3}. We define the
multiplier vector fields:
\begin{equation*}
	Y_j \ = \ (1+\chi_{<j}(u))\partial_u \ ,
\end{equation*}
which are indexed by $j\in\mathbb{Z}$ and
where $\chi_{<j}$ is a non-negative, uniformly bounded,
and monotone decreasing function 
with $\chi_{<j}'=-2^{-j}\chi_j$ supported where $\la u\ra\approx 2^j$.

\step{1}{Output of the ${}^{(Y_j)}A^{\alpha\beta}\partial_\alpha \phi\partial_\beta\phi$ contraction}
We compute using polar Bondi coordinates $(u,r,x^{A})$.
Since $Y^{r}\equiv 0$, identities \eqref{polar_bondi_defiden1}
and \eqref{polar_bondi_defiden2}
get us:
\begin{equation}
		d^\frac{1}{2}{}^{(Y_{j})} \widehat{\pi}^{\alpha\beta} \ = \ 
		2^{-j}\chi_j(u)  ( \eta^{\alpha\beta} - \eta^{\alpha u}\delta_u^\beta
		- \eta^{\beta u}\delta_u^\alpha) + \mathcal{R}^{\alpha\beta} \ , \notag
\end{equation}
where:
\begin{equation}
		\mathcal{R}^{\alpha\beta} \ = \ -\mathcal{L}_{Y_{j}}(d^\frac{1}{2}g^{-1}-\eta^{-1})
		+ 2^{-j}\chi_j(u)  (d^\frac{1}{2}g^{-1}-\eta^{-1}) 
		\ . \notag
\end{equation}
The vector field $\jap{r} \tau_0Y_{j}$ satisfies all of the assumptions \eqref{vect_sym_bnds} and \eqref{good_X_conds}.
So by the results of part \emph{I} of lemma \ref{basic_lie_lem},
$\mathcal{R}^{\alpha\beta}$ satisfies the pointwise bounds:
\begin{align*} 
		|(\mathcal{R}^{rr},r\mathcal{R}^{rA},r^{2}\mathcal{R}^{AB})| 
		\ &\lesssim \ \jap{r}^{-\delta}\tau_{-}^{-1} \tau^{1-\gamma } \ , 
		&|\mathcal{R}^{ru}| \ &\lesssim \ \jap{r}^{-\delta}\tau_{-}^{-1} \tau_0^\frac{1}{2} 
		\tau^{1-\gamma } \ ,  \\
		|r \mathcal{R}^{uA} | \ &\lesssim \
		\jap{r}^{-\delta}\tau_{-}^{-1} \tau_0 \tau^{1-\gamma } \ , 
		&|\mathcal{R}^{uu} | \ &\lesssim \
		\jap{r}^{-\delta}\tau_{-}^{-1} \tau_0^2 \tau^{1-\gamma }
		  \ .
\end{align*}
Thus a little bit of additional calculation involving the previous formulas shows that
$A^{\alpha\beta}$ satisfies the pointwise bounds:
\begin{align*}
		{}^{(Y_j)}A^{rr} \  &\gtrsim \   \tau_-^{-1}\chi_j(u) \!-\! 
		C\jap{r}^{-1-\delta}\tau_0^{-1} \tau^{-\gamma } \ , &|{}^{(Y_j)}A^{rA}| \  
		&\lesssim \  r^{-1}\jap{r}^{-1-\delta}\tau_0^{-1}\tau^{-\gamma }  \ , \\
		r^{2} \!\cdot\! {}^{(Y_j)}A^{AB}  \ &\gtrsim \  \tau_-^{-1}\chi_j(u)\delta^{AB} \!-\! C\jap{r}^{-1-\delta}\tau_0^{-1}
		\tau^{-\gamma } \ ,  & |{}^{(Y_j)}A^{uu}| \ & \lesssim \   \jap{r}^{-1-\delta}\tau^{-\gamma }  \ , \\
		   |{}^{(Y_j)}A^{ur}| \ &\lesssim \   \jap{r}^{-1-\delta}\tau_0^{-\frac{1}{2}}\tau^{-\gamma } 
				 &|r\cdot {}^{(Y_j)}A^{uA}| \ & \lesssim \   \jap{r}^{-1-\delta}\tau^{-\gamma }   , \notag
\end{align*}
where $\delta^{AB}$ denotes the standard inverse metric on $\mathbb{S}^2$. Integrating the
resulting contraction over the time slab $[t_0,t_1]$ and taking $\sup {}_j$ gives us,
with suitable constants $0< c \ll C$:
\begin{multline}
		c\lp{\td \nabla_x \phi}{\NLE^{0,-\frac{1}{2}}}^2  \ \leqslant \ 
		C\lp{\chi_{\frac{1}{4}t<r}r^{-\frac{1}{2}-\frac{\delta}{2}}(\nabla\phi, \tau_{0}
		^{-\frac{1}{2}}\td \nabla_{x}\phi)}{L^{2}}^2 \\
		+C\lp{\chi_{r<\frac{1}{2}t}\tau_{+}^{-\frac{\gamma}{2}}\jap{r}
		^{-\frac{1}{2}+\frac{\gamma}{2}- \frac{\delta}{2}}\nabla\phi}{L^{2}}^2  
		 +  \sup_j \dint_{t_0 \leqslant t\leqslant t_1} {}^{(Y_j)}A^{\alpha\beta}\partial_\alpha  \phi \partial_\beta \phi \
		 dV_{g}  \label{du_output_ext}    \ .
\end{multline}

\step{2}{Output of the boundary terms}
Proposition \eqref{T_timelike} implies that the vector fields
$Y_{j}$ are uniformly timelike for all $j\in\mathbb{Z}$. Thus, 
we have for the boundary terms on line \eqref{div_identity1_cut}:
\begin{equation*}
	\sup_{j}\big(\int_{t=t_0}\!\! {}^{(Y_{j})}P_\alpha N^{\alpha} \
	d^{\frac{1}{2}} dx  -  \int_{t=t_1}\!\! {}^{(Y_{j})}
	P_\alpha N^{\alpha} \ d^{\frac{1}{2}} dx \big)  \lesssim 
	\lp{\nabla\phi(t_0)}{L^2_x}^2  -  \lp{\nabla\phi(t_1)}{L^2_x}^2 \ .
\end{equation*}	
Combining this with \eqref{du_output_ext}, re-arranging terms, taking absolute value,
$\sup$ in $t$ then square roots gives the result.

\end{proof}

\begin{proof}[Proof of Proposition \ref{K0_prop}]
Here we use the conformal multiplier setup of Lemma \ref{conf_dividen_lem}
with weights $\Omega ={}^{I}\Omega, \ {}^{I\!I}\Omega$. For $j\in\mathbb{Z}$
we define the multiplier vector fields:
\begin{equation}
		X_{j} \ = \ (1+\chi_{<j}(u))K_0 \ , \qquad
		K_0 \ = \ \tau_-^2 \partial_u + 2(u+r)r \td \partial_r \ , \notag
\end{equation}
where the $\chi_j$ are the same as in the previous proof. Set $\chi=\chi_{(-\frac{1}{4},\frac{1}{4})}(r/t)$
which smoothly cuts off on $r<\frac{1}{4}t$. Using the divergence identity \eqref{div_identity1_cut} we need to estimate each
term on lines \eqref{td_div_iden_cut3} -- \eqref{AB_formulas3} as well as the
boundary errors on \eqref{boundary1} and \eqref{boundary2}. We do this
for each term separately:

\step{1}{Output of the $A^{\alpha\beta}$ contraction}
As before we compute using polar Bondi coordinates $(u,r,x^{A})$.
The key property of $ K_0$ is that:
\begin{equation*}
	\mathcal{L}_{K_0}\eta^{-1} +(\partial_u  K_0^u 
	+  \td \partial_r K_0^r)\eta^{-1} \equiv 0 \ .
\end{equation*}	
Therefore from 
\eqref{polar_bondi_defiden1} and \eqref{polar_bondi_defiden2} we have:
\begin{multline}
		d^\frac{1}{2}( \widehat{\mathcal{L}_{X_{j}} g} 
		+ 2 X_{j}  \ln \Omega g)^{\alpha\beta}  \  =   \ 2^{-j} \chi_j(u)\big(  \tau_-^2\eta^{\alpha\beta}
		- \eta^{\alpha u}K_{0}^\beta - \eta^{\beta u}K_{0}^\alpha\big)   \\
		 +  2(1  +  \chi_{<j}(u))(r^{-1} K_{0}^{r} -  K_{0} \ln \Omega) )\eta^{\alpha\beta}
		 +   \mathcal{R}^{\alpha\beta}  \ ,  \label{main_partial_u_1/2_line}
\end{multline}  
where,
\begin{equation*}
		\mathcal{R}   =  -\mathcal{L}_{X_{j}}(d^\frac{1}{2}g^{-1}-\eta^{-1}) 
		+ \big[ \tau_-^2 2^{-j}\chi_j(u)
		-(1+\chi_{<j}(u))\big(4(u+r)
		+2(r^{-1} K_0^{r}- K_0 \ln \Omega) \big)\big]
		 (d^\frac{1}{2} g^{-1}-\eta^{-1})  \ .
\end{equation*}
For the term $(r^{-1} K_{0}^{r}- K_{0} \ln \Omega)$ one can compute:
\begin{equation}
	\frac{K_{0}^{r}}{r}- K_{0} \ln {}^{I}\Omega \ = \ \frac{2(u+r)}{\jap{r}^2} \ , 
	\qquad \frac{K_{0}^{r}}{r}- K_{0} \ln {}^{I\!I}\Omega \ 
	= \ \frac{1+C^2}{\tau_{+}}-C \ . \label{K_on_ln_jap}
\end{equation}
These terms are $O(\tau^{-2})$ and $O(1)$, respectively,
this allows us to treat them
as lower order errors below.

To bound $\mathcal{R}$
we use the fact that on a dyadic scale $\tau_+\approx 2^k$ the vector field $2^{-k}K_{0}$ satisfies the 
symbol bounds \eqref{vect_sym_bnds} and all the conditions on line \eqref{good_X_conds} 
except we have $\partial_u X^r\neq 0$. Therefore we are in case \ref{R_alt2} of 
lemma \ref{basic_lie_lem} and so we have  estimates \eqref{R_sym_bnds} with the modification
\eqref{alt_ij_1} and \eqref{alt_ij_2}. In particular thanks to \eqref{mod_coords} the error
$\mathcal{R}$ satisfies:
\begin{align*}
		| \mathcal{R}^{rr} | 
		\ &\lesssim \ \jap{r}^{-\delta}\tau_+\tau_0^{-\frac{1}{2}}\tau^{1-\gamma }  \ , 
		&|(r \ \mathcal{R}^{rA},r^{2}\ \mathcal{R}^{AB})| 
		\ &\lesssim \ \jap{r}^{-\delta}\tau_+ \tau^{1-\gamma } \ , \notag \\ 
		|\mathcal{R}^{ru}| \ &\lesssim \ \jap{r}^{-\delta}\tau_+ \tau_0^\frac{1}{2}\tau^{1-\gamma }   \ , 
		&|r \ \mathcal{R}^{uA} | \ &\lesssim \
		\jap{r}^{-\delta}\tau_+\tau_0\tau^{1-\gamma } \ ,  \notag \\
		| \mathcal{R}^{uu} | \ &\lesssim \
		\jap{r}^{-\delta}\tau_+\tau_0^2\tau^{1-\gamma } 
		  \ . \notag
\end{align*}
Combining the last few lines we get:
\begin{align*}
		|A^{uu}|   \ &\lesssim \ \jap{r}^{-\delta}\tau_+ \tau_0^2\tau^{1-\gamma }  \ , 
		&A^{rr}   \  &\gtrsim \  \tau_+^2 \tau_-^{-1}\chi_j(u)  \! - C\big(\jap{r}^{-\delta}\tau_+\tau_0^{-\frac{1}{2}}
		\tau^{1-\gamma } \big) \ , \\ 
		 r^2A^{AB} \   &\gtrsim  \  \tau_-\chi_j(u) \delta^{AB} -C\big(\jap{r}^{-\delta}\tau_+
		 \tau^{1-\gamma } \big)  \ , 
		& |A^{ur}| \  &\lesssim \   \jap{r}^{-\delta}\tau_+ \tau_0^\frac{1}{2}\tau^{1-\gamma } \ , \\ 
		 |rA^{rA}| \  &\lesssim \   \jap{r}^{-\delta}  \tau_+ \tau^{1-\gamma }   \ ,
		 &|rA^{uA}|  \ &\lesssim \  \jap{r}^{-\delta} \tau_+\tau_0  \tau^{1-\gamma }  \ ,
\end{align*}
where $\delta^{AB}$ again denotes the standard inverse metric on $\mathbb{S}^2$. Integrating the
resulting contraction over the time slab $[t_0,t_1]$ and taking $\sup {}_j$ then gives us, with
suitable constants $0< c \ll C$:

\begin{multline}
	c \lp{\phi}{{}^{\Omega}C\!H [t_{0},t_{1}]}^2  \ \leqslant \ C\lp{\chi_{\frac{1}{4}t<r}
	\frac{\tau_{+}}{r^{\frac{1}{2}+\frac{\delta}{2}}} \big(\frac{\tau_{0}\nabla(\Omega\phi)}{\Omega},
	\frac{\td \nabla_{x}(\Omega\phi)}{\Omega}, 
	\frac{\tau_{0}^{-\frac{1}{2}}  \td \partial_r(\Omega \phi)}{\Omega} \big)}{L^{2}}^{2}    \\
	 + \ C\lp{\chi_{r<\frac{1}{2}t}\frac{\tau_{+}^{1-\frac{\gamma}{2}}}{\jap{r}^{\frac{1}{2}
	 -\frac{\gamma}{2}+\frac{\delta}{2}}}(\nabla\phi,\frac{\phi}{\jap{r}})}{L^{2}}^{2}  
	  \ + \ \sup_j \dint_{t_0 \leqslant t\leqslant t_1} \!\!  \frac{{}^{(X_{j})}A^{\alpha\beta}\partial_\alpha  (\Omega\phi)
	  \partial_\beta (\Omega\phi)}{\Omega^{2}} \
	 dV_{g}     \ ,  \label{spacetime1b}
\end{multline}
where $\Omega ={}^{I}\Omega, \ {}^{I\!I}\Omega$. Note that we used $| \Omega^{-1}\nabla(\Omega\phi)| 
\lesssim |\nabla\phi|+|\jap{r}^{-1}\phi|$ to remove the conjugation in the interior region. Next we observe that the
$\td \partial_{r}$ term in the exterior can be bootstrapped onto the LHS using null energies
since $r^{-\frac{\delta}{4}}\lesssim \epsilon$
and:
\begin{equation*}
	\lp{\chi_{\frac{1}{4}t<r}
	\frac{\tau_{+}}{r^{\frac{1}{2}+\frac{\delta}{4}}} \big(
	\frac{\tau_{0}^{-\frac{1}{2}}  \td \partial_r(\Omega \phi)}{\Omega} \big)}{L^{2}}^{2}
	\ \lesssim \ \lp{\phi}{{}^{\Omega}C\!H [t_{0},t_{1}]}^2 \ .
\end{equation*}



%

\step{2}{Output of the $B^\chi$ term}
Let $\td{\chi}$ be an auxiliary smooth cutoff function supported
inside $r<\frac{1}{2}t $. Multiplying the bounds we derived
for $A$ in the last paragraph times the bounds \eqref{V_bound}
for the conformal potentials ${}^{I}V$, ${}^{I\!I}V$ and dividing by
$({}^{I}\Omega)^2$ and $({}^{I\!I}\Omega)^2$ we get:
\begin{equation}
 		| {}^{I}\! B^\chi| \ \lesssim \ \tau_{+}^{2-\gamma } \jap{r}^{-3-\delta+\gamma}
		 \td{\chi}(r/t)  \ , \qquad \ | {}^{I\!I}\! B^\chi| \ \lesssim \ 
		  \tau_{+}^{2-\gamma}(\tau_{-}\tau_{+})^{-\frac{1}{2}}
		\jap{r}^{ -2-\delta+\gamma }  \td{\chi}(r/t)  \notag \ .
\end{equation}
Integrating these bounds:
\begin{equation}
		 \dint_{t_0\leqslant t\leqslant t_1}  |B^\chi | \phi^2 \ dV_{g}  \  
		 \lesssim \ \lp{\chi_{r<\frac{1}{2}t}\tau_{+}^{1-\frac{\gamma}{2}}\jap{r}
		 ^{-\frac{3}{2}+\frac{\gamma}{2}-\frac{\delta}{2}} \phi}{L^{2}}^{2}  \label{B_K0_output} \ .
\end{equation}


\step{3}{Output of the $C^\chi$ term}
Using estimates \eqref{V_bound} and the fact that $C^{\chi}$ is
supported where $r>\frac{1}{8} t$:

\begin{align*}
	| {}^{I}\! C^\chi| \ \lesssim \  \jap{r}^{-2-\delta} \phi \  \Omega^{-1} K_{0}(\Omega\phi)  \ , \qquad 
	| {}^{I\! I}\! C^\chi| \ \lesssim \ \jap{r}^{-1-\delta} (\tau_{-}\tau_{+})^{-\frac{1}{2}} 
	\phi \ \Omega^{-1} K_{0}(\Omega\phi) \ .
\end{align*}
Integrating these bounds and using Young's inequality
with $0<c\ll 1$:
\begin{align}
		 \dint_{ t_0\leqslant t\leqslant t_1 } \!\!\!\!
		 |C^\chi \phi \cdot \frac{K_0(\Omega \phi)}{\Omega}|\ dV_{g}  \ & \lesssim   \
		   \lp{\chi_{\frac{1}{8}t<r} r^{-\delta} \tau_0^\frac{3}{2}\phi  \partial_u\phi }{L^1}
		  + \lp{\frac{\chi_{\frac{1}{8}t<r}}{r^{\delta} \tau_0^{\frac{1}{2}}} \cdot \phi \frac{\td \partial_r(\Omega \phi)}{\Omega} }{L^1} 
		    \label{C_K0_output}\\
		 \ &\lesssim  \  c^{-1}\lp{\chi_{\frac{1}{8}t<r} r^{-\frac{\delta}{2}}\tau_+^\frac{1}{2}
		(\tau_0 \Omega^{-1}\nabla (\Omega\phi),r^{-1}\phi)}{L^2 }^2 
		\ + \ c \lp{\phi}{{}^{\Omega}C\!H [t_{0},t_{1}]}^2\ . \notag
\end{align}

\step{4}{Output of the boundary term}
Adding estimates \eqref{boundary1} + \eqref{boundary2}:
\begin{equation}
	\frac{{}^{(X_j)}\!\td{P}_\alpha^\chi N^{\alpha}}{\Omega^{2}}
		\  \gtrsim  \  \big(\tau_{-}^2  \! - \!  C K_{0}^r\jap{r}^{-\delta}\tau_0^2\big)
	|\frac{\nabla(\Omega\phi)}{\Omega}|^2  +  
	K_{0}^r | \frac{\td\nabla_x (\Omega \phi)}{\Omega}|^2  
		-   C( |{}^{\Omega}V|\chi (K_{0}^u\! + \! K_{0}^r)   )\big|\frac{\phi}{\Omega}\big|^2\ . 
		\label{boundary_K0_output_0}
\end{equation}
We claim that the following two bounds hold in the asymptotic region $I^{*}$:
\begin{align}
	\tau_{-}^2\ - \ C K_{0}^r\jap{r}^{-\delta}\tau_0^2
	\ \gtrsim \ |\tau_{+}\tau_{0}|^{2} \ , \qquad K_{0}^r
	\ \gtrsim \ 2tr \ . \label{claim_bdry}
\end{align}
To prove this we let $0<c \ll 1$ and treat the regions $r<ct$ and $ct< r$ separately.

\case{1}{Inside $ r< c t$}
Since $\tau_{0}\approx 1$ the first estimate in \eqref{claim_bdry} reduces
to proving $(t-r)^2+C_{0}tr  \gtrsim  |\tau_{+}\tau_{0}|^{2}$ for some new constant $C_{0}>0$.
This clearly holds by choosing $c$ sufficiently small. Since $u=t-r$ here, the second bound
in line \eqref{claim_bdry} holds trivially.

\case{2}{Exterior region $ c t<r$} The first bound in \eqref{claim_bdry} holds
since $K^{r}_{0} \tau_{0}^2 \approx \tau_{-}^{2}$ and the term containing the gain $r^{-\delta}$
is small. The second bound only needs to be proved inside a neighborhood of the wave zone.
Integrating estimate \eqref{grad_u_bnd} yields $K_{0}^r=2tr +O(r^{2-\delta})$ so we can bootstrap the term with
the gain $r^{-\delta}$ here. This finishes the proof of the estimates in line \eqref{claim_bdry}.

We now use this claim to control the boundary terms: integrating \eqref{boundary_K0_output_0} and using line \eqref{claim_bdry} gives us:
\begin{multline}
	\sup_{j}\big(\int_{t=t_0} {}^{(X_j)}\!\td{P}_\alpha^\chi N^{\alpha}
	\Omega^{-2}  \ d^{\frac{1}{2}} dx  -  \int_{t=t_1}
	{}^{(X_j)}\!\td{P}_\alpha^\chi N^{\alpha} \Omega^{-2}  \ d^{\frac{1}{2}} dx \big)   \\
	\ \lesssim \    \lp{\phi(t_0)}{{}^{\Omega}\CE}^2  -  \lp{\phi(t_1)}{{}^{\Omega}\CE}^2
	+ \sup_{t_0 \leqslant t \leqslant t_1} \lp{ \tau_{-}^{2}\chi_{r<\frac{1}{2}t}
	 \Omega^{-2} \  {}^{\Omega}V \phi^{2}(t)} {L^{1}_{x}} \ .  \label{boundary_K0_output} 
\end{multline}

We add equations \eqref{spacetime1b} + \eqref{B_K0_output} + \eqref{C_K0_output}
+ \eqref{boundary_K0_output} and choose $c$ sufficiently small to
bootstrap the ${}^{\Omega}C\!H$ energy term on RHS\eqref{C_K0_output}. Re-arranging terms,
taking absolute values, $\sup$ in $t$ then square roots finishes the proof of
estimate \eqref{K0_est}.

\end{proof}

\begin{rem}
We require the full decay for $g^{ui}-\omega^i\omega_j g^{uj}$
and $g^{uu}$ given by \eqref{mod_coords}
in order to produce estimate \eqref{boundary_K0_output}. If we assume weaker decay for these
components, estimate \eqref{boundary1} would be also be weaker and
the terms containing the gain $K_{0}^{r}r^{-\delta}$ in the exterior could no longer
be bootstrapped. In other words, the exterior proof of \eqref{boundary_K0_output} above would
no longer work.
\end{rem}







\subsection{Proof of the main estimates}
We now prove parts II and III of Theorem \ref{conf_energy_th}.
\begin{proof}[Proof of Proposition \ref{du_prop1}]
We bound each of the error terms on RHS\eqref{NLE_bound}.

\step{1}{Bounding the source term} By Young's inequality:
\begin{align}
	\sup_j \Big|\dint_{t_0\leqslant t\leqslant t_1} 
	 \Box_g\phi \cdot Y_j \phi \ dV_g \Big|^{\frac{1}{2}} \ \lesssim \  c^{-1}\lp{\Box_{g}\phi}{ {\LE}^{*}}
	 \ + \   c\lp{\phi}{\LE} \ . \label{du_output_int2}
\end{align}

\step{2}{Bounding the spacetime error  terms} To control $\nabla\phi$
in the exterior we use $r^{-{\frac{\delta}{4}}}\lesssim \epsilon$
together with $\LE   \subseteq r^{-\frac{1}{2}-\frac{\delta}{4}}L^2 $.
For $ \tau_{0}^{-\frac{1}{2}}\td \nabla_{x}\phi$ in the
exterior we use $r^{-{\frac{\delta}{4}}}\lesssim \epsilon$ together with
$\NLE    \subseteq r^{-\frac{\delta}{4}}\tau^{-\frac{1}{2}}_{-}L^2$. For the interior we use
$\gamma<\delta$ so that $ \LE  \subseteq r^{-\frac{1}{2}+\frac{\gamma}{2}-\frac{\delta}{2}}L^2$
and combine this with $t^{-\frac{\gamma}{4}}\lesssim \epsilon$.
This yields:
\begin{equation}
	\lp{\chi_{\frac{1}{4}t<r}r^{-\frac{1}{2}-\frac{\delta}{2}}(\nabla\phi, \tau_{0}
		^{-\frac{1}{2}}\td \nabla_{x}\phi)}{L^{2}}   \ +\lp{\chi
		_{r<\frac{1}{2}t}\tau_{+}^{-\frac{\gamma}{2}}\jap{r}^{-\frac{1}{2}
		+\frac{\gamma}{2}-\frac{\delta}{2}}\nabla\phi}{L^{2}}
		\ \lesssim \ \epsilon  \big(   \lp{\td \nabla_{x}\phi}
		{\NLE}  \ + \  \lp{\phi}{\LE} \big)\label{du_output_ext2} \ .
\end{equation}
We use this and estimate \eqref{du_output_int2} on the RHS\eqref{NLE_bound}. Adding estimate \eqref{LS1}
to the resulting bound and choosing $c \ll 1$ allows us to bootstrap $c\lp{\phi}{\LE}$ present on RHS\eqref{du_output_int2}.
Bootstrapping all the terms containing $\epsilon$ finishes the proof.


\end{proof}

\begin{proof}[Proof of Proposition \ref{K0_prop}]
We work on each error term on RHS\eqref{K0_est}.

\step{1}{Bounding the source term}
Let $c_1, c\ll 1$ be small positive constants. We bound the
timelike and null/spacelike regions separately. In the exterior
we use Young's inequality together with $\ell_{t}^{1}
\rightarrow\ell_{t}^{1}\times\ell_{t}^{\infty}$ Holder in time to get:
\begin{align}
	\sup {}_j \lp{ \chi_{\frac{1}{8}t<r}  \Box_g\phi\cdot
	\Omega^{-1}X_{j}(\Omega\phi) }{L^1}^{\frac{1}{2}}   \
	&\lesssim \  \lp{\chi_{\frac{1}{8}t<r}\jap{r}^{\frac{1}{2}}
	 \tau_{+}\tau_{0}^{\frac{1}{2}}\Box_{g}\phi \big( \tau_{-}^{\frac{1}{2}}
	\tau_0^{\frac{1}{2}}\partial_u\phi ,  \big|\frac{\tau_{+}}{\tau_0}\big|^{\frac{1}{2}}
	\frac{\td \partial_r(\Omega \phi)}{\Omega} \big)}{L^1}^{\frac{1}{2}}  \label{source_out1}  \\
	&\lesssim  \ c_{1}^{-1}\lp{\Box_{g}\phi}{\ell_{t}^{1}N}\ + \  c_{1}\big( \sup_{t_0\leqslant t\leqslant t_1}
	\lp{\phi(t)}{{}^{\Omega}\CE}+  \lp{\phi}{{}^{\Omega}C\!H [t_{0},t_{1}]}  \big) \ , \notag 
\end{align}
where we have applied the bound:
\begin{equation*}
	\lp{\chi_{\frac{1}{8}t<r} \tau_{-}^{\frac{1}{2}}
	\tau_0^{\frac{1}{2}} \Omega^{-1}\partial_{u}
	(\Omega\phi) }{\ell^\infty_t L^2 }\lesssim
	\sup_{t_0 \leqslant t\leqslant t_1}\lp{\phi(t)}{{}^{\Omega}\CE} \ ,
\end{equation*}
on the last line. In the interior
we use a similar argument in addition to
$| \Omega^{-1}\nabla(\Omega\phi)| \lesssim |\nabla\phi|+|\jap{r}^{-1}\phi|$
to remove the conjugation and get:
\begin{equation*}
	\sup_j \lp{\chi_{r<\frac{1}{2}t}\Box_g\phi \Omega^{-1}
	X_{j}(\Omega\phi)}{L^1}^{\frac{1}{2}}  	\ \lesssim \ 
	c^{-1}\lp{\Box_{g}\phi}{\ell^{1}_tN} \ + \ c \lp{\chi_{r<\frac{1}{2}t}\phi}{\ell^{\infty}_t \LE^1} \ .
\end{equation*}

\step{2}{Bounding the spacetime errors}
In the interior and the exterior we use a similar proof to \eqref{du_output_ext2} 
with the important caveat that we now must split the interior gain $\tau_{+}^{-\frac{\gamma}{2}}$: half goes to produce
$\tau_{+}^{-\frac{\gamma}{4}}\lesssim \epsilon $ while the other half is used to sum the
dyadic pieces in $t$. A bit of additional computation then yields:
\begin{multline} 
	\lp{\chi_{\frac{1}{8}t<r}\frac{ \tau_{+}}{r^{\frac{1}{2}+\frac{\delta}{2}} } \big(
	\frac{\tau_{0}\nabla (\Omega\phi)}{\Omega}, \frac{\td \nabla_{x}(\Omega\phi)}{\Omega}, \frac{\phi}{r}\big)}{L^{2}}
	+\lp{\chi_{r<\frac{1}{2}t} \frac{\tau_{+}^{1-\frac{\gamma}{2}}}{\jap{r}^{\frac{1}{2}
	 -\frac{\gamma}{2}+\frac{\delta}{2}}} (\nabla\phi,\frac{\phi}{\jap{r}})}{L^{2}} 
	\label{K0_sptime2} \\
	 \lesssim \epsilon \big(\sup_{t_0\leqslant t\leqslant t_1} \lp{\phi(t)}{ \CE}+\sup_{t_0\leqslant t\leqslant t_1} \lp{\phi(t)}{{}^{\Omega}\CE}
	+ \lp{\chi_{r<\frac{1}{2}t}\phi}{\ell^{\infty}_t \LE^1} \big) \ . 
\end{multline}

\step{3}{Bounding the undifferentiated boundary errors}

We apply the bounds \eqref{error_density_int1}
and \eqref{error_density_int2} for each undifferentiated
boundary error. Next we choose $c_1 \ll 1$ small in step 1 above
so we can bootstrap the small errors in \eqref{source_out1}. This,
together with the result of step 2 above yields the following
two bounds for the conjugated energies for ${}^{I}\Omega$ and
${}^{I\!I}\Omega$:
\begin{multline}
		\sup_{t_0\leqslant t\leqslant t_1} \lp{\phi(t)}{{}^{I}\CE} \ + \ 
		\lp{\phi}{{}^{I}\CH} \ \lesssim \   \sup_{t_0\leqslant t\leqslant t_1} 
		\lp{\phi(t)}{{}^{I\!I}\CE} +  \epsilon \sup_{t_0\leqslant t\leqslant t_1}\lp{\phi(t)}{\CE} \\ 
		\label{conf_int_w1}
		+\lp{\phi(t_0)}{{}^{I}\CE}+ c^{-1}\lp{\Box_{g}\phi}{\ell^{1}_tN }
		+c \lp{\chi_{r<\frac{1}{2}t}\phi}{\ell^{\infty}_t \LE^1}  \ ,
\end{multline}
\begin{multline}
		\sup_{t_0\leqslant t\leqslant t_1} \lp{\phi(t)}{{}^{I\!I}\CE} \ + \ 
		\lp{\phi}{{}^{I\!I}\CH}
		 \ \lesssim  \ \epsilon \sup_{t_0\leqslant t\leqslant t_1}\lp{\phi(t)}{\CE}
		 \label{conf_int_w2} \\
		  +\lp{\phi(t_0)}{{}^{I\! I}\CE}+ c^{-1}\lp{\Box_{g}\phi}{\ell^{1}_t N}
		+c  \lp{\chi_{r<\frac{1}{2}t}\phi}{\ell^{\infty}_t \LE^1}  \ . 
	\end{multline}
Note that $\lp{\phi(t)}{{}^{I\!I}\CE}$ on the RHS\eqref{conf_int_w1} does not include any smallness
constants. However, taking an appropriate linear combination of estimates \eqref{conf_int_w1},  \eqref{conf_int_w2}
and \eqref{unif_bound} allows us to absorb this error. Since all the fixed-time energies in the resulting estimate
are of the form:
\begin{equation*}
	\lp{\phi(t)}{{}^{I}\CE}  + \lp{\phi(t)}{{}^{I\!I}\CE}+\lp{\nabla\phi(t)}{L_{x}^{2}} \ ,
\end{equation*}
we apply estimate \eqref{conjlemma} to trade this for $\lp{\phi(t)}{\CE}$.
Bootstrapping the remaining term $\epsilon\cdot \sup_{t_0\leqslant t\leqslant t_1}\lp{\phi(t)}{\CE}$ finishes the proof.
\end{proof}


\section{Commutators} \label{comm_sect}

In this section we commute the equation once with
the Lie algebra $\mathbb{L}=\{ \partial_{u},   S ,\Omega_{ij}\}$
and use our conformal energy estimate \eqref{conf_energy_est} to produce the higher order bound
\eqref{conf_energy_est_vf}. By Lemma \ref{reduction} it suffices to prove our estimate for the asymptotic region
$t\in I^{*}$. Inside $I^{*}$ the corresponding estimate \eqref{conf_energy_est_vf} will follow from:


\begin{thm}\label{conf_energy_vf_th}
Assume the hypotheses of the main theorem and
estimate \eqref{conf_energy_est} hold. Then for all
$[t_0, t_1]\subset I^{*}$ we have the following uniform bounds:
\begin{enumerate}[I)]
	\item(CE estimate for $\partial_{u}\phi$ with error)
		\begin{multline}
			\sup_{t_0\leqslant t\leqslant t_1}  \lp{\partial_{u}\phi(t)}{\CE}
			+\lp{\partial_{u}\phi}{S} \   \lesssim  \ \epsilon
			\big( \sup_{t_0\leqslant t\leqslant t_1} \lp{\phi(t)}{\CE_{1}}
			+ \lp{\chi_{r<\frac{1}{2}t}\phi}{\ell_{t}^{\infty}\LE^{1}}  \big)\\
			+\lp{\partial_{u}\phi(t_0)}{\CE}+\lp{\Box_{g}\phi}{\ell_{t}^{1}N_{1}} 
			\ . \label{en_plus_ce}
		\end{multline}
	\item(CE estimate for $S\phi$ and $\Omega_{ij}\phi$ with error)
		\begin{multline}
			\sup_{t_0\leqslant t\leqslant t_1}\sum_{\Gamma= S,\Omega_{ij}}
			\!\!\!\!\lp{\Gamma\phi(t)}{\CE}+\lp{\Gamma\phi}{S} 
			\ \lesssim  \ \epsilon \big(\sup_{t_0\leqslant t\leqslant t_1}
			\!\!\lp{\phi(t)}{\CE_{1}}+ \lp{\chi_{r<\frac{1}{2}t}\phi}
			{\ell_{t}^{\infty}\LE_{1}^{1}}\big)\\
			+\sum_{\Gamma= S,\Omega_{ij}}\lp{ \Gamma\phi(t_0)}{\CE}
			+\lp{\Box_{g}\phi}{\ell_{t}^{1}N_{1}} \ . \label{conf_loss}
		\end{multline}
\end{enumerate}
\end{thm}
Let's show how the conformal energy estimate with vector fields follows
from this theorem:

\begin{proof}[Proof of estimate \eqref{conf_energy_est_vf}]
Adding estimates \eqref{en_plus_ce} + \eqref{conf_loss} + \eqref{conf_energy_est} yields:
	\begin{equation*}
			\sup_{t_0\leqslant t\leqslant t_1} \lp{\phi(t)}{\CE_1}
			 +  \lp{\phi}{S_{1}}   
			\ \lesssim  \
			\epsilon \big(\sup_{t_0\leqslant t\leqslant t_1}
			\lp{\phi(t)}{\CE_{1}}+ \lp{\chi_{r<\frac{1}{2}t}\phi}
			{\ell_{t}^{\infty}\LE_{1}^{1}}\big)
			+ \lp{ \phi(t_0)}{\CE_1}  +  \lp{\Box_{g}\phi}{\ell_{t}^{1}N_{1}} \ .
	\end{equation*}
We bootstrap the small error terms and close the estimate. Inside $[0,T^{*}]$, lemma
\eqref{reduction} gives us the analogous bound. A combination of these two estimates
then yields \eqref{conf_energy_est_vf}.
\end{proof}

The rest of this section is devoted to proving Theorem \ref{conf_energy_vf_th}.





\subsection{Preliminary estimates}


Our goal in this first part is to establish
some commutator bounds for all vector
fields $\Gamma\in \mathbb{L}$.
We start with the following:

\begin{lem}{(Pointwise bounds)}\label{comm_point} Let $\Gamma\in \mathbb{L}$ be in rectangular 
Bondi coordinates $(u,x^i)$. The following uniform estimate holds:
	\begin{multline}
		\sum_{\Gamma= S,\Omega_{ij}} |[\Box_{g},\Gamma]\phi| \  \lesssim \  \jap{r}^{-2-\delta}
		\tau^{1-\gamma } \big(\!\! \sum_{\substack{1\leqslant {k}+|J|\leqslant 2\\ k\neq 2,|J| \neq2} }
		|\tau_0^{-\frac{1}{2}} (\jap{r} \tau_0\partial_u)^k(\jap{r}\td \partial_{x})^J  \phi| 
		+   | (\jap{r}\tau_{0}\partial_{u})^{2}\phi|  \\
		 +  | (\jap{r}\td \partial_{x})^{2}\phi|  +  \jap{r}^{2} | \Box_{g}\phi| \big)
		 \label{Q_bondi_est2}  \ . 
	\end{multline}
In the case of $[\Box_{g},\partial_{u}]$ the same estimate holds
with the exponent $1-\gamma$ above replaced with $-\gamma$.
\end{lem}

\begin{proof}
By \eqref{X_comm_form}, \eqref{usual_R}, and
\eqref{S_comm_form} it suffices to bound the quantities:
\begin{align*}
	\mathcal{E}(u,\alpha)   &= d^{-\frac{1}{2}} \partial_u d^{\frac{1}{2}}\mathcal{R}_{\Gamma}^{u\alpha}\partial_{\alpha}\phi \  , 
	&\mathcal {E}(\alpha, u)     &= d^{-\frac{1}{2}} \partial_{\alpha}   d^{\frac{1}{2}} \mathcal{R}_{\Gamma}^{u\alpha}\partial_{u}\phi \ , \\
	\mathcal { E}(i,j)   &= d^{-\frac{1}{2}}  \td \partial_{i} d^{\frac{1}{2}}\mathcal{R}_{\Gamma}^{ij}\td\partial_{j}\phi \ ,  
	&\mathcal{E}(s)   &=  d^{-\frac{1}{2}}\Gamma(d^{\frac{1}{2}}) \ ,
\end{align*}
where:
\begin{align*}
		d^{\frac{1}{2}}\mathcal{R}_{\Omega_{ij}}  &= - \mathcal{L}_\Gamma
		(d^\frac{1}{2}g^{-1}-\eta^{-1}) \ , \\
		d^{\frac{1}{2}}\mathcal{R}_{S}   &=
		- \mathcal{L}_S(d^\frac{1}{2}g^{-1}-\eta^{-1}) 
		-(\partial_u S^u + \td \partial_r S^r + \td \partial_i\overline S^i )
		(d^\frac{1}{2} g^{-1}-\eta^{-1}) \ .
\end{align*}
To bound these quantities we use part I of lemma
\ref{basic_lie_lem} together with estimate \eqref{det_bondi_decay} for
$d=\det(g_{\alpha\beta})$. This yields:
 \begin{align*}
 	\big| \mathcal{E}(u,\alpha) |  +    \big | \mathcal {E}(\alpha, u) \big| 
	 \ &  \lesssim \  \jap{r}^{-\delta} \tau^{1-\gamma }  
	\big( \tau_{0}^{\frac{1}{2}} |\td \partial_{x} \partial_{u}\phi| + 
	\tau_{0}^{2} |\partial_{u}^{2}\phi|  + \jap{r}^{-1}\tau_{0}^{\frac{1}{2}}|\partial_{u}\phi|
	+  \jap{r}^{-1}\tau_{0}^{-\frac{1}{2}} |\td \partial_{x}\phi|\big) \ , \\
	\big | \mathcal { E}(i,j) \big| \  &\lesssim \ \jap{r}^{-\delta}\tau
	^{1-\gamma } \big(|\td \partial_{x}^{2} \phi|   +    \jap{r}^{-1}|\td \partial_{x}\phi |\big) \ , \\
	\big | \mathcal{E}(s)  \big|  \ &\lesssim \  \jap{r}^{-\delta}\tau^{1-\gamma } |\Box_{g}\phi | \ .
 \end{align*}
In the case of $\Gamma=\partial_{u}$ we use part III of lemma \ref{basic_lie_lem}
which yields the improvement in the interior by the same argument.
\end{proof}


\begin{rem}
A quick computation using \eqref{mod_coords} shows that the estimate derived in
lemma \ref{comm_point} above for the vector field $\partial_{u}$ is missing a
weight $\tau_{-}^{-1}$ in the exterior. We chose to omit this weight since
there's no real improvement in the estimates below if we were to include it.
Ultimately, this is a consequence of our weak decay for the $g^{ij}$ metric
components which, essentially, force us to treat $\partial_u$
at the same level of decay as $S$ and $\Omega_{ij}$ in the exterior.
\end{rem}

\begin{rem}
From this point on we will often make use of the parameter $\mu\ll 1$ satisfying the smallness condition
in definition \ref{bspara}.a.
\end{rem}

In the region $\{ \mu t< r\} \cap I^{*}$ we control all weighted combinations
of two derivatives in RHS\eqref{Q_bondi_est2}  via:

\begin{lem}[Exterior Klainerman-Sideris identity]
Let $\Gamma\in \mathbb{L}$ be in rectangular 
Bondi coordinates $(u,x^i)$. In the exterior region
$\{ \mu t < r \} \cap I^*$ one has the following uniform estimate:
\begin{align}
		&\sum_{1\leqslant k+|J|\leqslant 2} | (r\tau_0 \partial_u)^k(r\td\partial_x)^J \phi|
		\  \lesssim  \!\!\!\sum_{ {k}+|J|=1, \ |I|\leqslant 1}\!\!\!\!
		| (r \tau_0\partial_u)^k(r\td \partial_x)^J \Gamma^I \phi| + r^2\tau_0 |\Box_g\phi |
		\ .  \label{KS_iden} 
\end{align}
\end{lem}

\begin{proof}
Let $\mathcal{R}=d^\frac{1}{2}g-\eta$. Inside $\{ \mu t< r\} \cap I^{*}$
the estimates \eqref{mod_coords} imply the uniform bound:
\begin{equation}
		r^2\tau_0 | \partial_\alpha  \mathcal{R}^{\alpha\beta} \partial_\beta \phi|
		\ \lesssim \ r^{-\delta} 
		\!\!\!\! \sum_{1\leqslant k+|J|\leqslant 2} | (r \tau_0\partial_u)^k(r\td \partial_{x})^J \phi|
		 \ , \notag
\end{equation}
where all quantities are computed with respect to Bondi coordinates.
By the support property and estimate \eqref{tstar_cond} we have $r^{-\delta} \lesssim
(\mu t)^{-\delta}\lesssim \epsilon$. Thus, we may replace $\Box_g$ by the
Minkowski wave operator $\Box_\eta$
in estimate \eqref{KS_iden}. Next, we have the two identities:
\begin{align}
		\frac{r^2}{\tau_+} \td \partial_{r} S \ &= \ \frac{r^2 u}{\tau_+}\partial_u\td \partial_{r} + r^3\tau_+^{-1}
		\td \partial_{r}^2 +\frac{r^2}{\tau_+}\td \partial_{r} \ , \notag\\
		\frac{1}{2} \frac{r^2u}{\tau_+}\Box_\eta \ &= \ -\frac{r^2u}{\tau_+}\partial_u\td \partial_{r}
		+ \frac{1}{2} \frac{r^2u}{\tau_+}\td \partial_{r}^2 -\frac{ru}{\tau_+}\partial_u + \frac{ru}{\tau_+}\td \partial_{r}
		+ \frac{1}{2}\frac{u}{\tau_+} \sum_{i<j}(\Omega_{ij})^2
		\ , \notag
\end{align}
where we used equation \eqref{eta_idens} on the second line. Adding the two
operators on the LHS above applied to $\phi$ yields:
\begin{equation}
		|(r\td \partial_{r})^2\phi| \ \lesssim  \sum_{ {k}+|J|=1, \ |I|\leqslant 1}
		| (r \tau_0\partial_u)^k(r\td \partial_{x})^J \Gamma^I \phi| 
		+  r^2 \tau_0|\Box_\eta \phi | \ . \notag
\end{equation}
In addition to this we have the inequality:
\begin{equation*}
		|(r\tau_0 \partial_u)^2\phi| + |(r\tau_0 \partial_u)(r\td \partial_{r})\phi| 
		\ \lesssim \ |(r\td \partial_{r})^2\phi| +\!\!\!
		\sum_{ {k}+|J|=1 } | (r\tau_0\partial_u)^k(r\td \partial_{x})^J  S \phi| 
		+\sum_{ {k}+|J|=1 } | (r\tau_0\partial_u)^k(r\td \partial_{x})^J \partial_u \phi| \ . 
\end{equation*}
All other combinations of derivatives on LHS\eqref{KS_iden} are  controlled
by the sum on RHS\eqref{KS_iden}.
\end{proof}




As a consequence of these two lemmas we get:

\begin{prop}[Global commutator bounds] \label{global_commprop}
Let $\Gamma= S ,\Omega_{ij}$. The following uniform estimates hold on $I^*$.
\begin{equation}
		\lp{ [\Box_g, \Gamma]\phi }{\ell_{t}^{1}N} 
		\   \lesssim \   \lp{\chi_{r< \mu t} \jap{r}^{-\frac{1}{2}}\tau_{+}^{2-\frac{\gamma}{2}}(\nabla^{2}\phi,
		  \jap{r}^{-1}\nabla\phi)}{L^{2}} 
		 + \epsilon  \sup_{t_{0}\leqslant t\leqslant t_{1}}\lp{\phi(t)}{\CE_1}
		 +\lp{\Box_{g}\phi}{\ell_{t}^{1}N_{1}}   \label{comm_est_SR_int} \ ,
\end{equation}
we also have the following (interior) improvement in the case of
$\Gamma=\partial_{u}$:
\begin{multline}
		\lp{[\Box_g,\partial_{u}]\phi}{\ell_{t}^{1}N} \ \lesssim \
		\lp{ \chi_{r< \mu t} \jap{r}^{-\frac{1}{2}}\tau_{+}^{1-\frac{\gamma}{2}}(\nabla^{2}\phi
		,  \jap{r}^{-1}\nabla\phi)}{\ell^{\infty}_{r}L^{2}}\\
		  +   \epsilon  \sup_{t_{0}\leqslant t\leqslant t_{1}}\lp{\phi(t)}{\CE_1}
		 +\lp{ \jap{r}^{-1-\delta}\tau_{0}^{-1}\tau^{-\gamma}\Box_{g}\phi)}{\ell_{t}^{1}N} \ . \label{comm_est_der_ce} 
\end{multline}
\end{prop}

\begin{proof}[Proof of estimate \eqref{comm_est_SR_int}]
We treat the regions $r<\mu t$ and $\mu t< r$ separately.

\case{1}{The interior $ r< \mu t$}
We multiply estimate \eqref{Q_bondi_est2} by:
\begin{equation*}
	2^{\frac{1}{2}j}2^{k}\chi_{t\approx 2^{k}}\chi_{\jap{r}\approx 2^{j}}
	\chi_{r< \mu t} \  ,
\end{equation*}
square, integrate then use $\gamma<\delta$ together with the inclusion
$L^2 \subseteq \jap{r}^{\gamma-\delta}\! \tau_{+}^{-\frac{\gamma}{2}}\ell_{t,r}^{1}L^2 $. For the source term we
apply estimate \eqref{hardy0} with $a=1/2$. Combining all this:
\begin{equation}
		\sum_{\Gamma= S,\Omega_{ij}}\lp{ \chi_{r< \mu t}[\Box_g, \Gamma]\phi }{\ell_{t}^{1}N}
		\ \lesssim \ \lp{\chi_{r< \mu t} \jap{r}^{-\frac{1}{2}}\tau_{+}^{2-\frac{\gamma}{2}}(\nabla^{2}\phi,
		  \jap{r}^{-1}\nabla\phi)}{L^{2}}  
		  + \lp{\chi_{r< \mu t}\Box_{g}\phi}{\ell_{t}^{1}N_{1}}\ . \label{dt_int_ext1}
\end{equation} 
\case{2}{The exterior region $\mu t< r$} We multiply estimate 
\eqref{Q_bondi_est2} by:
\begin{equation*}
	\jap{r}^{\frac{1}{2}}\tau_{+}\tau_{0}^{\frac{1}{2}}
	\chi_{t\approx 2^{k}}\chi_{\jap{r}\approx 2^{j}}
	\chi_{\jap{u}\approx 2^{l}}\chi_{\mu t< r} \ ,
\end{equation*}	
then apply the Klainerman-Sideris estimate \eqref{KS_iden}. Squaring, integrating and
using the support property together with $L^2 \subseteq \jap{r}^{-\frac{\delta}{2}}\ell_{t,r}^{1}L^2$ yields:
\begin{align*}
	&\sum_{\Gamma= S,\Omega_{ij}} \lp{\chi_{\mu t< r} [\Box_{g},\Gamma]\phi}{\ell_{t}^{1}N}\\
	&\qquad \lesssim  \sum_{ {k}+|J|=1 \ ,  |I|\leqslant 1}
	\lp{\chi_{\mu t< r}\jap{r}^{\frac{1}{2}}\tau_{+}\tau_{0}^{\frac{1}{2}}
	\jap{r}^{-2-\frac{\delta}{2}}\tau_0^{-\frac{1}{2}}( (r \tau_0\partial_u)^k
	(r\td \partial_x)^J \Gamma^I \phi )}{L^{2}} + \lp{   \chi_{\mu t< r} \jap{r}^{-\frac{\delta}{2}} \Box_{g}\phi}{\ell_{t}^{1}N_{1}}\\
	&\qquad   \lesssim   \
	 \lp{ \chi_{\mu t< r} \jap{r}^{-\frac{1}{2}-\frac{\delta}{2}}\phi(t)}{L_{t}^{2}(\CE_{1})}
	 +\lp{\Box_{g}\phi}{\ell_{t}^{1}N_{1}} \ .
\end{align*}
To bound the $L_{t}^{2}(\CE_{1})$ term above we apply the support property then split the resulting
$t^{-\frac{\delta}{2}}$: half goes to produce $(\mu t)^{-\frac{\delta}{4}}\lesssim\epsilon$
thanks to \eqref{tstar_cond}, while the rest is used to integrate in time
since $t^{-\frac{1}{2}-\frac{\delta}{4}} \in L_{t}^{2}$. Taking $\sup_t$
and integrating the $t$ weight then finishes the proof of \eqref{comm_est_SR_int}.
\end{proof}

\begin{proof}[Proof of estimate  \eqref{comm_est_der_ce}]
This follows the same exact proof as the previous estimate. The improvement in
the interior is a direct result of using the version of estimate \eqref{Q_bondi_est2}
with the better interior weight  $\tau^{-\gamma}$ corresponding to
$\Gamma=\partial_{u}$, as well as an application of the inclusion
$\ell^{\infty}_{r}L^2\subseteq \jap{r}^{\gamma-\delta}\tau_{+}^{-\frac{\gamma}{2}}\ell_{t,r}^{1}L^{2}$ for the
interior.

\end{proof}



\subsection{Interior $L^{2}$-estimates for two derivatives}

In order to establish our main estimates we must bound
the error terms supported inside
$r< \mu t$ on RHS\eqref{comm_est_SR_int} and RHS\eqref{comm_est_der_ce}. This
will be achieved in the next lemma by establishing some
weighted $L^2$ bounds for $\lp{\chi_{r<  \mu t}\jap{r}^{\frac{1}{2}}
\tau_{+}^{a}\nabla^{2}\phi}{L^2}$ with the exponent ``$a$" depending on
the vector field $\Gamma$. In the case $\Gamma=\partial_{u}$,
we have $a=1-\gamma/2$ and it will suffice to commute
with cutoffs and use some standard elliptic estimates.

For the cases $\Gamma=S,\Omega_{ij}$, the $\tau_{+}^a$ weight comes with an
exponent of $a=2-\gamma/2$, a number which no longer allows us to commute
with cutoffs. To resolve this issue, we will use a Klainerman-Sideris type bound
together with the elliptic estimate \eqref{elliptic} to establish the necessary inequalities
in this case. 
In light of this discussion, the ellipticity of
the operator $P(t,x,\nabla_{x})$ (see Corollary \ref{P_cor}.II)
is of fundamental importance for the proof of the next lemma.
Consequently, all computations are in $(t,x^i)$ coordinates here.





\begin{lem}[Klainerman-Sideris type estimates for the interior] \label{lemma_two_der}
The following uniform bounds hold inside the region
$I^{*}$:

\begin{enumerate}[I)]
\item(Bounds for $\nabla\partial_{t}\phi$)
	\begin{align}
		\lp{  \frac{\chi_{r< \mu t}\tau_{+}^{1-\frac{\gamma}{2}}}
		{\jap{r}^{\frac{1}{2}}}\nabla\partial_{t}\phi}{\ell^{\infty}_{r}L^{2}}
		 \ &\lesssim  \  \epsilon  \lp{\chi_{r< \frac{1}{2} t}\partial_{t}\phi}
		{\ell_{t}^{\infty}\LE^{1}} \label{two_dt_round1} \ ,  \\
		\lp{ \frac{\chi_{r< \mu t} \tau_{+}^{2-\frac{\gamma}{2}}}
		{\jap{r}^{\frac{1}{2}-\frac{\gamma}{10}}}  \nabla\partial_{t}\phi}{ L^{2}}
		\  &\lesssim \ \epsilon  \lp{\chi_{r< \frac{1}{2} t}\phi}{\ell_{t}^{\infty}\LE_{1}^{1}}
		+\mu\lp{\chi_{r<  \mu t}\frac{\tau_{+}^{2-\frac{\gamma}{2}}}
		{\jap{r}^{\frac{1}{2}-\frac{\gamma}{10}}}\nabla_{x}^{2}\phi}{ L^{2}}  \label{two_dt_round2} \ .
	\end{align}
\item(Bounds for $\nabla^{2}_{x}\phi$)
	\begin{align}
		\lp{\frac{ \chi_{r< \mu t} \tau_{+}^{1-\frac{\gamma}{2}}}
		{ \jap{r}^{\frac{1}{2}}}\nabla_{x}^{2}\phi}{\ell^{\infty}_{r}L^{2}}
		\ &  \lesssim \  \epsilon \lp{\chi_{r<\frac{1}{2}t}(\partial_{t}\phi,\phi)}
		{\ell_{t}^{\infty}\LE^{1}}
		+\lp{\Box_{g}\phi}{N_{1}} 
		\label{two_dx_round1} , \\
		\lp{\frac{\chi_{r< \mu t}\tau_{+}^{2-\frac{\gamma}{2}}}
		{\jap{r}^{\frac{1}{2}-\frac{\gamma}{10}}}  \nabla_{x}^{2}\phi}{L^{2}}
		\ &  \lesssim \ \epsilon \big(\sup_{t_0\leqslant t\leqslant t_1} \lp{\phi(t)}{\CE_1}
			 +  \lp{\phi}{S_{1}}\big)+\lp{\Box_{g}\phi}{N_{1}} 
		 +\lp{\frac{\chi_{r< \mu t}\tau_{+}^{2-\frac{\gamma}{2}}}
		 {\jap{r}^{\frac{1}{2}-\frac{\gamma}{10}}}\nabla\partial_{t}\phi}{L^{2}}  	
		 \ .  \label{two_dx_round2}
	\end{align}
\item(Bounds for  $\nabla^2\phi$)
	\begin{align}
		\lp{ \frac{ \chi_{r< \mu t} \tau_{+}^{1-\frac{\gamma}{2}}}
		{ \jap{r}^{\frac{1}{2}}}\nabla^{2}\phi}{\ell^{\infty}_{r}L^{2}}
		\ &\lesssim \ \epsilon \lp{\chi_{r<\frac{1}{2}t}(\partial_{t}\phi,\phi)}
		{\ell_{t}^{\infty}\LE^{1}}
		  , \label{two_dx_tless_energy} \\
		\lp{\frac{\chi_{r< \mu t}\tau_{+}^{2-\frac{\gamma}{2}}}
		{\jap{r}^{\frac{1}{2}-\frac{\gamma}{10}}}  \nabla^{2}\phi}{L^{2}}
		\  & \lesssim \  \epsilon \big(\!\!\sup_{t_0\leqslant t\leqslant t_1} \lp{\phi(t)}{\CE_1}
			 +  \lp{\phi}{S_{1}}\big)
		+\lp{\Box_{g}\phi}{N_{1}}  \label{twodersclaim10}
	\end{align}
\end{enumerate}
\end{lem}

\begin{proof}[Proof of Lemma \ref{lemma_two_der}]
\Part{1}{Bounds for $\nabla\partial_{t}\phi$}
For estimate \eqref{two_dt_round1} we split the gain $t^{-\frac{\gamma}{2}}$:
half goes to $t^{-\frac{\gamma}{4}}\lesssim\epsilon$ and the other half
is used in the inclusion $\ell_{t,r}^{\infty}L^2 \subseteq    t^{-\frac{\gamma}{4}}\ell_{r}^{\infty}L^{2}   $.
The bound follows by using the definition of the norms and the support property.

For estimate \eqref{two_dt_round2} we start with the pointwise inequality:
\begin{align}
	 |\partial_t\phi | \ \lesssim \  |t^{-1} S\phi| \ + \ (r/t)| \partial_{r}\phi  |  \ , \label{Stodt}
\end{align}
which is valid inside $r<\frac{1}{2}t$. Applying this to
$\partial_{t}\phi$ and $\nabla_{x}\phi$ yields, respectively:
	\begin{align*}
		| \partial_t^{2}\phi| \ &\lesssim \  |t^{-1}\partial_{t} S\phi
		|+|t^{-1}\partial_{t}\phi|+(r/t)| \partial_{r}\partial_{t}\phi  |
		\ , \qquad 
		| \partial_t\nabla_{x}\phi| \ \lesssim \  |t^{-1}\nabla_{x} S\phi|
		+|t^{-1}\nabla_{x}\phi| +(r/t)| \partial_{r}\nabla_{x}\phi  | \ .
	\end{align*}
Multiplying these two bounds by $\jap{r}^{-\frac{1}{2}+\frac{\gamma}{10}}\tau_{+}^{2-\frac{\gamma}{2}}
\chi_{r<\mu t}$, squaring, integrating and using the support property $(r/t)<\mu$ we get:
	\begin{align}
		\lp{\frac{\chi_{r< \mu t} \tau_{+}^{2-\frac{\gamma}{2}}}{ \jap{r}^{\frac{1}{2}-\frac{\gamma}{10}}}
		\partial^{2}_{t}\phi}{L^{2}}  \ &\lesssim \ \epsilon  \lp{\chi_{r< \frac{1}{2} t}\phi}
		 {\ell_{t}^{\infty}\LE_{1}^{1}} + \mu\lp{\frac{\chi_{r< \mu t} \tau_{+}^{2-\frac{\gamma}{2}}}
		 {\jap{r}^{\frac{1}{2}-\frac{\gamma}{10}}}\nabla_{x}\partial_{t}\phi}{L^{2}} \ , \label{two_dtdt}\\
		\lp{\frac{\chi_{r< \mu t} \tau_{+}^{2-\frac{\gamma}{2}}}{ \jap{r}^{\frac{1}{2}-\frac{\gamma}{10}}}
		\nabla_{x}\partial_{t}\phi}{L^{2}} \ &\lesssim \ \epsilon \ \lp{\chi_{r< \frac{1}{2} t}\phi}
		 {\ell_{t}^{\infty}\LE_{1}^{1}}+ \mu\lp{\frac{\chi_{r< \mu t} \tau_{+}^{2-\frac{\gamma}{2}}}{ \jap{r}^{\frac{1}{2}-\frac{\gamma}{10}}}\nabla_{x}^{2}\phi}{L^{2}} \ , \label{two_dxdt}
	\end{align}
where we have split the gain $t^{\frac{\gamma}{10}-\frac{\gamma}{2}}$ as in the previous proof
to produce the $\epsilon\cdot  \ell_{t}^{\infty}\LE_{1}^{1}$ terms. By definition \ref{bspara}.a, the constant
$\mu$ is small enough that we can add \eqref{two_dtdt} + \eqref{two_dxdt}, bootstrap the term $\nabla_{x}\partial_t\phi$,
and get estimate \eqref{two_dt_round2}.

\Part{2}{Bounds for $\nabla^{2}_{x}\phi$}
To prove estimate \eqref{two_dx_round1} we let $P(t,x,\nabla_{x}) $
be as in Corollary \ref{P_cor}.II. By that result, the operator $P$ is uniformly
elliptic in $(t,x^i)$ coordinates. Thus, commuting $\chi_{r< \mu t}$ with $\nabla_{x}^{2}$
and using basic elliptic estimates:
	\begin{align*}
	\lp{ \chi_{r< \mu t}  \jap{r}^{-\frac{1}{2}}\tau_{+}^{1-\frac{\gamma}{2}}\nabla_{x}^{2}\phi}{\ell^{\infty}_{r}L^{2}}
	  \  &\lesssim   \lp{ \chi_{r< \mu t}  \jap{r}^{-\frac{1}{2}}\tau_{+}^{1-\frac{\gamma}{2}} P\phi}{\ell^{\infty}_{r}L^{2}}
	+\mu^{-1}\lp{ \chi_{r\sim \mu t}  \jap{r}^{-\frac{1}{2}}\tau_{+}^{1-\frac{\gamma}{2}}(r^{-1}\nabla\phi, \ r^{-2}\phi)}{\ell^{\infty}_{r}L^{2}} \\
	 & \lesssim \  \lp{ \chi_{r< \mu t}  \jap{r}^{-\frac{1}{2}}\tau_{+}^{1-\frac{\gamma}{2}} P\phi}{\ell^{\infty}_{r}L^{2}}
	+ \epsilon\lp{\chi_{r< \frac{1}{2}t} \phi}{\ell_{t}^{\infty}\LE^{1}}	\ , 
	\end{align*}
where $\mu^{-1}$ comes from terms where derivatives land on the
cutoff $\chi_{r<\mu t}$. We also used $\mu^{-1}t^{-\frac{\gamma}{4}}\lesssim \epsilon$ on the last line
 which follows since definition \ref{bspara} implies $\epsilon \ll \mu$. Applying estimate \eqref{pbound}
to the first term on the last line above:
		\begin{align*}
		\lp{ \chi_{r< \mu t}  \jap{r}^{-\frac{1}{2}}\tau_{+}^{1-\frac{\gamma}{2}} P\phi}{\ell^{\infty}_{r}L^{2}} 
		\ & \lesssim \ \lp{\chi_{r< \mu t}  \jap{r}^{-\frac{1}{2}}\tau_{+}^{1-\frac{\gamma}{2}}
		\big(\Box_{g}\phi, \ \jap{r}^{-\delta}\partial_{t}\partial_{x}\phi,\jap{r}^{-1-\delta}
		\tau^{-\gamma}\partial_{x}\phi\big)}{\ell^{\infty}_{r}L^{2}}\\
		&  \hspace{1in} +\lp{\chi_{r< \mu t}  \jap{r}^{-\frac{1}{2}}\tau_{+}^{1-\frac{\gamma}{2}}
		\big(\jap{r}^{-1-\delta}\partial_{t}\phi, \partial_{t}^{2}\phi\big)}{\ell^{\infty}_{r}L^{2}}  \\
		&  \lesssim \lp{\Box_{g}\phi}{\ell_{t}^{1}N_1}
		+ \epsilon \big( \lp{\chi_{r< \frac{1}{2}t} \partial_{t}\phi}{\ell_{t}^{\infty}\LE^{1}}
		 +  \lp{\chi_{r< \frac{1}{2}t} \phi}{\ell_{t}^{\infty}\LE^{1}} \big) \ ,
	\end{align*}
	on the last line we used the gain $\tau_{+}^{-\frac{\gamma}{2}}$ to put
	the $\Box_{g}\phi$ term in the form above. For all other terms we split
	$t^{-\frac{\gamma}{2}}$ as we did before and use the definition of the
	norms to finish the proof of \eqref{two_dx_round1}.

For estimate \eqref{two_dx_round2} we can apply the weighted $L^{2}$
estimate \eqref{elliptic} since the exponent of $\jap{r}^{-\frac{1}{2}+\frac{\gamma}{10}}$
is above the threshold value $a=1/2$. This introduces a term supported  where $ \mu t< r $
which we also need to control. After an application of estimate \eqref{pbound} we get:
\begin{align}
	\lp{\chi_{r< \mu t}\jap{r}^{-\frac{1}{2}+\frac{\gamma}{10}} & \tau_{+}^{2-\frac{\gamma}{2}}\big(
	\nabla_{x}^{2}\phi,\jap{r}^{-1}\nabla_{x}\phi\big)}{L^{2}} \notag \\
	&\lesssim \ \lp{\chi_{r< \mu t}\jap{r}^{-\frac{1}{2}+\frac{\gamma}{10}}\tau_{+}^{2-\frac{\gamma}{2}}P\phi}{L^{2}}
	+\lp{\chi_{\mu t<r}\jap{r}^{-\frac{1}{2}+\frac{\gamma}{10}}(\tau_{+}\tau_{0})^{2-\frac{\gamma}{2}}P\phi}{L^{2}} \notag \\
	&\lesssim  \ \lp{\chi_{r< \mu t}\jap{r}^{-\frac{1}{2}+\frac{\gamma}{10}}\tau_{+}^{2-\frac{\gamma}{2}}\big(\Box_{g}\phi,
	\jap{r}^{-\delta}\partial_{t}\partial_{x}\phi,\jap{r}^{-1-\delta}\tau^{ -\gamma}\partial_{x}\phi\big)}{L^{2}}   \notag \\
	&\qquad +\lp{\chi_{r< \mu t}\jap{r}^{-\frac{1}{2}+\frac{\gamma}{10}}\tau_{+}^{2-\frac{\gamma}{2}}\big(\jap{r}^{-1-\delta}\partial_{t}\phi,
	\partial_{t}^{2}\phi\big)}{L^{2}} \!+ \!\lp{\chi_{\mu t<r}\jap{r}^{-\frac{1}{2}
	+\frac{\gamma}{10}}|\tau_{+}\tau_{0}|^{2-\frac{\gamma}{2}}P\phi}{L^{2}}
	  \label{int_ellp1} 
\end{align}
We first bound all the terms above supported inside $r<\mu t $. For the source term
we use the gain $t^{\frac{\gamma}{10}-\frac{\gamma}{4}}$ to get:
\begin{equation*}
	\lp{\chi_{r< \mu t}\jap{r}^{-\frac{1}{2}+\frac{\gamma}{10}}\tau_{+}^{2-\frac{\gamma}{2}}\Box_{g}\phi}{
	L^{2}} \ \lesssim \ \lp{\chi_{r<\frac{1}{2}t} \Box_{g}\phi}{\LE^{2-\frac{\gamma}{4}}} \ \lesssim  \ \lp{\chi_{r<\frac{1}{2}t}\Box_{g}\phi}{N_{1}} \ .
\end{equation*}
For combinations of the form $\nabla\partial_{t}\phi$ we apply estimate \eqref{two_dt_round2}.
For $\tau^{-\gamma}\partial_{x}\phi$ we use $\jap{r}^{-\delta+\gamma}\lesssim 1$
and $t^{-\gamma}\lesssim \epsilon$ to get:
\begin{align*}
	\lp{\chi_{r< \mu t}\jap{r}^{-\frac{3}{2}-\delta+\frac{\gamma}{10}}\tau^{-\gamma}
	\tau_{+}^{2-\frac{\gamma}{2}}\nabla_{x}\phi}{L^{2}} \ \lesssim \ \epsilon
	\lp{\chi_{r< \mu t}\jap{r}^{-\frac{3}{2}+\frac{\gamma}{10}}
	\tau_{+}^{2-\frac{\gamma}{2}}\nabla_{x}\phi}{L^{2}} \ ,
\end{align*}
and bootstrap this term to LHS\eqref{int_ellp1}.
For the $\partial_{t}\phi$ term we multiply estimate \eqref{Stodt} times:
\begin{equation*}
	\jap{r}^{-\frac{3}{2}-\delta+\frac{\gamma}{10}}
	\tau_{+}^{2-\frac{\gamma}{2}}\chi_{r< \mu t} \ .
\end{equation*}	
Squaring the resulting bound, integrating, and using $t^{-\frac{\gamma}{4}}
\lesssim \epsilon$ together with the inclusion
$ \ell_{t,r}^{\infty}L^{2}  \subseteq    t^{-\frac{\gamma}{4}}\jap{r}^{-\delta+\frac{\gamma}{10}}L^{2}  $
then yields:
\begin{equation}
	\lp{\chi_{r< \mu t}\jap{r}^{-\frac{3}{2}-\delta+\frac{\gamma}{10}}\tau_{+}^{2-\frac{\gamma}{2}}
	\partial_{t}\phi}{L^{2}} \
	\lesssim \ \lp{\chi_{r< \mu t}\jap{r}
	^{-\frac{3}{2}-\delta+\frac{\gamma}{10}}\tau_{+}^{1-\frac{\gamma}{2}}\big( S\phi,  r\partial_{x}\phi\big)}{L^{2}}
	\ \lesssim \  \epsilon\lp{\chi_{r< \frac{1}{2}t}\phi}{\LE_{1}^{1}}  \ . \label{dttoS2}
\end{equation}
 This takes care of all the terms supported inside $r<\mu t $ on the last two lines of
 estimate \eqref{int_ellp1}.
 
For the term supported where $\mu t< r $ we go back to Bondi derivatives via the
simple identity: $\partial_{i}  =  \td \partial_{i}+u_i\partial_{u}$. After an application of the gradient bounds
\eqref{grad_u_bnd} and the Klainerman-Sideris estimate \eqref{KS_iden} this yields:
	\begin{align*}
		\lp{\chi_{\mu t<r}\jap{r}^{-\frac{1}{2}+\frac{\gamma}{10}}  (\tau_{+}\tau_{0})^{2-\frac{\gamma}{2}}P\phi}{L^{2}}\
		 & \lesssim  \ \lp{\chi_{\mu t<r}\jap{r}^{-\frac{1}{2}+\frac{\gamma}{10}}(\tau_{+}\tau_{0})^{2-\frac{\gamma}{2}}
		(\nabla_{x}^{2}\phi,\jap{r}^{-\delta}(\tau_{+}\tau_{0})^{-1}\nabla_{x}\phi)}{L^{2}}   \\
		&  \lesssim  \lp{ \chi_{\mu t<r}\jap{r}^{-\frac{1}{2}+\frac{\gamma}{10}}(\tau_{+}\tau_{0})^{2-\frac{\gamma}{2}}
		\big(\td \nabla_{x}^{2}\phi,\jap{r}^{-\delta}\td \nabla_{x}\partial_{u}\phi, \jap{r}^{-2\delta}\partial_{u}^2\phi)\big)}{L^{2}}\\
		&\qquad \qquad + \lp{ \chi_{\mu t<r}\jap{r}^{-\frac{1}{2}-\delta+\frac{\gamma}{10}}(\tau_{+}\tau_{0})^{1-\frac{\gamma}{2}}
		\big(\td \nabla_{x}\phi,  \partial_u\phi\big)}{L^{2}}\\
		 &  \lesssim \epsilon \sup_{t_0\leqslant t\leqslant t_1}\lp{\phi(t)}{\CE_{1}}+ \lp{\Box_{g}\phi}{N_{1}} \ ,
	\end{align*}
	where we have used the gain $\jap{r}^{-\frac{\gamma}{2}}$ on the last line.

\Part{3}{Bounds for $\nabla^{2}\phi$}
Estimate \eqref{two_dx_tless_energy} follows by taking an appropriate
linear combination of inequalities \eqref{two_dt_round1} and \eqref{two_dx_round1}
and using the fact that $\mu$ was chosen to be small enough so we can bootstrap
all terms containing this factor.

Estimate \eqref{twodersclaim10} follows similarly by taking an
appropriate linear combination of inequalities \eqref{two_dt_round2} and \eqref{two_dx_round2}
then absorbing all terms containing a factor of $\mu $.
\end{proof}


\subsection{Proof of the main estimates}

We now prove Theorem \ref{conf_energy_vf_th}.

\begin{proof}
\Part{1}{Proof of estimate \eqref{en_plus_ce}}
We commute equation \eqref{wave_eqn} with the vector field $\partial_{u}$ then use the conformal
energy estimate \eqref{conf_energy_est}
and apply the commutator bound \eqref{comm_est_der_ce} to get:
		\begin{multline}
			\sup_{t_0\leqslant t\leqslant t_1} \lp{\partial_{u}\phi(t)}{\CE}
			+\lp{\partial_{u}\phi}{S} \lesssim \ \lp{ \chi_{r<\mu t}\jap{r}^{-\frac{1}{2}} \tau_{+}
			^{1-\frac{\gamma}{2}}(\nabla^{2}\phi,  \jap{r}^{-1}\nabla\phi)}
			{\ell^{\infty}_{r}L^{2}}  \\
			+\epsilon  \sup_{t_{0}\leqslant t\leqslant t_{1}}\lp{\phi(t)}{\CE_1}
			+ \lp{\partial_{u}\phi(t_0)}{\CE} + \lp{\partial_u(\Box_{g}\phi)}{\ell_{t}^{1}N }
			+\lp{ \jap{r}^{-1-\delta}\tau_{0}^{-1}\tau^{-\gamma}\Box_{g}\phi)}{\ell_{t}^{1}N}  \ . \label{cool_int1}
		\end{multline}
	The last source term above is obviously bounded by $\lp{\Box_{g}\phi}{\ell_{t}^{1}N_1}$.
	Next we bound the terms containing $\chi_{r<\mu t}$: for the terms with only one derivative, we drop
	the weight $\jap{r}^{-1}$ and use $t^{-\frac{\gamma}{4}}\lesssim\epsilon$ together with the
	inclusion $  \ell_{t}^{\infty}\LE  \subseteq     \jap{r}^{-\frac{1}{2}}
	t^{-\frac{\gamma}{4}}\ell_{r}^{\infty}L^2  $. For the terms containing $\chi_{r<\mu t}$
	with two derivatives, we apply \eqref{two_dx_tless_energy} directly.
	Using these two bounds in succession on the RHS \eqref{cool_int1}
	yields:
	\begin{equation*}
			\lp{ \chi_{r< \mu t}\jap{r}^{-\frac{1}{2}} \tau_{+}
			^{1-\frac{\gamma}{2}}(\nabla^{2}\phi,  \jap{r}^{-1}\nabla\phi)}
			{\ell^{\infty}_{r}L^{2}} \
			\lesssim \  \epsilon \lp{\chi_{r<\frac{1}{2}t}(\partial_{t}\phi,\phi)}
		{\ell_{t}^{\infty}\LE^{1}}  \ . 
	\end{equation*}
Bootstrapping the highest order term above to LHS\eqref{cool_int1} finishes the proof
of estimate \eqref{en_plus_ce}.

\Part{2}{Proof of estimate \eqref{conf_loss}}
We commute equation \eqref{wave_eqn} with the vector fields $\{ S,  \Omega_{ij} \}$, use the conformal
energy estimate \eqref{conf_energy_est}, then apply estimate \eqref{comm_est_SR_int}
to the commutators. This yields:
\begin{multline}
	\sup_{t_0\leqslant t\leqslant t_1}\sum_{\Gamma= S,\Omega_{ij}}
	\lp{\Gamma\phi(t)}{\CE}+\lp{\Gamma\phi}{S} \ \lesssim \ 
	 \lp{\jap{r}^{-\frac{1}{2}}\tau_{+}^{2-\frac{\gamma}{2}}
	 \chi_{r< \mu t}\big(\nabla^{2}\phi, \jap{r}^{-1}\nabla\phi\big)}
	{L^{2}}  \label{twoders_sec_round1} \\
	+\epsilon\sup_{t_0\leqslant t\leqslant t_1}\lp{\phi(t)}{\CE_{1}}+
	\sum_{\Gamma= S,\Omega_{ij}}\lp{ \Gamma\phi(t_0)}{\CE}
	+\lp{\Box_{g}\phi}{\ell_{t}^{1}N_{1}} \ .
\end{multline}
We multiply the terms containing $\chi_{r<\mu t}$ above times the weight
$\jap{r}^{\frac{\gamma}{10}}$. For the ensuing terms with only one derivative,
the exponent of $\jap{r}^{a}$ is now above the
$a=3/2$ threshold. This allows us to use the Hardy estimate \eqref{hardy0}
and get:
	\begin{align}
		\lp{\chi_{r< \mu t}\jap{r}^{-\frac{3}{2}+\frac{\gamma}{10}}\tau_{+}^{2-\frac{\gamma}{2}}\nabla\phi}
		 {L^{2}} \ \lesssim \  
		\lp{\chi_{r< \mu t}\jap{r}^{-\frac{1}{2}+\frac{\gamma}{10}}\tau_{+}^{2-\frac{\gamma}{2}}\nabla_{x}\!\nabla\phi}
		 {L^{2}}  \ .
		\label{2dxdtinter} 
	\end{align}
Thus all terms containing $\chi_{r<\mu t}$ in RHS\eqref{twoders_sec_round1} are now
of the form LHS\eqref{twodersclaim10}. An application of the latter estimate
finishes the proof of Theorem \ref{conf_energy_vf_th}. \end{proof}


\section{Global $L^{\infty}$ Decay}

In this section we prove the pointwise bound \eqref{point1}.
Let's begin by showing some preliminary estimates.
\begin{lem}[Preliminary estimates] For test functions $\phi(t,x),f(t,x)$ the following
uniform bounds hold:
\begin{enumerate}[I)]
	\item (Global $L^{\infty} estimate$)
		\begin{equation}
			\lp{\jap{r}^\frac{3}{2}\tau_0^\frac{1}{2}\phi }{L^\infty_x} \ \lesssim \ 
			\sum_{k + |J|\leqslant 2}\lp{(\jap{r}\tau_0\partial_u)^k(\jap{r}\td \partial_x)^J \phi}{L^2_x}
			\ . \label{global_L00}
		\end{equation}
	\item (Interior $L^{\infty}$ estimate) Assume that $\phi$ is supported in $r<\frac{1}{2}t$. Then one has:
		\begin{equation}
			\lp{\phi }{L^\infty_x} \ \lesssim \ \lp{ \jap{r}^{\frac{1}{2}}(\nabla^{2}\phi,
			  \jap{r}^{-1}\nabla\phi,   \jap{r}^{-2}\phi)}{\ell^\infty_r L_{x}^{2}}\ . \label{int_L00}
		\end{equation}
\item (Average to uniform bounds via scalings) Let $f$ be supported
on the time interval $[T^{*},\infty)$. Then for all $a\in\mathbb{R}$ and $T>T^{*}$:
\begin{equation}
		\sup_{0\leqslant t\leqslant T}\!
		\lp{\tau_+^a f(t)}{L^2_x}
		\ \lesssim \ 
		\lp{ \tau_+^{a-\frac{1}{2}} (f,r\td \partial_r f, Sf)}{L^{2}[0,T]}
		\ . \label{ave_scaling_est}
\end{equation}
%
\end{enumerate}
\end{lem}

\begin{proof}
\step{1}{Proof of Estimate \eqref{global_L00}} Inside the set $r\leqslant 1$ this follows from
the standard $L^{\infty}-L^{2}$ Sobolev estimate. For the complement it suffices to
consider the region $\frac{1}{2}t<r<\frac{3}{2}t$ as the remainder is 
easier to handle because we have $u=t-r$ there. Using dyadic cutoffs we may
assume $\phi$ is supported where $\tau_-\approx 2^k$ and $r\approx 2^j$. Using
angular sector cutoffs in the $x$ variable we may further assume without loss of generality 
that $\phi$ is supported in a $\frac{\pi}{4}$ wedge about the $x^1$ axis.  
Now introduce new variables on $t=const$:
\begin{equation}
		y^1\ =\ 2^{-k}u \ , \qquad
		y^2\ =\ 2^{-j}x^2 \ , \qquad
		y^3\ =\ 2^{-j}x^3 \ . \notag
\end{equation}
There exists vector fields $e_\alpha$, $\alpha=0,1,2,3$ such that 
$\partial_{y^i}|_{t=const}=\sum_\alpha c^i_\alpha e_\alpha$ where $c^i_\alpha$ are uniformly
bounded and such that:
\begin{equation}
		2^{\frac{1}{2}k+j}\sum_{|I|\leqslant 2}\lp{e^I \phi}{L^2(dy)} \ \lesssim \ 
		\sum_{k + |J|\leqslant 2}\lp{(r\tau_0\partial_u)^k(r\td \partial_x)^J \phi}{L^2(dx)}
		\ . \notag
\end{equation}
Estimate \eqref{global_L00} follows from this last line by concatenating the Sobolev embeddings 
$H^1\subseteq L^6$ and $W^{1,6}\subseteq L^\infty$.

\step{2}{Proof of \eqref{int_L00}}
Once again by the $L^{\infty}-L^{2}$ Sobolev estimate it suffices to prove the result outside $r\leqslant 1$.
By the support property we have $\tau_0\approx 1$, therefore applying estimate \eqref{global_L00}
to $r^{-\frac{3}{2}}\chi_{ \jap{r} \approx 2^{k}}\phi$, taking $\sup_{k}$ and using the fact that $\partial_{u},\td \partial_x$
are bounded linear combinations of $\partial_t, \partial_x$ derivatives yields the claim.

\step{3}{Proof of \eqref{ave_scaling_est}}
Integrating the time derivative of $(t+r)^{2a}f^2$ over $0<t<T$ we have:
\begin{equation}
		 (T+r)^{2a}f^2(T) \ = \ r^{2a}f^2(0) + 2a \int_0^T (t+r)^{2a-1}f^2dt +
		 2\int_0^T (t+r)^{2a}f\partial_tf dt \ . \notag
\end{equation}
 Using Cauchy-Schwartz on the last RHS term, integrating in $x$ over $0<r<\infty$,
 then using the support property of $f$ we get:
\begin{equation}
		\lp{\tau_+^a f(T)}{L_x^2} \ \lesssim \ 
		\lp{\tau_+^{a-\frac{1}{2}}(f, r\td \partial_rf, Sf)}{L^2[0,T]} \  . \notag
\end{equation}
This yields \eqref{ave_scaling_est}.

\end{proof}

Now we demonstrate the main $L^{\infty}$ bound.

\begin{proof}[Proof of estimate \eqref{point1}]
By Lemma \ref{reduction} and by using an appropriate cutoff function we may
assume $\phi$ is supported on $[T^{*},\infty)$ and that $T>T^{*}$. We then estimate the
timelike and null/spacelike regions separately.
Let $\mu$ be as in definition \ref{bspara}.a.

\step{1}{Estimate for $\mu t < r$} Applying \eqref{global_L00}   to $\chi_{\mu t < r}\phi $ followed
by  \eqref{KS_iden}:
\begin{align*}
	\lp{\chi_{\mu t < r}r^\frac{3}{2}\tau_0^\frac{1}{2}\phi (t)}{L_x^\infty}
	\ &\lesssim  \sum_{k + |J|\leqslant 2}\lp{(r\tau_0\partial_u)^k(r\td 
	\partial_x)^J \phi (t) }{L^2_x}\ \lesssim \  \lp{\phi(t)}{\CE_{1}}  +  \lp{r^2 \tau_{0}\Box_{g}\phi (t)}{ \ell_r^1 L^2_x} \ .
\end{align*}
We take $\sup_t$ then apply the bound \eqref{conf_energy_est_vf} to control
$\sup_{0\leqslant t\leqslant T}\lp{\cdot}{\CE_{1}}$. For the
source term we use the trace estimate \eqref{ave_scaling_est}.

\step{2}{Estimate for $r<\mu t$} Let $\chi_{r<\mu t}$ be a smooth
cutoff to the region $r< \mu t$. Applying estimate
\eqref{int_L00} to $\chi_{r< \mu t}\tau_{+}^{\frac{3}{2}}  \phi$ we get:
\begin{equation}
	\sup_{0\leqslant t\leqslant T}\lp{ \chi_{r< \mu t} \tau_{+}^{\frac{3}{2}}\phi (t) }{L^\infty_x} 
	\ \lesssim \  \sup_{0\leqslant t\leqslant T}\lp{\td \chi_{r<\mu t} \tau_{+}^{\frac{3}{2}}
	\jap{r}^{\frac{1}{2}}(\nabla^{2} \phi,   \jap{r}^{-1}\nabla \phi,  \jap{r}^{-2} \phi)(t)}{\ell^\infty_r L_{x}^{2}}
	\label{int_prelim_pt1}  \ ,
\end{equation}
where $\td \chi_{r<\mu t}$ is a smooth cutoff with slightly larger support. For the undifferentiated
terms we apply \eqref{ave_scaling_est} with $a=3/2$ to $\jap{r}^{-\frac{3}{2}}\td \chi_{r<\mu t}\chi_{r\approx 2^{k}}\phi$
to produce:
\begin{align}
	\sup_{0\leqslant t\leqslant T}\lp{\jap{r}^{-\frac{3}{2}}\tau_{+}
	^{\frac{3}{2}}\td \chi_{r<\mu t}\phi}{\ell^\infty_r L_{x}^{2}}
	\ \lesssim \ \lp{\chi_{r<\frac{1}{2}t}\phi}{\ell_{t}^{\infty}\LE_{1}^{1}} \ .
	\label{int_prelim_pt2}
\end{align}
In remains for us to control the rest of the terms on RHS\eqref{int_prelim_pt1}.
For terms with one derivative we apply the Hardy estimate \eqref{hardy0} with $a=1/2$ .
Thus it suffices to estimate the terms with two derivatives. For this we claim:
\begin{align}
	\lp{\td \chi_{r<\mu t}
	\tau_{+}^{\frac{3}{2}}\jap{r}^{\frac{1}{2}}\nabla^{2}\phi(t)}{L^2_x}
	\ \lesssim \ \lp{\phi(t)}{\CE_{1}}  \ + \ \lp{\tau_{+}^{\frac{3}{2}}\jap{r}^{\frac{1}{2}}
	\Box_{g}\phi (t)}{L^2_x}  \ . \label{P_l2}
\end{align}
By using a similar proof as inequality \eqref{twodersclaim10}, this estimate can be further reduced to
proving:
\begin{align}
	\lp{\td \chi_{r< \mu t} \tau_{+}^{\frac{3}{2}}\jap{r}^{\frac{1}{2}}\nabla\partial_{t}\phi(t)}{L^2_x}
	\ &\lesssim \   \lp{\phi(t)}{\CE_{1}}  +  \mu \lp{\td \chi_{r<\mu t} \tau_{+}^{\frac{3}{2}}\jap{r}^{\frac{1}{2}}
	\nabla_{x}^{2}\phi(t)}{L^2_x} \ , \label{dtdx_l2}\\
	\lp{\td \chi_{r<\mu t} \tau_{+}^{\frac{3}{2}}\jap{r}^{\frac{1}{2}}\nabla_{x}^{2}\phi(t)}{L^2_x}
	\  &\lesssim  \ \lp{\phi(t)}{\CE_{1}}  +   \lp{\tau_{+}^{\frac{3}{2}}\jap{r}^{\frac{1}{2}}\Box_{g}\phi(t)}{L^2_x}
	+  \lp{\td \chi_{r<\mu t} \tau_{+}^{\frac{3}{2}}\jap{r}^{\frac{1}{2}}\nabla\partial_{t}\phi(t)}
	{L^2_x} \label{dxdx_l2} \ . 
\end{align}
Estimate \eqref{dtdx_l2}, in turn, follows by a slightly simpler version of the proof of \eqref{two_dt_round2} and
using the support property. Likewise, estimate \eqref{dxdx_l2} follows by applying the
elliptic estimate \eqref{elliptic} with $a=1/2$ then following a similar argument to estimate \eqref{two_dx_round2}
and using the support property. Taking an appropriate linear combination of estimates \eqref{dtdx_l2} and \eqref{dxdx_l2}
then using the smallness of $\mu$ to bootstrap errors proves \eqref{P_l2}.
To get estimate \eqref{point1} from this, we take $\sup$ in $t$ and use the bound \eqref{conf_energy_est_vf}
to control $\sup_{0\leqslant t\leqslant T}\lp{\cdot}{\CE_{1}}+\lp{\chi_{r<\frac{1}{2}t}\phi}{\ell_{t}^{\infty}\LE_{1}^{1}}$
along with the trace estimate \eqref{ave_scaling_est} for the source term. This, together with
Lemma \ref{reduction} finishes the proof.
\end{proof}


\section{Appendix: Weighted $L^2$-Elliptic Estimates}
Our main result in this section is the following:

\begin{thm}[Global elliptic estimate]
Let the operator $P(t,x,\nabla_{x})=\partial_i h^{ij}\partial_j$ with $h^{ij}=d^\frac{1}{2}g^{ij}$
be uniformly elliptic. Suppose that for all $t > 0$, the $h^{ij}$ satisfy the uniform bounds:
\begin{equation}
		|(\jap{r}\partial_x)^\alpha (h^{ij}-\delta^{ij})| \ \lesssim \ \jap{r}^{-\delta} \ . \label{subtract}
\end{equation}
Then for all $t> 0$ the operator $P$ satisfies the fixed-time estimates:
\begin{equation}
		\lp{\jap{r}^{a}\big(\nabla^2_{x}\phi ,  \jap{r}^{-1}
		\nabla_{x}\phi\big)(t)}{ L_{x}^2} 
		\ \lesssim \  \lp{\jap{r}^a P \phi(t)}{ L_{x}^2} \ , \qquad 
		-\frac{1}{2} <  a <  \frac{3}{2} \ , \label{elliptic}
\end{equation}
where the implicit constants are independent of $t$.
\end{thm}

\begin{proof}
Let $\Delta, \Delta^{-1}$ denote the standard 3D Laplacian and its inverse, respectively.
Write $P\phi=F$. We approximately solve for $F$ in terms of a Neumann series:
\begin{equation}
		\td{\phi} \ = \ \Delta^{-1} \sum_{i=0}^k  R^i F \ ,
		\qquad R \ = \ I-P \Delta^{-1} \ .
		 \notag
\end{equation}
Then we have: $P (\td{\phi}-{\phi}) \ = \ R^{k+1}F$.
Therefore, setting $L_{x}^{2,a}$ for the norms on line \eqref{elliptic}
it suffices to show:
\begin{equation}
		\nabla^2\Delta^{-1}
		: L_{x}^{2,a}  \to  L_{x}^{2,a} \ , \qquad
		\jap{r}^{-1}\nabla\Delta^{-1}: L_{x}^{2,a}  \to 
		L_{x}^{2,a} \ , \qquad
		R: L_{x}^{2,a}  \to  L_{x}^{2,a+\delta} \ ,  
		\label{kernel_bounds}
\end{equation}
for the range $\frac{1}{2} < a$ and  $a+\delta < \frac{3}{2} \ ,$
followed by the non-perturbative estimate:
\begin{equation}
		\lp{\jap{r}\nabla^2_{x}\phi(t)}{ L_{x}^2} \ +  \
		\lp{\nabla_{x}\phi(t)}{ L_{x}^2} 
		\ \lesssim \ \lp{\jap{r} P \phi(t)}{L_{x}^2} \ . \label{exact_elliptic}
\end{equation}
\end{proof}

\begin{proof}[Proof of \eqref{kernel_bounds}]
We decompose into dyadic scales $|x|\sim 2^{i}$, $|y|\sim 2^{j}$ with $|x|\lesssim 1$ when $i=0$
since the weights are non-singular.

\step{1}{$\nabla^2\Delta^{-1}: L_{x}^{2,a} \ \to \ L_{x}^{2,a}$ is bounded} To
establish this it suffices to show:
\begin{align}
	\sum_{i,j} \big|\dint_{|x|\sim 2^{i} , |y|\sim 2^{j}}
	 F(x)K_{1}(x-y)G(y)\ dxdy\big| \ \lesssim  \lp{F}{
	 L_{x}^{2,-a}}\lp{G}{L_{x}^{2,a}} \ , \label{est_ap_1}
\end{align}
where $K_{1}$ is the convolution kernel for $\nabla^2\Delta^{-1}$.
We break up the proof into cases:

\case{1} {$|i-j|=O(1)$} The operator defined above is a singular integral operator.
In this case the weights $2^{-a j}\approx 1$ and $2^{-a i}\approx 1$ balance
since they are both approximately of size one. By Cauchy-Schwarz: 
\begin{align*}
	\sum_{i+j=O(1)} \big|\dint_{|x|\sim 2^{i} , |y|\sim 2^{j}}
	 F(x)K_{1}(x-y)G(y)\ dxdy\big| \  & \lesssim \ \sum_{i,j} \lp{\chi_{i}F}{L^{2,-a}}\lp{\chi_{j}G}{L^{2,a}} \ , 
\end{align*}
with $\chi_{i},\chi_{j}$ smooth cutoff functions supported where $|x|\sim 2^{i}$, $|y|\sim 2^{j}$, respectively.

\case{2} {$i>j+c$} We now have $|K_{1}(x-y)|=O(|x|^{-3})$ and since convolution with an $L^{1}$ function  is a bounded 
operator in any $L^{p}$ space with $p\geq 1$:
\begin{equation*}
	\sum_{i>j+c} \big|\dint_{|x|\sim 2^{i} , |y|\sim 2^{j}}
	 F(x)K_{1}(x-y)G(y)\ dxdy\big| \
	 \lesssim \ \sum_{i>j+c}  2^{-\frac{3}{2}(i-j)}\cdot 2^{a(i-j)} \lp{\chi_{i}F}{L^{2,-a}} 
	 \lp{\chi_{j}G}{L^{2,a}} \ ,
\end{equation*}
which forces $-\frac{3}{2}+a\leqslant 0$.

\case{3} {$j>i+c$} by switching the roles of $x$ and $y$ and using the same argument as case 2:
\begin{equation*}
	\sum_{j>i+c} \big|\dint_{|x|\sim 2^{i} , |y|\sim 2^{j}}
	 F(x)K_{1}(x-y)G(y)\ dxdy\big| \ 
	 \lesssim \ \sum_{j>i+c}  2^{\frac{3}{2}(i-j)}\cdot 2^{a(i-j)} \lp{\chi_{i}F}{L^{2,-a}} 
	 \lp{\chi_{j}G}{L^{2,a}} \ ,
\end{equation*}
which is convergent for $-\frac{3}{2}-a < 0$. This proves \eqref{est_ap_1}.

\step{2}{$\jap{x}^{-1}\nabla\Delta^{-1}: L_{x}^{2,a} \ \to \ L_{x}^{2,a}$ is bounded} We again aim to show:
\begin{equation}
	\sum_{i,j} \big|\dint_{|x|\sim 2^{i} , |y|\sim 2^{j}}
	 F(x)K_{2}(x-y)G(y)\ dxdy\big|  \lesssim  \lp{F}{L_{x}^{2,-a}}\lp{G}
	 {L_{x}^{2,a}} \ , \label{est_ap_2}
\end{equation}
with the kernel $K_2(x,y)$ the convolution kernel for $\jap{r}^{-1}\nabla\Delta^{-1}$.

\case{1}{$|i-j|=O(1)$} Here we have $K_{2}(x,y)=O(\jap{x}^{-1}|x-y|^{-2})$. By the Hardy-Littlewood-Sobolev inequality:
\begin{multline*}
	\sum_{i+j=O(1)} \big|\dint_{|x|\sim 2^{i} , |y|\sim 2^{j}}
	 F(x)K_{2}(x-y)G(y)\ dxdy\big| \  \lesssim \ \sum_{i+j=O(1)} \big|\dint_{|x|\sim 2^{i} , |y|\sim 2^{j}}
	 \frac{F(x)G(y)}{\jap{x}|x-y|^{2}}\ dxdy\big| \\
	  \lesssim \ \sum_{i,j} 2^{-i}\cdot2^{\frac{1}{2}i}
	 \cdot 2^{\frac{1}{2}j}\lp{\chi_{i}F}{L^{2}}\lp{\chi_{j}G}{L^{2}} \ \lesssim \ \sum_{i,j} \lp{\chi_{i}F}{L^{2,-a}}\lp{\chi_{j}G}{L^{2,a}} \ , 
\end{multline*}
where we've used $2^{-i}\approx2^{\frac{1}{2}i}\cdot 2^{\frac{1}{2}j}\approx 1$ on the last line.

\case{2} {$i>j+c$} Here we have $|K_{2}(x-y)|=O( \jap{x}^{-1}|x|^{-2})$. Therefore:
\begin{equation*}
	\sum_{i>j+c} \big|\dint_{|x|\sim 2^{i} , |y|\sim 2^{j}}
	 F(x)K_{2}(x-y)G(y)\ dxdy\big| \ 
	 \lesssim \ \sum_{i>j+c}  2^{-\frac{1}{2}i}\cdot2^{-(i-j)}\cdot 2^{a(i-j)} \lp{\chi_{i}F}{L^{2,-a}} 
	 \lp{\chi_{j}G}{L^{2,a}} \ ,
\end{equation*}
and since $i>j+c$ the extra $2^{-\frac{1}{2}i}$ helps us get $-\frac{3}{2}+a< 0$ once again.

\case{3}{$j>i+c$} Here we have $|K_{2}(x-y)|=O( \jap{x}^{-1}|y|^{-2})$. Hence:
\begin{equation*}
	\sum_{j>i+c} \big|\dint_{|x|\sim 2^{i} , |y|\sim 2^{j}}
	 F(x)K_{2}(x-y)G(y)\ dxdy\big| \ 
	 \lesssim \ \sum_{j>i+c}  2^{\frac{1}{2}i}\cdot2^{(i-j)}\cdot 2^{a(i-j)} \lp{\chi_{i}F}{L^{2,-a}} 
	 \lp{\chi_{j}G}{L^{2,a}} \ ,
\end{equation*}
and since $j>i+c$ the extra $2^{\frac{1}{2}i}$ gives us the restriction $-\frac{1}{2}-a< 0$.

\step{3}{$R: L^{2,a} \ \to \ L^{2,a+\delta}$ is bounded}
By estimate \eqref{subtract} we have:
\begin{equation*}
	P-\Delta \ = \ \partial_{i}(h^{ij}-\delta^{ij})\partial_{j}
	\ = \ O(\jap{r}^{-\delta})\partial_x^2+O(\jap{r}^{-1-\delta})\partial_x \ .
\end{equation*}
This observation together with the results above finish the proof.
\end{proof}

\begin{proof}[Proof of \eqref{exact_elliptic}]
Let $D$ denote the Levi-Civita connection for $h$ and let $dV_{h}$
be the corresponding volume form. We have the estimate: 
\begin{align}
	\int_{\mathbb{R}^{3}} D^{i}\phi D_{i}\phi \ dV_{h} \ \lesssim \ \int_{\mathbb{R}^{3}} \jap{r}^{2}|P\phi|^{2} \ dV_{h} \ ,
	\label{claim_app}
\end{align}
which follows from Green's identity:
\begin{equation*}
	- \int_{\mathbb{R}^{3}} D^{i}\phi D_{i}\phi \ dV_{h}
	\ = \ \int_{\mathbb{R}^{3}} P\phi\cdot \phi \ dV_{h} \ , 
\end{equation*}	
by taking absolute value, applying Young's inequality and
using the Hardy estimate:
\begin{equation*}
	\int_{\mathbb{R}^{3}} \big|r^{-1}\phi\big|^{2} \ dV_{h} 
	\ \lesssim \ \int_{\mathbb{R}^{3}} |D\phi|^{2} \ dV_{h} \ .
\end{equation*}	
To prove the estimate for two derivatives we integrate by parts twice 
then take absolute value, apply \eqref{subtract}  
together with estimate \eqref{claim_app} and Young's inequality to produce:
\begin{align*}
	\int_{\mathbb{R}^{3}} \big|\jap{r}^{2}(D_{i}D_{j}\phi) (D^{i}D^{j}\phi)\big| \ dV_{h}\ 
	 &  \lesssim \ \int_{\mathbb{R}^{3}} \big( \jap{r}^{2}|P\phi|^{2} +
	 \jap{r}|P\phi| |D\phi| +\jap{r}^{-1-\delta}|\phi| |D\phi| + D_{i}
	\phi D^{i}\phi \big) \ dV_{h}\\
	& \lesssim \   \int_{\mathbb{R}^{3}} \jap{r}^{2}|P\phi|^{2} dV_{h} \ .
\end{align*}
This finishes the proof of \eqref{exact_elliptic}.
\end{proof}


\section{Acknowledgements}

The author would like to thank his advisor Jacob Sterbenz for suggesting
the problem, for sharing many valuable insights, and for providing unconditional
support throughout the writing of this work.


\bibliographystyle{abbrv}  

\end{document}